\documentclass[11pt]{article}
\usepackage{amssymb}
\usepackage{amsmath}
\usepackage{graphicx}
\usepackage{soul}
\usepackage{subcaption}
\usepackage{float}
\usepackage{fancyvrb}
\usepackage{xcolor}
\usepackage{subcaption}
\usepackage{hyperref}
\usepackage{chngcntr}
\usepackage{makeidx}

\newtheorem{theorem}{Theorem}
\newtheorem{definition}{Definition}
\newtheorem{lemma}[theorem]{Lemma}

\usepackage{pst-plot}
\usepackage{tikz}
\usetikzlibrary{arrows}
\usetikzlibrary{er,positioning}
\usetikzlibrary{patterns}
\usepackage[shortlabels]{enumitem}

\newcommand{\inte }{{\rm int}\,}
\newcommand{\diam }{\,{\rm diam}\,}
\newcommand{\dom }{\,{\rm dom}\,}

\def\reals{\mathbb{R}}

\newcommand{\comment}[1]{\mbox{}}

\newcommand{\config}[3]{
        \begin{minipage}[c]{.4\textwidth}
        \centering
        \includegraphics[width=3.3cm]{#1}\\
        (\ref{fig:#1})
        #2
        \end{minipage}
        \begin{figure}\phantomcaption\label{fig:#1}\end{figure}
    \begin{minipage}[c]{.6\textwidth}
    \tiny #3
\end{minipage}
}

\def\qed{{\hfill{\vrule height5pt width3pt depth0pt}\medskip}}

\title{Central configurations in planar  $n$-body problem for $n=5,6,7$ with equal masses}

\author{Ma{\l}gorzata Moczurad and Piotr Zgliczy\'nski\footnote{Partially supported by the NCN grant 2015/19/B/ST1/01454}\\
   \{malgorzata.moczurad, piotr.zgliczynski\}@ii.uj.edu.pl \\
Faculty of Mathematics and Computer Science, Jagiellonian University,\\
ul. prof. Stanis\l awa \L ojasiewicza 6,
30-348 Krak\'ow, Poland
}

\tikzstyle{io} = [fill=black,inner sep=2pt,circle]

\begin{document}

\maketitle

\begin{abstract}
We give a computer assisted proof of the full listing of central configuration for $n$-body problem for Newtonian potential on the plane for $n=5,6,7$ with equal masses. We show all these central configurations have a reflective symmetry with respect to some line. For $n=8,9,10$ we establish the existence of central configurations without any reflectional symmetry.
\end{abstract}

\textbf{keywords:} central configurations, symmetries, interval arithmetic, Kraw\-czyk operator

\tableofcontents

\section{Introduction}
A central configuration\index{central configuration}, denoted as CC, is an initial configuration $(q_1, \ldots, q_n)$ in the Newtonian $n$-body problem, such that if the particles were all released with zero velocity, they would collapse toward the center of mass $c$ at the same time.
In the planar case CCs are initial positions for periodic solutions which preserve the shape of the configuration. CCs play also an important role in the study of the topology of integral manifolds in the $n$-body problem (see~\cite{MSch} and the references given there).


In this paper, for $n=4,5,6,7$, we consider two questions:
\begin{itemize}[--]
\item finding all CCs in the $n$-body problem on the plane ($d = 2$) with equal masses and
\item showing that each CC has a line of reflection symmetry.
\end{itemize}
For $n=8,9,10$ we establish the existence of some non-symmetric CCs previously found numerically in \cite{MNum,F02,Si}.

\subsection{The state of the art}
The listing (apparently full for $n \leqslant 9$ ) of central configurations with equal masses was given by Ferrario in unpublished notes~\cite{F02} for $n \in \{3,\dots,10\}$ and by Moeckel \cite{MNum} for $n \leqslant 8$.
For $n=4$ it was shown by Albouy that all CCs have some reflectional symmetry \cite{A95} and later in \cite{A96} with computer assistance
 the full list of central configurations was given.
 From numerical simulations (see for example \cite{MNum, F02}) it is apparent that all CCs with equal masses have some reflectional symmetry for $n=5,6,7$. Moeckel \cite{MNum} have found numerically some CCs without any symmetry for $n=8$. Also for $n=9$ Sim\'o \cite{Si} has found
2 families, non-equivalent, and without any symmetry. Some CCs without symmetry for $n=10$ can be found also in \cite{F02}.

The investigations  of central configurations for equal masses is a subcase of more general problem of central configurations with arbitrary positive masses.
The general conjecture of finiteness of central configurations (relative equilibria) in the $n$-body problem is stated in \cite{Wintner} and appears as the sixth
problem of Smale's eighteen problems for the 21st century \cite{SmNext}.
There are many works on the existence of some particular  central configurations. Here we discuss only those papers which aim to more general statement about all CCs.
The two most important works are \cite{HM} and \cite{AK}. In \cite{HM} the finiteness of CC for $n=4$ for any system of positive masses was proved with computer assistance. In \cite{AK} for $n=4$ problem the finiteness of CC was proven without computer assistance. In the same paper, the finiteness for $n=5$ was proven for arbitrary positive masses, except  perhaps if the 5-tuple of positive masses belongs to a
given codimension 2 subvariety of the mass space. It is interesting to notice that the equal masses case treated in our paper belongs to this subvariety. For the spatial  $5$-body problem Moeckel in~\cite{M01}  established the generic finiteness of  Dziobek's CCs (CCs which are not planar). A computer assisted work by Hampton and Jensen~\cite{HJ}  strengthens this result by giving an explicit list of conditions for exceptional values of masses.  A common feature of these works
is that they give a quite poor estimate for the maximum number of central configurations. In this context it is worth to mention the work of Simo,  based on extensive numerical studies. In~\cite{Si78} he gives  the number of CCs for all possible masses for $n=4$.

In~\cite{LS} the spatial 5-body problem with equal masses was considered. A complete classification of the isolated central
configurations of the 5-body problem was given (this includes also a planar isolated CCs). The approach has a numerical component, hence it cannot be claimed to be fully rigorous. Also
the proof does not exclude the possibility that a higher dimensional
set of solutions exists. On the other hand the existence of identified isolated CC, has been proven using the Krawczyk's operator, i.e.\ a tool from interval arithmetic we also use.
Kotsireas (see \cite{Kotsireas} and references given there)  considers the 5-body problem with equal masses. He gives computer assisted proof of a full list of all such configurations and shows that each of them posses some reflectional symmetry.

The above  mentioned works study the polynomial equations derived from the equations for CC  using the (real or complex) algebraic geometry tools.
In contrast, we take a different approach: we use standard interval arithmetic tools, hence in principle we can treat also other potentials which cannot be reduced to polynomial equations.

\subsection{The main results}
\begin{theorem}
\label{thm:main}
There exist only a finite  number of various types of CCs, for $n=5,6,7$ the planar $n$-body Newtonian problem with equal masses. They are listed in Section~\ref{sec:listing}.
Any CC can be obtained from one of them by suitable composition of translation, rescaling, rotation, reflection
and permutation of bodies.
Moreover, each of these central configurations has some reflectional symmetry.
\end{theorem}

\begin{theorem}
\label{thm:non-sym}
For $n=8,9,10$ in the planar $n$-body Newtonian problem with equal masses, there exist CCs without any line of reflectional symmetry.  They are listed in Section~\ref{sec:asymm-cc}.
\end{theorem}

In the case of equal masses
   one can consider equivalences in two different ways: either one
   passes from a solution to another one by rotation (scaling is already
   taken into account) or one can also add permutations and reflections. For instance,
   for 4 bodies in the first criterion of equivalence there are 50
   classes (see numerical work by \cite{Si78}), while in the second only 4 classes.   In this paper we use this second criterion for equivalence.

Let us briefly describe our method. This is basically a brute force approach using standard interval arithmetic tools. Throughout the paper we will use often \emph{box} or $\emph{cube}$ to describe a set which is a product of intervals (some of them can be degenerate). The interval arithmetic allows to evaluate elementary functions on the box in a single call, i.e. the box is returned containing the true result for all points in the argument box (see for example \cite{Mo,N}).  When looking for CCs we explore the whole configuration space (modulo some a priori bounds), and it is surprising
that the most demanding part is to exclude the possibility of the existence of CC in a given box. Once we are `very close' to an isolated CC, it is relatively easy to establish its existence and local uniqueness using the Krawczyk's operator~\cite{K}.  The additional difficulty is that the potential
contains singularity, which introduces some non-compactness in the domain to be covered.  Our algorithm, which is more or less a binary
search algorithm, scales poorly with $n$ --- this is the \emph{dimensionality curse} (see~\cite{TWW}), which means that the complexity of our algorithm grows exponentially with $n$. For example, assume that we can decide if a box in the configuration space contains  some CC, only when its diameter is  less than $10^{-2}$ in each direction. Then adding a new body in  $[-1,1]\times [-1,1]$ multiplies the number of boxes to be examined by $(2/10^{-2})^2 = 4\cdot 10^4$. For this reason we were not able to obtain a rigorous listing of CCs for $n=8$. Note that for $n=5$ the computations were done in 24 seconds, for $n=6$ it took
about one hour to get the result,  while for $n=7$ we needed almost a hundred hours (see Sec.~\ref{subsec:tech-data} for more technical data).

For any CC from the listing in~\cite{MNum} or \cite{F02} for $n\geqslant 8$ we have found no difficulty proving its existence and local uniqueness.
In particular we confirmed the existence of non-symmetric  planar CCs for $n=8,9,10$ (see Theorem~\ref{thm:non-sym}).

The paper is organized as follows. In Section~\ref{sec:cc-eq} we recall the equations for the central configurations and their basic properties. In Section~\ref{sec:ap-bounds} we derive several a priori bounds for CC, so that we obtain a compact domain for our search algorithm.
In Section~\ref{sec:excl-tests} we discuss various tests which are used to show that a given box does not contain any CC. In Section~\ref{sec:red-sys-equiv} we derive a reduced set of equations for CC.  This is necessary because  to apply the Krawczyk's method we need to ensure that the system of equations does not contain any degeneracies, which are due to  symmetries of the original system of equations for CCs.
In Section~\ref{sec:cap} we give assumptions and basic ideas concerning the computer assisted proofs of main Theorems~\ref{thm:main} and~\ref{thm:non-sym} and we explain the Krawczyk's method. Details of the algorithm are described in Section\ref{sec:alg}.
In Section~\ref{sec:dependency} we present an attempt to minimize the dependency problem in interval arithmetic in the evaluation of the gravitational force.
In the Appendix we give an output of the program establishing Theorems~\ref{thm:main} and~\ref{thm:non-sym} (also for $n= 3,4$) and pictures of all CCs found.

\section{Equations for central configurations}
\label{sec:cc-eq}

Let $q_i \in \mathbb{R}^d$, $i=1,\dots,n$ and $d\geqslant 1$\index{$q_i$} (the physically interesting cases are $d=1,2,3$), where $q_i$ is a position of $i$-th body
with mass $m_i \in \mathbb{R}_+$\index{$m_i$}.
Let us set
\begin{equation}
  M=\sum_{i=1}^n m_i.\index{M} \label{eq:def-mass-sum}
\end{equation}
\emph{Central configurations}\index{central configurations} are solutions of the following system of equations (see~\cite{MSch}):
\begin{equation}
  \lambda (q_i-c) = \sum_{\substack{j=1\\
j\neq i}}^n \frac{m_j}{r_{ij}^3}(q_i - q_j)=:\frac{1}{m_i}f_i(q_1,\dots,q_n), \quad i=1,\dots,n \label{eq:cc-with-lambda}
\end{equation}
where $\lambda\in\reals$ is a constant, $c=\left(\sum_{i=1}^n m_i q_i\right)/M$ is center of mass,  $r_{ij}=r_{ji}=|q_i - q_j|$\index{$r_{ij}$} is the Euclidean distance between $i$-th and $j$-th bodies and
$(-f_i)$\index{$f_i$} is the force which acts on $i$-th body resulting from the gravitational pull of other bodies.
The system of equations (\ref{eq:cc-with-lambda}) has the same symmetries as the $n$-body problem. It is invariant with respect to group of
isometries of $\mathbb{R}^d$ and the scaling of variables.

In the planar case if we consider the bodies in a rotating system (with the center of mass at the origin)
   with constant angular velocity $\omega=\sqrt{\lambda}$, the physical meaning of (\ref{eq:cc-with-lambda}) is
   obvious: the gravitational attraction is compensated by the
   centrifugal force and the central configurations are fixed points in the rotating frame (see~ \cite{MSch} and the references given there).

The system~(\ref{eq:cc-with-lambda}) has $dn$ equations and $dn + 1$ unknowns: $q_i \in \mathbb{R}^d$ for $i = 1, \ldots, n$ and $\lambda \in \mathbb{R}_+$.
The system has a $O(d)$ and scaling symmetry (with respect to $q_i$'s and $m_i$'s).
If we demand that $c=0$ (which is obtained by a suitable translation) and $\lambda=1$ (which can be obtained by rescaling $q_i$'s or $m_i$'s) we obtain the equations (compare \cite{MSch,Mlect2014,AK})
\begin{equation}
  q_i= \sum_{j,j\neq i} \frac{m_j}{r_{ij}^3}(q_i - q_j)=:\frac{1}{m_i}f_i(q_1,\dots,q_n), \quad i=1,\dots,n.
  \index{$q_i$} \label{eq:cc-kart}
\end{equation}
It is easy to see that if (\ref{eq:cc-kart}) is satisfied, then $c=0$ (see Sec.~\ref{subsec:eq-conservation}) and  (\ref{eq:cc-with-lambda}) also holds for $\lambda=1$.
A $q = (q_1,\dots,q_n) \in \left(\mathbb{R}^d\right)^n$ is called a \emph{configuration}\index{configuration}.  If $q$ satisfies (\ref{eq:cc-kart}) then it is called
a \emph{normalized central configuration}\index{normalized central configuration} (abbreviated as CC)\index{CC}.
For the future use we introduce the function $F:\Pi_{i=1}^n\mathbb{R}^{d } \to \Pi_{i=1}^n\mathbb{R}^{d }$ given by
\begin{equation}\label{eq:vector-field}
  F_i(q_1,\dots,q_n) =  q_i - \sum_{j,j\neq i} \frac{m_j}{r_{ij}^3}(q_i - q_j), \quad i=1,\dots,n.
  \index{$F_i$}
\end{equation}
Then the system (\ref{eq:cc-kart}) becomes
\begin{equation}\label{eq:cc-abstract}
  F(q_1,\dots,q_n)=0.
  \index{$F$}
\end{equation}

\subsection{Some identities and conservation laws}
\label{subsec:eq-conservation}
It is well know that for any $(q_1,q_2,q_3,\dots,q_n)\in (\reals^d)^n$ holds
\begin{eqnarray}
\sum_{i=1}^n f_i&=&0, \label{eq:n-mom-con} \\
\sum_{i=1}^n f_i \wedge q_i & = & 0, \label{eq:n-angular-mom-con}
\end{eqnarray}
where $v \wedge w$\index{$v \wedge w$} is the exterior product of vectors, the result being an element of exterior algebra. If $d=2,3$ it can be interpreted as the vector product of $v$ and $w$ in dimension $3$.  The identities (\ref{eq:n-mom-con}) and (\ref{eq:n-angular-mom-con})
are easy consequences of the third Newton's law (the action equals reaction) and the requirement that the mutual forces between bodies are in direction of the other body.

But (\ref{eq:n-mom-con}) and (\ref{eq:n-angular-mom-con})
  can be seen also as the consequences of the symmetries of Newtonian $n$-body problem. According to Noether's Theorem, by the translational symmetry we have a conservation
of momentum, which is equivalent to (\ref{eq:n-mom-con}), while the rotational symmetry implies   the conservation of angular momentum, which is implied by  (\ref{eq:n-angular-mom-con}).

 Note that the components of $v \wedge w$  are given by  determinants. In any dimension in the presence of the rotational symmetry, for any direction of rotation identified by $v_1 \wedge v_2$ ( $v_1$ and $v_2$ are perpendicular unit vectors) the following quantity must be zero (as a consequence of the Noether Theorem and the invariance with respect to the rotation in the plane $v_1,v_2$)
\begin{equation}\label{eq:conserved-ang-momentum}
  \sum_{i=1}^{n} \det \left[\begin{array}{cc}
                        (f_i|v_1) & (q_i|v_1) \\
                        (f_i|v_2) & (q_i|v_2)
                      \end{array}
                      \right]=0.
\end{equation}

Consider system (\ref{eq:cc-kart}). After multiplication of  $i$-th equation by $m_i$ and addition of all equations using (\ref{eq:n-mom-con}) we obtain (or rather recover)
the center of mass equation
\begin{eqnarray}
 \left(\sum_{i=1}^n m_i\right) c=\sum_i m_i q_i = 0. \label{eq:cc-cofmass}
\end{eqnarray}
We can take the equations for $n$-th body and replace it with (\ref{eq:cc-cofmass}) to obtain an equivalent system.
\begin{eqnarray}
  q_i&=& \sum_{j,j\neq i} \frac{m_j}{r_{ij}^3}(q_i - q_j)=:\frac{1}{m_i}f_i(q_1,\dots,q_n), \quad i=1,\dots,n-1, \label{eq:cc-kart-1n-1} \\
   q_n&=&-\frac{1}{m_n}\sum_{i=1}^{n-1} m_i q_i. \label{eq:cc-kart-n-th}
   \index{$q_i$}
\end{eqnarray}

Later in Section~\ref{sec:red-sys-equiv} we will use (\ref{eq:n-angular-mom-con}) to define a reduced system of equations for CCs which will
not have the degeneracies present in system (\ref{eq:cc-kart}).

\subsection{Moment of inertia of central configurations}
\label{subsec:mom-iner}

The important role of the moment of inertia in the investigation of central configurations is well known. In our context it plays a crucial role in stating some a priori bounds for central configurations.  We present, with proofs, some well known results on moment of inertia taken from  the notes by Moeckel \cite{Mlect2014} and the paper of Albouy and Kaloshin \cite{AK}.

\begin{definition}
For a configuration $q$  let the moment of inertia $I(q)$\index{$I(q)$} and the potential function $U(q)$\index{$U(q)$} be given by
  \begin{equation}
    I(q)=\sum_i m_i q_i^2, \quad U(q)=\sum_{i<j} \frac{m_i m_j}{r_{ij}}.
  \end{equation}
\end{definition}

\begin{lemma}
\label{lem:Irij}
Assume that $\sum_i m_i q_i=0$ and  $M=1$.
Then
\begin{equation}
I(q)= \sum_{i<j}m_i m_j (q_i - q_j)^2.  \label{eq:I-diff}
\index{$I(q)$}
\end{equation}
\end{lemma}
\textbf{Proof:}
Let us denote $q_{i,j}=q_i - q_j$.
Since
$$
\begin{array}{lrlp{5cm}}
q_i^2 = (q_i|q_i) & = & \left(\left(q_i - \sum_j m_j q_j \right)|q_i\right) & \hfill (since $\sum_i m_i q_i=0$)\\
  & = & \left(\left( q_i\sum_j m_j -  \sum_j m_j q_j\right)|q_i\right) & \hfill (since $\sum_j m_j = 1$)\\
  & = & \left(\left(\sum_j m_j(q_i - q_j)\right)|q_i\right) =  \sum_j m_j \left(q_{i, j}|q_i\right)
\end{array}
$$
we have
\begin{eqnarray*}
I(q)=\sum_i m_i q_i^2= \sum_{i,j} m_i m_j \left(q_{i,j}|q_{i} \right).
\end{eqnarray*}
Observe that
\begin{equation*}
  \sum_{i,j} m_i m_j \left(q_{i,j}|q_{i} \right) = - \sum_{i,j} m_i m_j \left(q_{i,j}|q_{j} \right),
\end{equation*}
hence
\begin{eqnarray*}
  I(q)=\frac{1}{2}\left(\sum_{i,j} m_i m_j \left(q_{i,j}|q_{i} \right)  - \sum_{i,j} m_i m_j \left(q_{i,j}|q_{j} \right) \right) =
  \frac{1}{2}\sum_{i,j} m_i m_j (q_{i,j}|q_{i,j})
\end{eqnarray*}
\qed

\begin{lemma}
\label{lem:I=U}
If $q \in \left(\mathbb{R}^d\right)^n$ is a (normalized) central configuration, then
\begin{equation}\label{eq:I=U}
  I(q)=U(q).
  \index{$I(q)$}
  \index{$U(q)$}
\end{equation}
\end{lemma}
\textbf{Proof:}
We take the scalar product of  $i$-th equation in (\ref{eq:cc-kart}) by $m_i q_i$ and add these equations to obtain
\begin{eqnarray*}
  I(q)=\sum_{i} m_i q_i^2 = \sum_{i,j; i \neq j} \frac{m_i m_j}{r_{ij}^3}(q_i - q_j|q_i) = \\
  \sum_{i<j} \frac{m_i m_j}{r_{ij}^3}(q_i - q_j|q_i-q_j)=  \sum_{i<j} \frac{m_i m_j}{r_{ij}}=U(q).
\end{eqnarray*}
\qed

\section{A priori bounds for central configurations}
\label{sec:ap-bounds}

From the point of view of CAP (computer assisted proofs)\index{CAP} in the problem of finding and counting all CCs  the issue of compactness of the search domain is fundamental.
The  lack of compactness arises for the following two reasons:
\begin{itemize}
\item two or more bodies might be arbitrary close to a collision,
\item some bodies might be arbitrary far from the origin.
\end{itemize}
The goal of this section is to deal with these issues. We will give  a priori bounds (depending on $m_i$'s) on the minimal distance of the closest bodies and for the size of the central configuration.

\subsection{Lower bound on the distances}
\label{subsec:lower-bnd-distances}

It is well known that central configurations are away from the collision set (see \cite{Sh}  or \cite[Prop. 15]{Mlect2014}).
However in these works no  quantitative statement directly applicable to system (\ref{eq:cc-kart}) has been given. Here we develop explicit a priori bounds.

The main idea is to use $I(q)=U(q)$ (see Lemma~\ref{lem:I=U}) to show that some term(s) $m_im_j/r_{ij}$ entering $U(q)$ dominate  and cannot be balanced, when bodies are very close.  Observe that using $I(q)=U(q)$ and positivity of all terms entering $U(q)$ allows us to escape the discussion of large terms on the rhs in the system (\ref{eq:cc-kart}), which
might cancel out or not etc. This is not the case in the framework of Albouy and Kaloshin \cite{AK}, where complex configurations and even complex masses have been considered, hence the positivity of $I$ and $U$ is lost.

\begin{lemma}\label{lm:iq-estm}
Assume that a (normalized) central configuration $q \in \left(\mathbb{R}^d\right)^n$ satisfies $|q_i| \leqslant  R$ for $i=1,\dots,n$. Then
$$I(q) \leqslant MR^2.$$
\end{lemma}
\textbf{Proof:\ }
Since for any $1 \leqslant i \leqslant n\colon |q_i| \leqslant  R$, thus
$
I(q) = \sum_{i} m_iq_i^2 \leqslant  \sum_{i} m_iR^2 = MR^2.
$
\qed

\begin{theorem}
\label{thm:l-bnd-dist}
Assume that a (normalized) central configuration $q \in \left(\mathbb{R}^d\right)^n$ satisfies $|q_i| \leqslant  R$ for $i=1,\dots,n$. Then
\begin{equation}
   r_{ij} > \frac{m_i m_j}{M R^2}, \quad 1\leqslant  i < j \leqslant  n. \label{eq:lb-dist}
\end{equation}
\end{theorem}
\textbf{Proof:}
From Lemma~\ref{lem:I=U} and Lemma~\ref{lm:iq-estm}, for any $1\leqslant  i < j \leqslant  n$, we have
\begin{eqnarray*}
M  R^2 \geqslant  I(q) = U(q) = \sum_{i < j}\frac{m_i m_j}{r_{ij}} > \frac{m_i m_j}{r_{ij}}.
\end{eqnarray*}
\qed

Below we establish a lower bound on the radius of ball centered at $0$ and containing a central configuration in the case of equal masses.
\begin{theorem}
\label{thm:cc-em-size}
  Assume that all masses are equal  and $q$ is a normalized central configuration such that $|q_i| \leqslant  R$.
Then
\begin{equation}
   R^3 \geqslant \frac{n-1}{4n}M.
\end{equation}
\end{theorem}
\textbf{Proof:}
Let $m=\frac{M}{n}$.
Since for any $1 \leqslant i \leqslant n\colon |q_i| \leqslant  R$, thus $r_{ij}\leqslant  2R$ and we obtain the following bound
\begin{eqnarray*}
  U(q) & = & \sum_{i < j} \frac{m_i m_j}{r_{ij}} \geqslant \sum_{i < j} \frac{m^2}{2R} \\
  & = &  \frac{m^2}{2R}\cdot \frac{n(n-1)}{2} =\frac{(n-1)M^2}{4nR}.
\end{eqnarray*}
Hence from Lemma~\ref{lem:I=U} and Lemma~\ref{lm:iq-estm} we obtain
\begin{eqnarray*}
  \frac{(n-1)M^2}{4nR} \leqslant  U(q)=I(q) \leqslant M R^2.
\end{eqnarray*}
\qed

If $M=1$, then $\lim_{n\to\infty} \sqrt[3]{\frac{n-1}{4n}} = 4^{-1/3}\approx 0.629961.$
This estimate appears to be reasonably good, as shown in Table~\ref{tab:minR}.
\begin{table}[htb]
  \centering
\begin{tabular}{c|l|l}
  $n$ & $\sqrt[3]{\frac{n-1}{4n}}$ & $R$ \\
  \hline
  $3$ & 0.550321 & $0.577350$ \\
  $4$ & 0.572357 & $0.620813$ \\
  $5$ & 0.584804 & $0.650513$ \\
  $6$ & 0.592816 & $0.672798$ \\
\end{tabular}
  \caption{The size of a minimal ball containing all normalized central configurations with $M=1$ for several $n$'s for equal masses case.  The minimum $R$ is realized for regular $n$-gon. }\label{tab:minR}
\end{table}

In the next theorem we do not assume that all masses are equal.
\begin{theorem}
\label{thm:cc-dist>1}
Assume $M=1$ and $q$ is a normalized central configuration. Then there exists a pair $i \neq j$ such that
\begin{equation}
  r_{ij} \geqslant  1.
\end{equation}
\end{theorem}
\textbf{Proof:}
For the proof by contradiction assume that $r_{ij} <1$ for all pairs of bodies. Hence we have $r_{ij}^2 < 1/r_{ij}$
for all pairs. From Lemma~\ref{lem:Irij} it follows
\begin{eqnarray*}
 I(q)=\sum_{i<j} m_i m_j r_{ij}^2 < \sum_{i<j} m_i m_j \frac{1}{r_{ij}}=U(q).
\end{eqnarray*}
From Lemma~\ref{lem:I=U} it follows that $q$ is not a central configuration.
\qed

Observe that the above estimate is optimal, because it is realized for the equilateral triangle for $n=3$ and a tetrahedron (non planar CC) for $n=4$.
From the above theorem we obtain the following lower bound for the size of a central configuration. Contrary to Theorem~\ref{thm:cc-em-size}
we do not assume that all masses are equal.
\begin{theorem}
\label{thm:cc-size-lowerbnd}
Assume that $M=1$ and $q$ is a normalized central configuration and $|q_i| \leqslant  R$ for $i=1,\dots,n$. Then $R \geqslant  1/2$.
\end{theorem}
\textbf{Proof:}
For the  proof by contradiction assume that $R < 1/2$. Then for all pairs $r_{ij} <1$. From Theorem~\ref{thm:cc-dist>1} it follows that $q$ is not central configuration.
\qed

\subsection{The upper bound on the size of central configuration}
\label{subsec:ap-up-bnds}

The goal of this section is to give the upper bounds for the size of the central configuration. This time we exploit the fact that
if the forces are bounded, then large  $q_i$'s on the left hand side of the system~(\ref{eq:cc-kart}) cannot be balanced. The obvious difficulty with the realization of this idea
is: we  can have a group of bodies with large norms  which are close to each other in the central configuration, which produce large terms on rhs of the system (\ref{eq:cc-kart}). To overcome this we consider clusters of
points far from the origin and the resulting force on it. In such situation mutual interactions between bodies in the cluster cancel out.

\begin{lemma}
\label{lem:upp-bnd}
Assume $q \in \left(\mathbb{R}^d\right)^n$ is a  normalized central configuration.
Let  $R=|q_{i_0}|=\max_{i=1,\dots,n} |q_i|$. Then  for all  $\varepsilon  \in \left(0,R/(n-1)\right)$ holds
\begin{equation}
  R - (n-2) \varepsilon  < \frac{M}{\varepsilon ^2}.
\end{equation}
\end{lemma}
\textbf{Proof:} For simplicity let's assume that $d = 2$ and $q_i = (x_i, y_i)$.
Let us fix any $\varepsilon \in \left(0,R/(n-1)\right)$.
After a  suitable rotation of coordinate system  we can assume that
\begin{equation}
q_{i_0}=(R,0).  \label{eq:maxqi-on-OX}
\end{equation}
Let $\mathcal{C}$\index{$\mathcal{C}$} be a minimal subset (cluster)\index{cluster} of indices of bodies satisfying the following conditions
\begin{itemize}
\item $i_0 \in \mathcal{C}$
\item if $j \in \mathcal{C}$ and $|q_k - q_j| \leqslant  \varepsilon $, then $k \in \mathcal{C}$
\end{itemize}
The cluster $\mathcal{C}$ can be constructed as follows: We  start with $i_0 \in \mathcal{C}$. Then we add all bodies which are not farther than $\varepsilon $
from the bodies already in $\mathcal{C}$. We repeat this until the set $\mathcal{C}$ stabilizes, which will happen after at most $n-1$ steps.
From  assumption about $\varepsilon$ and $R$ it follows that
\begin{equation}
  R > (n-1)\varepsilon . \label{eq:R>neps}
\end{equation}
Observe that (\ref{eq:R>neps}) implies that
$\mathcal{C} \neq \{1,\dots,n\}$. Indeed (\ref{eq:R>neps}) and (\ref{eq:maxqi-on-OX}) imply that $x_i >0$ for all $i \in \mathcal{C}$. This and the center of mass condition (\ref{eq:cc-cofmass}) implies that $\mathcal{C}$ cannot contain all bodies.
This implies that  the process of building $\mathcal{C}$ must stop after at most $n-2$ steps.
Therefore we obtained a cluster $\mathcal{C}$ with the following properties
\begin{eqnarray}
  q_i &\in& \overline{B(q_{i_0},(n-2)\varepsilon )}, \quad \forall i \in \mathcal{C}, \label{eqn:cluster-ball}\\
  |q_i - q_j| &>& \varepsilon, \quad i \in \mathcal{C}, j \notin \mathcal{C}.
\end{eqnarray}

\begin{figure}[htb]
  \centering
\begin{tikzpicture}
\begin{scope}
\clip (3.5,2.0) circle(1.0);           \clip (2.0,2.0) circle(1.5);
\fill[gray!20] (2.5,0.0) rectangle (3.5,4.0);
\end{scope}

\draw[->] (0.0, 2.0) -- (5.0, 2.0);
\draw[->] (2.0, 0.0) -- (2.0, 4.0);

\draw[] (2.0, 2.0) circle (1.5cm);
\draw[->] (2.0, 2.0) -- (2.75, 3.3);

\draw[fill] (3.5, 2.0) circle (0.4mm);
\node[] at (3.75,1.8) {\small $q_{i_0}$};

\draw[] (3.5, 2.0) circle (1.0cm);
\draw[->] (3.5, 2.0) -- (3.95, 2.85);

\node[] at (2.25,2.8) {\small $R$};
\node[] at (3.9,2.4) {\small $r$};

\end{tikzpicture}
\caption{$R > (n-1)\varepsilon$ and $r = (n-2)\varepsilon$; the darkened area is the region where all the bodies from cluster are located.}\label{fig:cluster}
\end{figure}
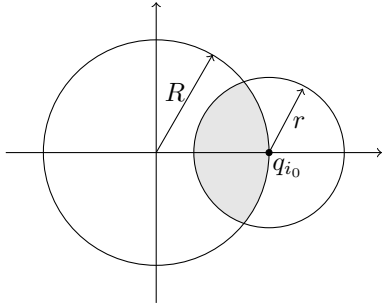

Without any loss of the generality we can assume that $i_0=1$ and $\mathcal{C}=\{1,\dots,s\}$.
Note that for any $k\geqslant 2$ there is (compare~(\ref{eq:n-mom-con}))
\begin{equation*}
\sum_{i = 1}^k\sum_{j = 1, j\neq i}^k \frac{m_im_j}{r_{ij}^3}(q_i - q_j) = 0
\end{equation*}
thus by adding equations~(\ref{eq:cc-kart}) multiplied by $m_i$ for $i=1,\dots,s$ we obtain
\begin{eqnarray}
  m_1 q_1 + \dots + m_s q_s & =  & \sum_{\substack{i=1,\dots,s,\\
j=1,\dots, s,\\ j\neq i}}\frac{m_i m_j (q_i - q_j)}{r_{ij}^3}\ +   \sum_{\substack{i=1,\dots,s,\\
j=s+1,\dots,n}}\frac{m_i m_j (q_i - q_j)}{r_{ij}^3}\nonumber\\
\nonumber\\
& = &
  \sum_{\substack{i=1,\dots,s,\\
j=s+1,\dots,n}}\frac{m_i m_j (q_i - q_j)}{r_{ij}^3}. \label{eq:sum-s-eq}
\end{eqnarray}

Let
\begin{eqnarray*}
  M_s&=&\sum_{i=1}^s m_i, \\
  c_{s}&=& \frac{1}{M_s}\sum_{i=1}^s m_i q_i, \\
  F_{s}&=& \frac{1}{M_s}\sum_{i\in \mathcal{C}; j \notin \mathcal{C}} \frac{m_i m_j (q_i - q_j)}{r_{ij}^3}.
\end{eqnarray*}

Observe that (\ref{eq:sum-s-eq}) could be now rewritten as
\begin{equation}
   c_s=F_s. \label{eq:cc-cluster-eq}
\end{equation}

It is easy to see (cf.~(\ref{eqn:cluster-ball})) that for $i\in\mathcal{C}$: $x_i \geqslant  R - (n-2) \varepsilon > 0$, hence
\begin{eqnarray*}
  |c_s| & = & \frac{1}{M_s}\left|\sum_{i=1}^s m_i q_i\right| \geqslant  \frac{1}{M_s}\left|\sum_{i=1}^s m_i(R - (n-2)\varepsilon, y_i)\right|\nonumber\\
  & \geqslant & \frac{1}{M_s}(R - (n-2)\varepsilon)\sum_{i=1}^s m_i\nonumber\\
  & = &  R - (n-2) \varepsilon ,   \label{eq:estm-Cs}
\end{eqnarray*}
and
\begin{eqnarray*}
|F_s| &\leqslant &
 \frac{1}{M_s}\sum_{i\in \mathcal{C}; j \notin \mathcal{C}}  \frac{m_i m_j }{r_{ij}^2} \leqslant
 \frac{1}{\varepsilon ^2} \frac{1}{M_s}\sum_{i\in \mathcal{C}; j \notin \mathcal{C}}  m_i m_j\\
  & = & \frac{1}{\varepsilon ^2} \frac{1}{M_s}\left(\sum_{i\in \mathcal{C}} m_i \right) \cdot \left(\sum_{j \notin \mathcal{C}}   m_j \right)=
  \frac{1}{\varepsilon ^2} \sum_{j \notin \mathcal{C}}   m_j <  \frac{M}{\varepsilon ^2}.
\end{eqnarray*}
Hence from the above and (\ref{eq:cc-cluster-eq}) we obtain
\begin{eqnarray*}
  R - (n-2) \varepsilon  \leqslant  |c_s| = |F_{s}| < \frac{M}{\varepsilon ^2}.
\end{eqnarray*}
\qed

From Lemma~\ref{lem:upp-bnd} we obtain the following estimate on the size of any central configuration.
\begin{theorem}
\label{thm:upp-bnd}
Assume that $M=1$ and $q \in \left(\mathbb{R}^d\right)^n$ is a normalized central configuration. Then
\renewcommand{\arraystretch}{1.6}
\begin{equation}
  \max_{i} |q_i| \leqslant  \left\{
                                   \begin{array}{ll}
                                     n-1, & \hbox{$n\geqslant  2$;} \\
                                     \left(2^{1/3}+2^{-2/3} \right) (n-2)^{2/3}, & \hbox{$n\geqslant  4$.}
                                   \end{array}
                                 \right.
\end{equation}
\renewcommand{\arraystretch}{1.2}
\end{theorem}
\textbf{Proof:}
Let  $R=\max_{i=1,\dots,n} |q_i|$.
From Lemma~\ref{lem:upp-bnd} it follows that for any $\varepsilon >0$  holds
\begin{equation}
  R \leqslant  \max \left((n-1)\varepsilon , (n-2)\varepsilon  + \frac{1}{\varepsilon ^2} \right). \label{eq:upp-bound-max}
\end{equation}
Indeed, if $\varepsilon  < R/(n-1)$ then we apply Lemma~\ref{lem:upp-bnd}, otherwise we have $R \leqslant  (n-1)\varepsilon $.
Let us optimize the bound (\ref{eq:upp-bound-max}) with respect to $\varepsilon $.
Let us denote
\begin{equation}
  B(\varepsilon )=\max \left((n-1)\varepsilon , (n-2)\varepsilon  + \frac{1}{\varepsilon ^2} \right).
\end{equation}
It is easy to see that
\begin{eqnarray*}
  (n-1)\varepsilon  &<&  (n-2)\varepsilon  + \frac{1}{\varepsilon ^2}, \quad \mbox{for} \, \varepsilon  < 1, \\
  (n-1)\varepsilon  &=&  (n-2)\varepsilon  + \frac{1}{\varepsilon ^2}, \quad \mbox{for} \, \varepsilon  = 1, \\
  (n-1)\varepsilon  &>&  (n-2)\varepsilon  + \frac{1}{\varepsilon ^2}, \quad \mbox{for} \, \varepsilon  > 1.
\end{eqnarray*}
Hence
\renewcommand{\arraystretch}{1.6}
$$
B(\varepsilon) = \left\{
\begin{array}{ll}
(n-2)\varepsilon  + \frac{1}{\varepsilon ^2}, & \mbox{for} \, \varepsilon  < 1, \\
(n-1)\varepsilon, & \mbox{for} \, \varepsilon  \geqslant  1.
\end{array}
\right.
$$
\renewcommand{\arraystretch}{1.2}
Therefore
\begin{equation*}
  \inf_{\varepsilon  >0} B(\varepsilon )= \inf_{\varepsilon  \in (0,1]}\left((n-2)\varepsilon  + \frac{1}{\varepsilon ^2}\right).
\end{equation*}
The function $g(\varepsilon )=(n-2)\varepsilon  + 1/\varepsilon ^2$ is decreasing for $\varepsilon  < \varepsilon _0=\left( \frac{2}{n-2}\right)^{1/3}$ and
is increasing  for $\varepsilon > \varepsilon _0$.  Observe that $\varepsilon _0 \in (0,1]$ iff $n \geqslant  4$.
 Therefore we obtain
\begin{equation}
  \inf_{\varepsilon  >0} B(\varepsilon )=\left\{
                                   \begin{array}{ll}
                                     n-1, & \hbox{$n\leqslant 4$;} \\
                                     g(\varepsilon _0), & \hbox{$n\geqslant  4$.}
                                   \end{array}
                                 \right.
\end{equation}
We have
\begin{equation*}
  g(\varepsilon _0) = (n-2)\left( \frac{2}{n-2}\right)^{1/3} + \left( \frac{n-2}{2}\right)^{2/3}=\left(2^{1/3}+2^{-2/3} \right) (n-2)^{2/3}
\end{equation*}
\qed

For $n\leqslant  10$ and the equal masses  the above estimate appears to be an overestimation. The largest possible size of CC found was slightly above $1$ and is considerably smaller than the one established in the above theorem.

\section{Exclusion tests for CCs}
\label{sec:excl-tests}

Assume that we have an interval set $D$ (i.e. a box, a product of intervals) in which we would like to exclude the existence of CC. We do not assume that $D \subset \dom F$ (see (\ref{eq:cc-abstract})) and this is an important point.  The a priori estimates discussed in Section~\ref{sec:ap-bounds} allow to exclude $D$
iff there is no point in $D$ which satisfies the obtained  bounds.

In the following we discuss other exclusion tests.

\subsection{Checking for zeros}
\label{subsec:check-for-zeros}
One obvious test is to check whether in the interval evaluation of $F(D)$ (see (\ref{eq:cc-abstract})) at least one of the component does not contain zero.
Observe that a partial tests of this type are possible even if $D$ admits  the collision, i.e. formally it is not contained in $\dom F$, but we can do this verification whenever $D \subset \dom F_i$.
The test takes the following form:
\begin{itemize}
\item[0.] given box $D$ in the configuration space, such that $D \subset \dom F_i$ (i.e. the $i$-th does not have a collision with other bodies). Let $D_i=\{q_i \,|\, q \in D\}$.

\item[1.] compute  the interval enclosure of $f_i(q)$ for $q \in D$, denoted by $\langle f_i(D)\rangle$,
\item[2.] if $D_i \cap \frac{1}{m_i}\langle f_i(D)\rangle=\emptyset$  (compare equation (\ref{eq:cc-kart})), then $D$ does not contain any normalized central configuration
\item[3.] if $D_i \cap \frac{1}{m_i}\langle f_i(D)\rangle\neq \emptyset$, then we define $D^{\mathrm{ref}}=\{q \in D\, |\, q_i \in \frac{1}{m_i}\langle f_i(D)\rangle \}$.
\end{itemize}
Observe that point 2 allow us exclude the box $D$, but if this is impossible then we can replace $D$ by $D^{\mathrm{ref}}$ obtained in point 3, as any CC
contained in $D$ must belong to $D^{\mathrm{ref}}$. This is one of several places in our algorithm, where we attempt to do better than doing naive binary subdivision in order to relieve the curse of dimensionality.

\subsection{The cluster of bodies - checking for zeros test }
\label{subsec:cluster-check-for-zeros}

In the case of colliding bodies or a cluster of close bodies, say with indices $i=1,\dots, s$, after adding first $s$ equations multiplied by $m_i$ we obtain
\begin{equation*}
  \sum_{i=1}^s m_i q_i = \sum_{i=1}^s \sum_{j=s+1}^n \frac{m_i m_j}{r_{ij}^3}(q_i - q_j).
\end{equation*}
Therefore for the cluster of bodies $\mathcal{C} \subset \{1,\dots,n\}$ we will check whether
\begin{equation}
   (0,0) \in \sum_{i \in \mathcal{C}}^s m_i q_i - \sum_{i\in \mathcal{C}}^s \sum_{j\notin \mathcal{C}}^n \frac{m_i m_j}{r_{ij}^3}(q_i - q_j), \label{eq:zero-for-claster}
\end{equation}
where the expression on the rhs of (\ref{eq:zero-for-claster}) is evaluated in the interval arithmetic on $D$.

If it is not satisfied then we conclude that $D$ does not contain any CC. Again let us stress that the set $D$ might contain some collisions,
and this test is still applicable.

\subsection{The cluster of bodies - test based on moment of inertia and potential }
\label{subsec:cluster-mom-pot}
We will  exploit $I(q)=U(q)$ for CC (compare Lemma~\ref{lem:I=U}), but our focus will be on  the subsets (clusters) of bodies.
Let us fix $\mathcal{C} \subset \{1,2,\dots,n\}$ and $Z \subset \left(\mathbb{R}^d\right)^n $.
Let us define
\begin{eqnarray*}
  U_{\mathcal{C},Z}&=&\inf_{q \in Z}   \sum_{i<j, i,j \in \mathcal{C}} \frac{m_i m_j}{r_{ij}}, \\
  I_{\mathcal{C},Z}&=&\sup_{q \in Z}   \sum_{i \in \mathcal{C}} m_i q_i^2, \\
  F_{\mathcal{C},Z}&=&\inf_{q \in Z}   \sum_{ i \in \mathcal{C}, k \notin \mathcal{C}} \frac{m_i m_k}{r_{ik}^3} (q_i-q_k| q_i).
\end{eqnarray*}
In the case when $\mathcal{C}=\{1,\dots,n\}$ we set $F_{\mathcal{C},Z}=0$.

The important point is that we can compute the infimum in $U_{\mathcal{C},Z}$ even if the set $Z$ contains collisions.  It makes sense to take as $\mathcal{C}$ a cluster of close points (containing possible collisions and near collisions), so that there is no collision between bodies in $\mathcal{C}$ and its complement. In such case $F_{\mathcal{C},Z}$ will be finite.

We have the following criterion for nonexistence of CC in $Z$:
\begin{lemma}\label{lm:FUI}
Assume that $ U_{\mathcal{C},Z},  I_{\mathcal{C},Z},  F_{\mathcal{C},Z} \in \mathbb{R}$ and
\begin{equation*}
  I_{\mathcal{C},Z} < U_{\mathcal{C},Z} + F_{\mathcal{C},Z},
\end{equation*}
then there is NO central configuration in $Z$.
\end{lemma}
\textbf{Proof:}
Without any loss of the generality we can assume that $\mathcal{C}=\{1,\dots,s\}$ for some $1 \leqslant  s \leqslant n$.
Consider system (\ref{eq:cc-kart}). We multiply  $i$-th equation by  $m_i q_i$ and we add first $s$ equations.
We obtain (compare the proof of Lemma~\ref{lem:I=U})
\begin{eqnarray}
\label{eq:cluster-I=U+cos}
  \sum_{i=1}^s m_i q_i^2 = \sum_{1\leqslant  i< j \leqslant  s} \frac{m_i m_j}{r_{ij}} + \sum_{i=1}^s \sum_{j>s} \frac{m_i m_j (q_i - q_j|q_i)}{r_{ij}^3}.
\end{eqnarray}
Let us stress that (\ref{eq:cluster-I=U+cos}) must hold for any central configuration.
Now for $q \in Z$ holds
\begin{eqnarray*}
  \sum_{1\leqslant  i< j \leqslant  s} \frac{m_i m_j}{r_{ij}} + \sum_{i=1}^s \sum_{j>s} \frac{m_i m_j (q_i - q_j|q_i)}{r_{ij}^3} \geqslant  \\
  U_{\mathcal{C},Z} + F_{\mathcal{C},Z} >   I_{\mathcal{C},Z} =  \sum_{i=1}^s m_i q_i^2.
\end{eqnarray*}
Hence  (\ref{eq:cluster-I=U+cos}) is not satisfied. Therefore  we do not have any central configuration in $Z$.
\qed

 Observe that if $\mathcal{C}=\{1,\dots,n\}$ the above lemma is reduced to checking whether $\inf_{q \in Z} U(q) > \sup_{q \in Z} I(q) $.

\section{The reduced system of equations for CC on the plane}
\label{sec:red-sys-equiv}

\subsection{Non-degenerate solutions of full and reduced systems of equations}

Following Moeckel \cite{Mlect2014} we state the definition.
\begin{definition}
  \label{def:non-deg-cc} We will say that a normalized central configuration $q=(q_1,\dots,q_n)$ is \emph{non-degenerate}
  if the rank of $D\!F(q)$ is equal to $dn-\dim SO(d)$. Otherwise the configuration will be called \emph{degenerate}.
\end{definition}
The idea of the above notion of degeneracy is to allow only for the  degeneracy related to the rotational symmetry of the problem, because
by setting $\lambda=1$ in (\ref{eq:cc-with-lambda}) and keeping the masses fixed we removed the scaling symmetry.

The  system (\ref{eq:cc-kart-1n-1}--\ref{eq:cc-kart-n-th}) obtained from (\ref{eq:cc-abstract}) after removing the $n$-th body using the center of mass (condition (\ref{eq:cc-cofmass})) we  write as
\begin{equation}\label{eq:cc-abstract-red}
  F_{\mathrm{red}}(q_1,\dots,q_{n-1})=0
  \index{$F_{\mathrm{red}}$}
\end{equation}
where $F_{\mathrm{red}}: \Pi_{i=1}^{n-1}\mathbb{R}^{d} \to \Pi_{i=1}^{n-1}\mathbb{R}^{d} $.
Then it is easy to see that $q=(q_1,\dots,q_{n-1},q_n)$ is a non-degenerate central configuration iff the rank
of $D\! F_{\mathrm{red}}(q_1,\dots,q_{n-1})$ is $d(n-1)-\dim SO(d)$.

\subsection{The reduced system on the plane}

We consider a planar case here, i.e. $d=2$.
The fact that the system of equations (\ref{eq:cc-kart}) is degenerate (each solution give rise to a circle of solutions) make this system not amenable for the use of standard interval arithmetic methods (see for example the Krawczyk operator discussed in Section~\ref{subsec:Krawczyk})) to rigorously count all possible solutions.
We need to kill the $SO(2)$-symmetry and then hope that all solutions will be non-degenerate. In this section we show how to reduce the system (\ref{eq:cc-kart}) to an equivalent system amenable to the interval analysis tools.

Let us fix $k \in \{1,\dots,n-1\}$ and consider the following set of equations
\begin{eqnarray}
  q_i &=&\frac{1}{m_i}f_i(q_1,\dots,q_n(q_1,\dots,q_n)), \quad i\in \{1,\dots,n-1\}, i \neq k \label{eq:cc-red-i} \index{$q_i$}\\
  x_k &=& \frac{1}{m_k}f_{k,x}(q_1,\dots,q_n(q_1,\dots,q_n)), \label{eq:cc-red-x2}\\
  q_n&=&-\frac{1}{m_n}\sum_{i=1}^{n-1} m_i q_i, \label{eq:cc-red-com}
\end{eqnarray}
where $f_i = (f_{i,x}, f_{i, y})$.\index{$f_{k,x}$}
Observe that the system (\ref{eq:cc-red-i}--\ref{eq:cc-red-com}) coincides with the system (\ref{eq:cc-kart-1n-1}--\ref{eq:cc-kart-n-th}) with the
equation for $y_k$ dropped.

The next theorem addresses   the question: whether from the reduced system  (\ref{eq:cc-red-i}--\ref{eq:cc-red-com})  we obtain the solution of (\ref{eq:cc-kart})?

\begin{theorem}
\label{thm:red-to-full}
  Assume $d = 2$. If $(q_1,\dots,q_n)$ is a solution of the reduced system (\ref{eq:cc-red-i}--\ref{eq:cc-red-com})  and $x_k \neq x_n$, then
    it is a normalized central configuration, i.e. it satisfies (\ref{eq:cc-kart}).
\end{theorem}
\textbf{Proof:}
Let $q_i$ be as in our assumptions.
We need to show that  $m_k y_k=f_{k,y}$.
From (\ref{eq:cc-red-com}), (\ref{eq:n-mom-con}) and (\ref{eq:n-angular-mom-con}) it follows that
\begin{eqnarray*}
 0=\sum_{i=1}^n f_i \wedge q_i =
 \sum_{i=1}^{n-1}f_i \wedge q_i +  \left(-\sum_{i=1}^{n-1}f_i \right) \wedge \left(\frac{-\sum_{i=1}^{n-1}m_i q_i}{m_n} \right) = \\
\nonumber   \sum_{i=1}^{n-1}f_i \wedge q_i \left(1 + \frac{m_i}{m_n} \right)  + \frac{1}{m_n} \sum_{i,j=1, i\neq j}^{n-1}m_j f_i \wedge q_j
\end{eqnarray*}
Since from (\ref{eq:cc-red-i})
\begin{eqnarray*}
 q_i \wedge f_i&=&0, \quad i=1,\dots,n-1; \quad i \neq k
 \end{eqnarray*}
 we obtain
 \begin{eqnarray}
 0= f_k \wedge q_k \left(1 + \frac{m_k}{m_n} \right)  + \frac{1}{m_n} \sum_{i,j=1, i\neq j}^{n-1}m_j f_i \wedge q_j. \label{eq:n-an-mom-red}
 \end{eqnarray}
Let us take the look at
 $\sum_{i,j=1, i\neq j}^{n-1}m_j f_i \wedge q_j $. We have from (\ref{eq:cc-red-i})
 \begin{equation*}
   \sum_{i,j=1, i\neq j, i \neq k, j\neq k}^{n-1}m_j f_i \wedge q_j = \sum_{i,j=1, i\neq j, i \neq k, j\neq k}^{n-1} f_i \wedge f_j =0,
 \end{equation*}
 hence again from (\ref{eq:cc-red-i}) it follows that
\begin{eqnarray*}
 \sum_{i,j=1, i\neq j}^{n-1}m_j f_i \wedge q_j =\sum_{i=1,i \neq k}^{n-1} \left(m_k f_i \wedge q_k + m_i f_k \wedge q_i\right) = \\
 \sum_{i=1,i \neq k}^{n-1} \left(m_k m_i q_i \wedge q_k - m_i q_i \wedge f_k\right)=\left(\sum_{i=1,i \neq k}^{n-1} m_i q_i \right)\wedge  \left(m_k  q_k - f_k\right).
\end{eqnarray*}
 From the above and (\ref{eq:n-an-mom-red}) we obtain
 \begin{eqnarray*}
 0=  f_k \wedge q_k \left(1 + \frac{m_k}{m_n} \right) + \frac{1}{m_n}\left(\sum_{i=1,i \neq k}^{n-1} m_i q_i \right)\wedge  \left(m_k  q_k - f_k\right).
 \end{eqnarray*}

From (\ref{eq:cc-red-x2})  we have
 \begin{eqnarray*}
     0=f_k \wedge q_k \left(1 + \frac{m_k}{m_n} \right) + \frac{1}{m_n}\left(\sum_{i=1,i \neq k}^{n-1} m_i q_i \right)\wedge  \left(m_k  q_k - f_k\right)= \\
     (f_{k,x}y_k - f_{k,y}x_k)\left(1 + \frac{m_k}{m_n} \right) + \frac{1}{m_n}\left(\sum_{i=1,i \neq k}^{n-1} m_i x_i \right)  \left(m_k  y_k - f_{k,y}\right)= \\
     (m_k x_k y_k - f_{k,y}x_k)\left(1 + \frac{m_k}{m_n} \right) + \frac{1}{m_n}\left(\sum_{i=1,i \neq k}^{n-1} m_i x_i \right)  \left(m_k  y_k - f_{k,y}\right)= \\
     \left(m_k y_k - f_{k,y}\right) \left(x_k\left(1 + \frac{m_k}{m_n} \right) +   \frac{1}{m_n}\left(\sum_{i=1,i \neq k}^{n-1} m_i x_i \right)  \right)= \\
     \frac{1}{m_n} \left(m_k y_k - f_{k,y}\right) \left(m_n x_k + \left(\sum_{i=1}^{n-1} m_i x_i \right) \right) = \\
      \frac{1}{m_n} \left(m_k y_k - f_{k,y}\right) \left(m_n x_k -m_n x_n \right)= \left(m_k y_k - f_{k,y}\right) \left( x_k - x_n \right).
 \end{eqnarray*}
 Now if $x_k - x_n \neq 0$, then $m_k  y_k = f_{k,y} $.
\qed

The system (\ref{eq:cc-red-i}--\ref{eq:cc-red-x2}) contains $2(n-1)-1$ equations in $2(n-1)$ variables and has  $O(2)$ symmetry (i.e. rotations
around origin and reflection symmetries with respect to lines passing through the origin map solutions of this system into itself). In order to obtain a system with the same number of equations and variables we can impose additional condition leading to the removal of $y_k$ variable, so that the rotational symmetry will be broken. Obviously in the above setting we could drop the equation for $x_k$ and we will obtain an analogous result.

 We can think of a general reduced system as follows:
\begin{itemize}
\item we fix some hyperplane $H$, in the reduced (by the center of mass condition) configuration space $\mathbb{R}^{2(n-1)}$, so that $H$ is transversal to the action $SO(2)$ and  $k$  is such that $v_k \in \{x_k,y_k\}$ can be computed in terms of other variables. This will induce an embedding, $J_{k}:\mathbb{R}^{2(n-1)-\dim SO(2)} \to H$.
\item in the system (\ref{eq:cc-kart-1n-1}--\ref{eq:cc-kart-n-th}) we remove the equation for  $v_{k}$.
Then the reduced system can be written as
\begin{equation}\label{eq:gen-red-system}
  R_{k} F_{\mathrm{red}} (J_{k}z)=0,
\end{equation}
where $R_k$ is a projection which removes $v_k$ variable in the image.
\end{itemize}
The system (\ref{eq:cc-red-i}--\ref{eq:cc-red-com}) supplemented by substitution $y_k=y_k(\dots)$ is an example of (\ref{eq:gen-red-system}).
We present the following obvious result
\begin{theorem}
Assume for simplicity that $H$ is given by
\begin{equation}
y_{k}=\sum_{i \neq k, i\leqslant   n-1} a_i y_i + \sum_{i\leqslant   n-1} b_i y_i, \label{eq:yk-norm}
\end{equation}
Assume that $q$ is a solution of reduced system  (\ref{eq:gen-red-system}) with a substitution (\ref{eq:yk-norm}), such that $x_k\neq x_n$.  Then $q$ is a normalized central configuration.
If $q$ is  non-degenerate solution of the reduced system, then
this is non-degenerate solution of (\ref{eq:cc-kart}).
\end{theorem}
\textbf{Proof:}
The first part is obvious in view of Theorem~\ref{thm:red-to-full} and condition $x_k \neq x_n$ implies that it is a central configuration. The maximum rank in the reduced system
gives the non-degeneracy of the configuration in the sense of Definition~\ref{def:non-deg-cc}.
\qed

Following \cite{AK} we have tried two possibilities
\begin{itemize}
\item we set $k=2$ and  we eliminate variable $y_2$ by setting
 \begin{equation}
  y_2=y_1,  \label{eq:yk=y1-norm}
\end{equation}
\item we set $k=n-1$ and we eliminate variable $y_{n-1}$ by setting
\begin{equation}
  y_{n-1}=0. \label{eq:yn-1=0}
\end{equation}
\end{itemize}
Observe that for both normalizations defined above for any CC $q$ there is a rotation  $R$  such that $Rq$ satisfies this normalization. Hence we can safely impose any of those conditions without losing any CC.
In both cases we will need
\begin{equation}
x_k \neq x_n.  \label{eq:xk-neq-xn}
\end{equation}
First consider (\ref{eq:yk=y1-norm}). Condition (\ref{eq:xk-neq-xn}) does not hold for some CC in the case of equal masses. For example, for $n=4$ and an equilateral rectangle, such that $x_2=x_4$ we obtained numerically (and also symbolically using Mathematica) that the jacobian matrix for the reduced system has a zero eigenvalue. Hence the solution is degenerate for the reduced system.

Now, consider the condition (\ref{eq:yn-1=0}). If we setup our computations so that  $q_{n-1}$ body maximizes the distance from the origin for all bodies, then we have (\ref{eq:xk-neq-xn}) satisfied, otherwise $q_n$ will be further from zero.  This observation does not prove that if $q$ is
a non-degenerate CC in the sense of Definition~\ref{def:non-deg-cc}, then it is also a non-degenerate solution of the reduced system, but this appears to happen in our rigorous computation of central configurations so far.

\section{On the computer assisted proof}
\label{sec:cap}

We normalized masses so that $M=\sum_i m_i=1$.
In this section we index bodies from $0$ to $n-1$ to be in the agreement with our program.
In the sequel we study the following reduced system
\begin{eqnarray}
  q_i &=&\frac{1}{m_i}f_i(q_0,\dots,q_{n-2},q_{n-1}(q_0,\dots,q_{n-2})), \quad i\in \{0,\dots,n-3\},  \label{eq:cc-red-final-i} \\
  x_{n-2} &=& \frac{1}{m_{n-2}}f_{n-2,x}(q_0,\dots,q_{n-2},q_{n-1}(q_0,\dots,q_{n-2})), \label{eq:cc-red-final-x2}
\end{eqnarray}
where we set
\begin{eqnarray}
  y_{n-2}&=&0, \label{eq:cc-red-final-yn-1}  \\
  q_{n-1}(q_0,\dots,q_{n-2})&=&-\frac{1}{m_{n-1}}\sum_{i=0}^{n-2} m_i q_i, \label{eq:cc-red-final-com}.
\end{eqnarray}

\subsection{Equal masses case, the reduction of the configuration space for CCs}
\label{subsec:cap-equal-masses}

In the case of equal masses, after a suitable rotation and permutation of the bodies, we can assume that
\begin{equation}
  |x_{n-2}| \geqslant |q_i| , \quad i=0,\dots,n-1. \label{eq:xn>qi}
\end{equation}
Condition (\ref{eq:xn>qi}) guarantees that $x_{n-2} \neq x_{n-1}$, hence by Theorem~\ref{thm:red-to-full} the solution of a reduced system (\ref{eq:cc-red-final-i}--\ref{eq:cc-red-final-x2}) is CC.
In view of symmetry and Lemma~\ref{thm:cc-size-lowerbnd} we impose some more conditions
\begin{equation}
 n-1 \geqslant x_{n-2} \geqslant 0.5. \label{eq:xn-1bnds}
\end{equation}
After a suitable permutation of bodies and a reflection with respect to $0X$-axis it is easy to see that each CC has its equivalent
in the set of the configurations satisfying the following conditions
\begin{itemize}
\item $q_{n-2}=(R,0)$ is the furthermost body from the origin
\item $q_0$ is the leftmost with non-negative $y$-coordinate
\item $q_1$ has the smallest $y$ coordinate
\item all other bodies have their $x$ coordinates in the order of increasing/decreasing indices.
\end{itemize}
This, combined with  Lemma~\ref{lem:upp-bnd}, shows that it is enough  to consider the following set in which we look for the central configuration
 \begin{eqnarray}
  0.5 \leqslant x_{n-2} \leqslant (n-1), \label{eq:con0} \\
  -(n-1) \leqslant x_0 <0, \label{eq:con1}\\
  x_0 \leqslant x_i \leqslant x_{n-2}, \quad i=0,\dots,n-1, \label{eq:con2}\\
  y_0 \geqslant 0, \label{eq:con3}\\
  -(n-1) \leqslant y_1 \leqslant  0, \label{eq:con4}\\
  y_1 \leqslant y_i \leqslant (n-1) , \quad i=0,\dots,n-1 \label{eq:con5}\\
  x_2 \leqslant x_3 \leqslant \dots \leqslant x_{n-3} \leqslant x_{n-1}.\label{eq:con6}
 \end{eqnarray}
We call this order {\em increasing} due to the requirement~(\ref{eq:con6}). In the computation we use analogous {\em decreasing} ordering in which we state the opposite, i.e.\
\begin{equation}
  x_2 \geqslant x_3 \geqslant \dots \geqslant x_{n-3} \geqslant x_{n-1}.\label{eq:con7}
\end{equation}
For now, we do not know why it is better to use the decreasing order, but the difference is significant (see Table~\ref{tab:diffOrder}).
\begin{table}[htb]
  \centering
\begin{tabular}{c|r|r}
  $n$ & increasing  & decreasing \\
  \hline
  4 & 0.027170 & 0.026153 \\
  5 & 3.369141 & 2.615399 \\
  6 & 1103.138988 & 924.083085 \\
\end{tabular}
  \caption{Times of asynchronous computations in minutes for different orderings (the computations  were carried out on the computer
Intel® Core™ i7-5500U CPU @ 2.40GHz x 4 with 8GB RAM;
a single thread was used).}\label{tab:diffOrder}
\end{table}

 \subsection{Outline of the approach}
 In the algorithm we look for all zeros of the reduced system (\ref{eq:cc-red-final-i},\ref{eq:cc-red-final-x2}), which under assumption $x_{n-1}\neq x_{n-2}$ by  Theorem~\ref{thm:red-to-full} is equivalent to (\ref{eq:cc-kart}).  For our algorithm proving an existence of locally unique solution in some box is as important as proving that in a given box there is no solution.

 For proving of the existence of the locally unique solution we use the Krawczyk operator applied to the system (\ref{eq:cc-red-final-i},\ref{eq:cc-red-final-x2}). To rule out the  existence of solution we use the exclusion tests discussed in Section~\ref{sec:excl-tests} and also the Krawczyk operator.

 The symmetry of CCs is established by proving the uniqueness in a suitable symmetric box (see~Section~\ref{sec:testing-stage}).

 \subsection{The Krawczyk operator}
\label{subsec:Krawczyk}

The Krawczyk operator\cite{A,K,N} is an interval analysis tool to establish the existence of unique zero for the system of $n$ nonlinear equations in $n$ variables. Below we briefly explain how the Krawczyk operator is derived, as it appears mysterious and little known outside the interval arithmetic community.

\subsubsection{Motivation, heuristic derivation}

Let $F: \mathbb{R}^n \to \mathbb{R}^n$ be a $C^1$-function. We
would like to solve the equation
\begin{equation}
  F(x)=0 \label{eq:F=0}
\end{equation}
We begin by explaining the basic idea of the Krawczyk method. The
Newton method is given by
\begin{equation}
  N(x)=x - dF(x)^{-1}F(x). \label{eq:N}
\end{equation}
It is well known that if $F(x^*)=0$ and $dF(x^*)$ is nonsingular,
then $x^*$ is an attracting fixed point for $N(x)$.
It turns out that the same is true if we replace $dF(x)^{-1}$ by a
fixed matrix $C$, which is sufficiently close to $dF(x^*)^{-1}$.
The modified Newton operator is given by
\begin{equation}
  N_m (x)= x - C F(x). \label{eq:Nm}
\end{equation}
Now let us turn the things around and ask how can we use $N_m$ as
a way to prove the existence of solution of (\ref{eq:F=0}).
This is quite obvious. Namely, if $U$ is homeomorphic to a closed
finite-dimensional ball and if
\begin{equation}
  N_m(U) \subset U,
\end{equation}
then from the Brouwer Theorem it follows, that there exists $x_0
\in U$ such that $N_m(x_0)=x_0$. Since $C$ is invertible we obtain
that $F(x_0)=0$.  To obtain the uniqueness it is enough show that
$N_m$ is a contraction on $U$.
Observe that  it is impossible to verify in a  single interval
evaluation of the formula (\ref{eq:Nm}), that for some interval
set $[x]$ holds $N_m([x]) \subset [x]$, because the computed
diameter of $[x] - CF[x]$ is greater than or equal to $\diam([x])
+ \diam(CF([x]))$.
It turns out the mean value form of $N_m$ can cure this
deficiency. If $x_0 \in [x]$, then
\begin{eqnarray*}
  N_m([x]) \subset N_m(x_0) + [d N_m([x])]_I \cdot ([x] -  x_0) = \\
   x_0 - C F(x_0) + (Id - C [df([X])]_I)([x] - x_0)=K(x_0,[x],F).
\end{eqnarray*}
This explains why the requirement $K(x_0,[x],F) \subset [x]$ has
something to do with zeros of $F(x)$.

\subsubsection{The Krawczyk method}
\label{subsubsec:Krawczyk}
A method proposed by Krawczyk for finding zero's of $F$:
\begin{itemize}
\item  $[x] \subset \mathbb{R}^n$ be an interval set (i.e.
product of intervals),
\item $x_0 \in [x]$. Typically $x_0$ is chosen to be midpoint of $[x]$, we will denote this by $x_0=mid([x])$.
\item $C \in \mathbb{R}^{n \times n}$ be a linear isomorphism
\end{itemize}
The Krawczyk operator is given by
\begin{equation}
  K(x_0,[x],F):= x_0 - C F(x_0) + (Id - C \left[ dF([x])
  \right]_I)([x] - x_0).
\end{equation}
\begin{theorem}
\label{thm:Kr-met}
\begin{description}
\item[1.]If $x^* \in [x]$ and $F(x^*) =0$, then $x^* \in
K(x_0,[x],F) $.

\item[2.] If $K(x_0,[x],F) \subset \inte [x]$, then there  exists in
$[x]$ exactly one solution of equation $F(x)=0$.

\item[3.] If $K(x_0,[x],F) \cap  [x]=\emptyset$ that for all $x \in [x]$ $F(x) \neq 0$.
\end{description}
\end{theorem}

Observe that point 2. in the above theorem gives us the way to establish the existence of unique zero of $F$ in $[x]$, while point 3. rules out the existence of zero in $[x]$ i.e. in the terminology of previous section this is the exclusion test.  When $[x]$ is close to a zero of $F$ then $<F([x])>$ the evaluation of $F([x])$ in the interval arithmetic might produce such overestimates that $0 \in <F([x])>$, while the Krawczyk operator will rule out the existence of $0$ of $F$ in $[x]$. This is in fact quite common phenomenon.

The Krawczyk operator is used as a part of iteration process
\begin{description}
  \item[0.] given $[x]_0 \subset \mathbb{R}^n$
  \item[1.] compute $[y] = K(mid([x]_k),[x]_k,F)$
  \item[2.] if $[y] \subset \inte [x]_k$, then return \textbf{success, a unique solution in $[x]_0$ was found}\\
           elseif $[x]_k \subset [y]$, then return \textbf{failure} \\
           elseif set $[x]_{k+1}:=[y] \cap [x]_k$ and goto \textbf{1.}
\end{description}
The above loop can be executed several times and even if no success is obtained the last computed $[y]$ may give us useful information, because
we know from Theorem~\ref{thm:Kr-met} that all possible zeros are contained in $[y]$. This set instead of $[x]_0$ can be used in further
computations and while $[x]_0 \setminus [y]$ can be discarded. In the next section we will discuss what is essentially a binary search algorithm, which
scales poorly with the number of bodies due to \emph{the curse of dimensionality} \cite{TWW} and the reduction obtained by the Krawczyk method, i.e. replacing $[x]_0$ by $[y]$ for further subdivision  leads to significant speed improvements, because the Krawczyk method on sufficiently small scales appear to work in time polynomial with respect to the number of bodies.

In our context the only weakness of the Krawczyk operator is that it requires the sets of the diameter in each coordinate directions to be smaller than $10^{-2}$ to give us something. Above that threshold we usually have $[x] \subset K(mid([x]),[x],F)$ and the Krawczyk method is useless.

\section{The algorithm}
\label{sec:alg}
The algorithm runs in the reduced configuration space which is a subset of $\reals^{2(n-1)-1}$, i.e.\ a configuration is represented by a point $(x_0, y_0,\ldots, x_{n-3}, y_{n-3}, x_{n-2})$. Physically, we interpret such a configuration as $n-1$ bodies with $q_i=(x_i, y_i)$ for $i = 0, \ldots, n-3$ and $q_{n-2}=(x_{n-2},0)$. From (\ref{eq:cc-red-final-com}) we obtain $q_{n-1}$ the position of the last body.

\noindent
\textit{Input:}
The input of the algorithm consists of
\begin{enumerate}
\item $n$ -- the number of bodies
\item some cube in the reduced configuration space.
\end{enumerate}

\noindent
\textit{Output:}\
All different (up to reflections and rotations) central configurations in the full system for a given input cube. Since we use interval arithmetic for calculations, central configurations are also cubes containing the exact solution in their interior.

Program is divided into two stages: \textit{searching} finds solutions and \textit{testing} tests them to distinguish different CC and to find the kind of symmetry (if any exists).

\subsection{Searching stage}
In this stage we cover the configuration space with cubes. To fulfill requirements of Krawczyk's method (see Theorem~\ref{thm:Kr-met})  we must ensure that every point is in the interior of some cube.
The algorithm runs for any initial cube, however if our goal is to find all the central configurations (for fixed $n$ and equal masses) the reasonable cube is as follows (compare Sec.~\ref{subsec:cap-equal-masses}):
\begin{eqnarray*}
x_0 & \in & [-(n-1), 0]\\
y_0 & \in & [0, n-1]\\
x_1, \ldots, x_{n-3} & \in & [-(n-1), n-1]\\
y_1 & \in & [-(n-1), 0]\\
y_2, \ldots, y_{n-3} & \in & [-(n-1), n-1]\\
x_{n-2} & \in & [0.5, n-1]
\end{eqnarray*}
Simple recursive algorithm works as follows:\begin{enumerate}[label=(\Roman*)]
\item if there is no solution in the cube return 0;
\item if there is unique solution in the cube return 1;
\item otherwise
  bisect the longest edge and recursively run the procedure for both parts;
\item return {\em result}.
\end{enumerate}

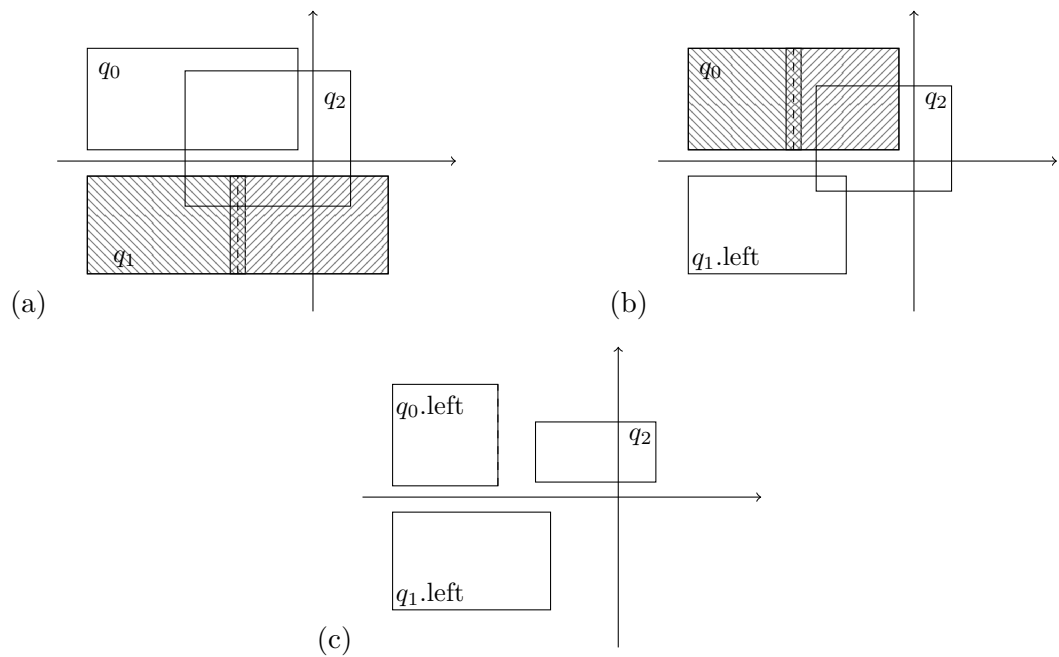
\begin{figure}[htb]
\begin{minipage}[t][25mm][t]{0,49\textwidth}
(a)
\begin{tikzpicture}

\draw[pattern=north west lines, pattern color=gray] (0.5,-1.5) rectangle (2.6,-0.2);

\draw[pattern=north east lines, pattern color=gray] (2.4,-1.5) rectangle (4.5,-0.2);

\draw[draw=black ] (0.5,0.15) rectangle (3.3,1.5);
\node[] at (0.8,1.2) {\small $q_0$};

\draw[draw=black ] (1.8,-0.6) rectangle (4.0,1.2);
\node[] at (3.8,0.8) {\small $q_2$};

\draw[draw=black ] (0.5,-1.5) rectangle (4.5,-0.2);
\node[] at (1.0,-1.3) {\small $q_1$};
\draw[color = black, dashed] (2.5,-1.5) -- (2.5,-0.2);

\draw[color = black, ->] (0.1,0) -- (5.4,0);
\draw[color = black, ->] (3.5,-2.0) -- (3.5,2.0);

\end{tikzpicture}

\end{minipage}
\begin{minipage}[t][25mm][t]{0,49\textwidth}
(b)
\begin{tikzpicture}

\draw[pattern=north west lines, pattern color=gray] (0.5,0.15) rectangle (2.0,1.5);

\draw[pattern=north east lines, pattern color=gray] (1.8,0.15) rectangle (3.3,1.5);

\draw[draw=black ] (0.5,0.15) rectangle (3.3,1.5);
\node[] at (0.8,1.2) {\small $q_0$};
\draw[color = black, dashed] (1.9,0.15) -- (1.9,1.5);

\draw[draw=black ] (2.2,-0.4) rectangle (4.0,1.0);
\node[] at (3.8,0.8) {\small $q_2$};

\draw[draw=black ] (0.5,-1.5) rectangle (2.6,-0.2);
\node[] at (1.0,-1.3) {\small $q_1$.left};

\draw[color = black, ->] (0.1,0) -- (5.4,0);
\draw[color = black, ->] (3.5,-2.0) -- (3.5,2.0);

\end{tikzpicture}

\end{minipage}

\begin{center}
\begin{minipage}[t][25mm][t]{0,49\textwidth}
(c)
\begin{tikzpicture}

\draw[draw=black ] (0.5,0.15) rectangle (1.9,1.5);
\node[] at (1.0,1.2) {\small $q_0$.left};
\draw[color = black, dashed] (1.9,0.15) -- (1.9,1.5);

\draw[draw=black ] (2.4,0.2) rectangle (4.0,1.0);
\node[] at (3.8,0.8) {\small $q_2$};

\draw[draw=black ] (0.5,-1.5) rectangle (2.6,-0.2);
\node[] at (1.0,-1.3) {\small $q_1$.left};

\draw[color = black, ->] (0.1,0) -- (5.4,0);
\draw[color = black, ->] (3.5,-2.0) -- (3.5,2.0);

\end{tikzpicture}

\end{minipage}
\end{center}
\caption{Two steps (a)$\stackrel{1}{\longrightarrow} $ (b)$\stackrel{2}{\longrightarrow} $ (c) of the algorithm on the  cube with each wall projected on the plane: the longest edge is divided and new coordinates of the last body,  $q_2$, are calculated. Notice that after subdivision the new cubes have an intersection with non-empty interior. }\label{fig:idea}
\end{figure}

In the more detailed version of the algorithm the cube is represented by a vector of bodies. The class {\tt Bodies} contains this vector and some methods to manipulate these bodies. An instance of {\tt Bodies} is a cube {\tt bodies} with some additional information.

{\tt
\begin{tabular}{rl}
 1 & {\tt int\ } \textbf{search} (Bodies bodies) \{\\
 2 & \mbox{}\hspace*{5mm} {\tt if\ } (\textbf{thereIsNoSolution}(bodies)) \ {\tt return\ } 0;\\
\\
 3 & \mbox{}\hspace*{5mm} newtonRes\ =\ \textbf{krawczykMethod}(bodies);\\
 4 & \mbox{}\hspace*{5mm} {\tt switch\ } (newtonRes) \{\\
5 & \mbox{}\hspace*{5mm} {\tt case\ } methodFailed:\\
6 & \mbox{}\hspace*{10mm} {\tt break}; \\
\\
 7 & \mbox{}\hspace*{5mm} {\tt case\ } uniqueZero:\\
8 & \mbox{}\hspace*{10mm} {\tt return\ } 1;\\
9 & \mbox{}\hspace*{10mm} {\tt break};\\
\\
10 & \mbox{}\hspace*{5mm} {\tt case\ } noZeroInSet:\\
11 & \mbox{}\hspace*{10mm} {\tt return\ } 0;\\
12 & \mbox{}\hspace*{10mm} {\tt break};\\
\\
13 & \mbox{}\hspace*{5mm} {\tt default}:\\
14 & \mbox{}\hspace*{10mm} {\tt break};\\
15 & \mbox{}\hspace*{8mm} \}\\
\\
16 & \mbox{}\hspace*{5mm} $\slash\slash$ we\ are\ looking\ for\ the\ longest\ interval:\\
17 & \mbox{}\hspace*{5mm} $\slash\slash$  maxiI\ =\ [leftPoint,\ rightPoint]\\
18 & \mbox{}\hspace*{5mm} maxi = bodies.maxDiam();\\
 & \\
19 & \mbox{}\hspace*{5mm}
leftPoint\ =\ bodies.maxi.leftBound();\\
20 & \mbox{}\hspace*{5mm}
rightPoint\ =\ bodies.maxi.rightBound();\\
21 & \mbox{}\hspace*{5mm} midPoint\ =\ (leftPoint\ +\ rightPoint) / 2.0;\\
22 & \mbox{}\hspace*{5mm} margin\ =\ (rightPoint\ -\ leftPoint) *  overlap;\\
\\
23 & \mbox{}\hspace*{5mm} $\slash\slash$ recursively\ call\ search\ for\ the\ left\ half\ of\ maxiI\\
24 & \mbox{}\hspace*{5mm} bodies.set(maxi,\ MyInterval(leftPoint,\ midPoint\ +\ margin));\\
25 & \mbox{}\hspace*{5mm} {\tt int\ } leftCount\ =\ \textbf{search}(bodies,\ resultFile);\\
\\
26 & \mbox{}\hspace*{5mm} $\slash\slash$ recursively\ call\ search\ for\ the\ right\ half\ of\ maxiI\\
27 & \mbox{}\hspace*{5mm} bodies.set(maxi,\ MyInterval(midPoint\ -\ margin,\ rightPoint));\\
28 & \mbox{}\hspace*{5mm} {\tt int\ } rightCount\ =\ \textbf{search}(bodies,\ resultFile);\\
 & \\
29 & \mbox{}\hspace*{5mm} {\tt return\ } leftCount\ +\ rightCount;\\
30 & \mbox{}\hspace*{5mm}\}
\end{tabular}
}\\[1ex]

\subsubsection{Details and optimizations}
\label{subsubsec:det-opty}
The crucial function {\tt thereIsNoSolution(bodies)} contains series of tests (the exclusion tests); if these tests are not satisfied, then we know that there  is no solution in {\tt bodies}:
\begin{enumerate}
\item {\tt checkAprioriBounds(bodies)} --- tests if {\tt bodies} satisfy a priori bounds (see Theorem~\ref{thm:upp-bnd});

\item {\tt checkUEqI(bodies)} --- if there is no collision in {\tt bodies}, tests if $U(q) == I(q)$ (see Lemma~\ref{lem:I=U});

\item {\tt clusterTest(bodies)} --- see Lemma~\ref{lm:FUI} and Subsec.~\ref{subsec:cluster-check-for-zeros}

\item {\tt distanceTest(bodies)} --- tests the order of bodies (see conditions~(\ref{eq:con0}) --- (\ref{eq:con6}))

\item {\tt checkZero(bodies, i)} --- computes functions of vector field (see equations~(\ref{eq:vector-field})) and tests if it is possible to have $F(q_0, \ldots, q_{N-1}) = 0$ as discussed in Sec.~\ref{subsec:check-for-zeros}.
\end{enumerate}
To break (or to rather to relieve) the dimensionality curse we are looking for the possibility of restricting {\tt bodies} before bisecting them (line 16).  First place we are able to do this is the function {\tt thereIsNoSolution(bodies)} (see Section~\ref{subsec:check-for-zeros}). Another place is in {\tt krawczykMethod(bodies)}.
If Krawczyk's method fails (line 5), bodies will have been restricted by intersection with the operator. Hence it is now (from line 17 onwards) the restricted cube that is being processed. This gives a large growth of efficiency.

Since the Krawczyk's  method  is  costly and usually fails for  large sets, we introduced a parameter  {\tt bias}. If  the size (diameter) of all variables is not greater than {\tt bias} then the Krawczyk's method
is run.
There is a big difference in execution time of the program depending on the value of the {\tt bias} parameter; in Table~\ref{tab:bias} we present computation times for 5 bodies . For the same initial  data program finds 8 solutions (some of them are later identified to be the same solutions), but numbers of failed and `no-zero-inside' cubes differs.

\begin{table}[htb]
\begin{center}
\begin{tabular}{l|r|r|r}
\tt bias & time & failed & no zero \\
  & m:s.d & & \\
\hline
1e-4 &  4:27.65 & 0 & 32095 \\
1e-3 &  3:09.13 & 0 & 25846 \\
5e-3 & 2:15.54 & 232 & 32580\\
1e-2 & 1:59.44 & 12886 & 48446\\
1e-1 & 3:29.29 & 585151 &  98946
\end{tabular}
\caption {Comparison of execution times for 5 bodies for different 	thresholds, where we start Krawczyk's method (the computations  were carried out on the computer
Intel® Core™ i7-5500U CPU @ 2.40GHz x 4  with 8GB RAM;
a single thread was used).}\label{tab:bias}
\end{center}
\end{table}

\subsection{Testing stage}\label{sec:testing-stage}
The main goal of this stage is to identify   distinct solutions. Additionally, we check the symmetry of solutions.
In this stage we consider solutions in the full system.

The solutions obtained in the search stage are given as a list of boxes in which we have a unique solution. Some of these boxes may overlap
and can in fact contain the same solution. Because we consider the equal masses case we also do not want to distinguish solutions which differ by the indexes of the bodies.
Hence two solutions produced in the searching stage can in fact be equivalent for two reasons:
\begin{enumerate}[label=(\arabic*)]
\item the only difference is the ordering of the bodies,
\item the boxes defining them have non-empty intersection, having been obtained in different series of partitions.
\end{enumerate}

\begin{figure}[htb]
\centering
\begin{tikzpicture}
\draw[step=1.5cm,gray,very thin] (-0.5,-0.5) grid (6.5,6.5);

\draw[pattern=north west lines, pattern color=gray] (-0.3,2.3) rectangle (0.5, 3.3) node[below = 10mm] {$q_0$};

\draw[pattern=north west lines, pattern color=gray] (1.2,2.7) rectangle (1.9, 3.5) node[above] {$q_2$};

\draw[pattern=north west lines, pattern color=gray] (2.5,2.5) rectangle (3.6, 3.3) node[right] {$q_3$};

\draw[pattern=north west lines, pattern color=gray] (4.2,2.6) rectangle (4.8, 3.6) node[above] {$q_4$};

\draw[pattern=north west lines, pattern color=gray] (5.4,2.2) rectangle (6.2, 3.2) node[below = 10mm] {$q_5$};

\draw[pattern=north west lines, pattern color=gray] (2.2,5.5) rectangle (3.3, 6.3)node[left = 1cm] {$q_6$};

\draw[pattern=north west lines, pattern color=gray] (2.3,-0.4) rectangle (3.2, 0.8) node[left = 1cm] {$q_1$};


\draw[pattern=north east lines, pattern color=gray] (-0.2,2.6) rectangle (0.6, 3.5) node[right, color = gray] {$q_0$};

\draw[pattern=north east lines, pattern color=gray] (1.3,2.5) rectangle (2.2, 3.2) node[below = 7mm, color = gray] {$q_2$};

\draw[pattern=north east lines, pattern color=gray] (2.7,2.2) rectangle (3.3, 3.5) node[above, color = gray] {$q_6$};

\draw[pattern=north east lines, pattern color=gray] (4.3,2.4) rectangle (5.0, 3.2) node[below = 8mm, color = gray] {$q_4$};

\draw[pattern=north east lines, pattern color=gray] (5.6,2.6) rectangle (6.4, 3.6) node[right, color = gray] {$q_5$};

\draw[pattern=north east lines, pattern color=gray] (2.7,5.5) rectangle (3.8, 6.3)node[right, color = gray] {$q_3$};

\draw[pattern=north east lines, pattern color=gray] (2.7,-0.4) rectangle (3.5, 0.8) node[right, color = gray] {$q_1$};


\draw[dashed, thick, color=black] (-0.4,2.3) rectangle (0.6, 3.5);

\draw[dashed, thick, color=black] (1.2,2.5) rectangle (2.2, 3.5);

\draw[dashed, thick, color=black] (2.5,2.2) rectangle (3.6, 3.5);

\draw[dashed, thick, color=black] (4.2,2.4) rectangle (5.0, 3.6);

\draw[dashed, thick, color=black] (5.4,2.2) rectangle (6.4, 3.6);

\draw[dashed, thick, color=black] (2.2,5.5) rectangle (3.8, 6.3);

\draw[dashed, thick, color=black] (2.3,-0.4) rectangle (3.5, 0.8);


\draw[fill] (0,3) circle (0.06);
\draw[fill] (1.5,3) circle (0.06);
\draw[fill] (3,3) circle (0.06);
\draw[fill] (4.5,3) circle (0.06);
\draw[fill] (6,3) circle (0.06);
\draw[fill] (3,6) circle (0.06);
\draw[fill] (3,0) circle (0.06);
\end{tikzpicture}
\caption{Reasons for two solutions (projected onto the plane) to be equivalent: the exact unique central configuration (marked with dots) is inside both. Interval hull of these solutions is marked by dashed boxes. }\label{fig:idi}
\end{figure}
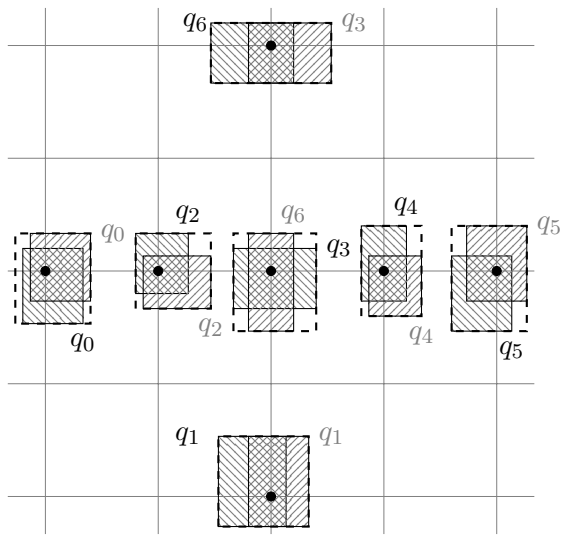

For any pair of solutions we treat the first one as a `model', whilst bodies in the second  solution are permuted in an attempt to match the model. When trying to match the solutions, we create a set containing both of them (an interval hull) and we prove the uniqueness within it. The rough idea is shown below:

{\tt
\begin{tabular}{rl}
1 & bool\ theSameSolution(Bodies\ \&\!\!\! sol1,\ Bodies\ \&\!\!\! sol2) \{\\
2 & \mbox{}\hspace*{5mm} Bodies\ unionBodies; \\
3 & \mbox{}\hspace*{5mm}
    intervalHullBodies(sol1,\ sol2,\ unionBodies);\\
  & \\
4 & \mbox{}\hspace*{5mm} int\ cZeros\ =\ search(unionBodies);\\
5 & \mbox{}\hspace*{5mm} if\ (cZeros\ ==\ 1)\\
6 & \mbox{}\hspace*{10mm} return\ true;\ \slash\slash sol1\ and\ sol2\ are\ the\ same\ solution\\
7 & \mbox{}\hspace*{5mm} else\ \{\\
8 & \mbox{}\hspace*{10mm} cZeros\ =\ blowUp(unionBodies);\\
9 & \mbox{}\hspace*{10mm} return\ (cZeros\ ==\ 1);\\
10 & \mbox{}\hspace*{5mm} \}\\
11 & \}

\end{tabular}
}

It may happen that there exists a solution in the set {\tt unionBodies} but the set is too small to prove this using the Krawczyk operator, thus we `inflate' it and retry the proof in the bigger set (the function {\tt blowUp(unionBodies)}).

Establishing a (reflectional) symmetry of CC is very similar to testing the uniqueness of solutions.
However, it requires an additional step to calculate a symmetric image. We take an interval hull of CC and its symmetric image. The possible lines of reflection are an axis $OX$ or  the bisector of the angle with the vertex at $(0,0)$ and the  rays passing through  $q_{n-2}$ (the body furthest from $(0,0)$) and through $q_i$ (different bodies are tested). The sketch of the function establishing this symmetry is given below.

{\tt
\begin{tabular}{rl}
1 & bool\ checkSym(Bodies \&\!\!\! bodies) \{\\
2 & \mbox{}\hspace*{5mm}    calculateBisectorOfAngle(cosB,\ sinB);\\
 & \mbox{}\hspace*{3cm} \slash\slash the\ line\ is\ t(cosB,\ sinB)\\
3 & \mbox{}\hspace*{5mm} Bodies\ symBodies,\ unionBodies;\\
4 & \mbox{}\hspace*{5mm} calculateSymmetricImageOfTheSolution(bodies,\  cosB,\  sinB,\ symBodies));\\
5 & \mbox{}\hspace*{5mm}
    intervalHullBodies(bodies,\ symBodies,\ unionBodies);\\
  & \\
6 & \mbox{}\hspace*{5mm} int\ cZeros\ =\ search(unionBodies);\\
7 & \mbox{}\hspace*{5mm} if\ (cZeros\ ==\ 1)\\
8 & \mbox{}\hspace*{10mm} return\ true;\ \slash\slash bodies\ and\ symBodies\ are\ symmetric\\
9 & \mbox{}\hspace*{5mm} else\ \{\\
10 & \mbox{}\hspace*{10mm} cZeros\ =\ blowUp(unionBodies);\\
11 & \mbox{}\hspace*{10mm} return\ (cZeros\ ==\ 1);\\
12 & \}
\end{tabular}
}

In line~2, we calculate the parameters of the possible reflection symmetry line, but the symmetry tested contains also a permutation of bodies, we construct a configuration {\tt symBodies} considering all possible permutations of bodies.
Note that lines 5--12 in the function {\tt checkSym} are identical, up to the variable names, to lines 3--11 in {\tt theSameSolutions}.

\subsection{Technical data}
\label{subsec:tech-data}
The main computations  were
carried out in parallel using the template function {\tt std::async} (from the standard C++ library) which runs the function  asynchronously (potentially in a separate thread which may be part of a thread pool) on Dell R930 4x Intel Xeon E7-8867 v3 (2,5GHz, 45MB), 1024 GB RAM. The compiler  is gcc version 4.9.2 (Debian 4.9.2-10+deb8u2).
The best obtained times for different number of bodies with ${\tt bias}=10^{-2}$ are presented in Table~\ref{tab:cc}.

\begin{table}[h]
\begin{center}
\begin{tabular}{r|r|r| r | r}
no bodies & no CCs & total no of & elapsed time & average percentage \\
 & & CPU-seconds & h:m:s.d & of the CPU\\
\hline
3 & 2 &  0.05 & 0:00.06 & 240\\
4 & 4 &  3.07 & 0:00.69 & 666\\
5 & 5 &   203.82 & 0:24.84 & 1023\\
6 & 9 & 42430.04 & 59:51.59 & 1203\\
7 & 14 & 8490959.77 & 98:56:00.00 & 2531
\end{tabular}
\caption {Comparison of execution times for different number of bodies.}\label{tab:cc}
\end{center}
\end{table}

\section{Minimizing dependency problem in gravitational force evaluation}\label{sec:dependency}

In this section we describe a method of calculations of functions and their derivatives used in the program.
Let us denote $F_i = (f_i^1, f_i^2)$, $r_{ij} = \sqrt{(x_i - x_j)^2 + (y_i - y_j)^2}$. Then functions $f_i^{[1,2]}$ and their derivatives are (with analogs for $y$'s):
\begin{eqnarray}
\displaystyle{f_i^1} & = & \displaystyle{x_i - \sum_{\substack{j = 1,\\
j\neq i}}^{N} m_j \frac{x_i - x_j}{r_{ij}^{3}}}\nonumber\\
\nonumber\\
\displaystyle{f_i^2} & = & \displaystyle{y_i - \sum_{\substack{j = 1,\\
j\neq i}}^{N} m_j \frac{y_i - y_j}{r_{ij}^3}}\nonumber\\
\nonumber \\
\displaystyle{\frac{\partial}{\partial x_k}f_i^1} & = & \displaystyle{-m_k \left(3\frac{(x_i - x_k)^2}{r_{ik}^5} - \frac{1}{r_{ik}^3}\right) + m_k\left(3\frac{(x_i - x_N)^2}{r_{iN}^5} - \frac{1}{r_{iN}^3}\right)} \quad \mbox{for $k\neq i$}\nonumber\\
\nonumber \\
\nonumber \\
\displaystyle{\frac{\partial}{\partial x_i}f_i^1} & = & \displaystyle{1 + \sum_{\substack{j = 1,\\
j\neq i}}^{N-1} m_j \left(3\frac{(x_i - x_j)^2}{r_{ij}^5} - \frac{1}{r_{ij}^3}\right) + m_N\left(1 + \frac{m_i}{m_N}\right)\left(3\frac{(x_i - x_N)^2}{r_{iN}^5} - \frac{1}{r_{iN}^3}\right)} \nonumber\\
\nonumber \\
\nonumber \\
 \displaystyle{\frac{\partial}{\partial x_k}f_i^2} & = & \displaystyle{3m_k \frac{(x_i - x_k)(y_i - y_k)}{r_{ik}^5} + 3m_k\frac{(x_i - x_N)(y_i - y_N)}{r_{iN}^5}}\label{eqn:df1i-dxk} \\
\nonumber \\
\nonumber \\
 \displaystyle{\frac{\partial}{\partial x_i}f_i^2} & = & \displaystyle{ 3\sum_{\substack{j = 1,\\
j\neq i}}^{N-1} m_j \frac{(x_i - x_j)(y_i - y_j)}{r_{ij}^5} + 3m_N\left(1+\frac{m_i}{m_N}\right)\frac{(x_i - x_N)(y_i - y_N)}{r_{iN}^5}} \label{eqn:df2i-dxi}\\
\displaystyle{\frac{\partial}{\partial y_k}f_i^1} & = & \displaystyle{\frac{\partial}{\partial x_k}f_i^2}\quad \mbox{for all $k$.}\nonumber
\end{eqnarray}
In the program we have to evaluate $f_i^{[1,2]}$ on a box $D$ in configurations space. The naive interval evaluation of $f_i^{[1,2]}$, where we just plug-in
the interval arguments might lead to severe overestimation due to the dependency problem \cite{Mo,N}. The best solution would be
a cheap but rigorous estimate of $\sup$ and $\inf$ of $f_{i}^1$  and $f_{i}^2$ over $D$; this however appears to be a difficult and costly task.
For the Krawczyk method we also need good estimates for $\frac{\partial f_i}{\partial x_j}$ and  $\frac{\partial f_i}{\partial y_j}$ and we face the same problem.
Thus we decided to optimize the computation of the following expressions
$$
\frac{(x_j - x_i)^a}{r_{ij}^b}\quad\mbox{and}\quad\frac{(y_j - y_i)^a}{r_{ij}^b}
$$
over $q = (x,y) = (x_j - x_i, y_j - y_i)\in D \subset \mathbb{R}^2$, since such components appear in all above equations. In equations~(\ref{eqn:df1i-dxk}) and~(\ref{eqn:df2i-dxi}) we evaluate expressions in the form $xy/r^5$ as a product of $x/r^3$ and  $y/r^2$ which are treated as separate expressions.

\subsection{Estimates for \texorpdfstring{$x^a / r^b$}{Lg} and \texorpdfstring{$y^a / r^b$}{Lg}}
Assume we want to calculate
\begin{eqnarray}
f_x(x,y) & = & \frac{x^a}{r^b} = \frac{x^a}{(x^2 + y^2)^\frac{b}{2}}\\
f_y(x,y) & = & \frac{y^a}{r^b} = \frac{y^a}{(x^2 + y^2)^\frac{b}{2}},
\end{eqnarray}
where $a < b$, $(x,y) = (x_j - x_i, y_j - y_i)$. Let us define $D = (x_L, x_R)\times (y_L,y_R)$. We want to estimate $f_x$ and $f_y$ on $D$. We always assume that $(0,0)\not\in D$.
To minimize the overestimation of these calculations we look for the possible local extrema in $D$.

\begin{figure}[htb]
  \centering
\begin{tikzpicture}

\draw[draw=black ] (1.8,-0.4) rectangle (4.2,1.0);

\draw[dashed] (1.8, -1.0) -- (1.8, -0.4);
\draw[dashed] (4.2, -1.0) -- (4.2, -0.4);
\node[] at (2.0,-1.2) {\small $x_i$};
\node[] at (4.15,-1.2) {\small $x_j$};
\node[] at (1.55,-0.4) {\small $y_i$};
\node[] at (1.55,1.0) {\small $y_j$};

\draw[color = black, ->] (0.1,-1.0) -- (5.4,-1.0);
\draw[color = black, ->] (2.8,-2.0) -- (2.8,2.0);
\node[] at (5.5, -1.2) {\small $x$};
\node[] at (3.0,2.0) {\small $y$};

\draw[color = blue] (2.4,-2.0) -- (3.8,1.5);
\draw[color = blue] (1.8,1.5) -- (3.2,-2.0);
\draw[color = blue] (1.9,-2.0) -- (4.6,1.0);
\draw[color = blue] (3.7,-2.0) -- (1.0,1.0);

\fill[black] (2.8,-0.4) circle (0.5mm);
\fill[black] (2.8,1.0) circle (0.5mm);

\fill[blue] (3.6,1.0) circle (0.5mm);
\fill[blue] (3.05,-0.4) circle (0.5mm);

\fill[blue] (2.0,1.0) circle (0.5mm);
\fill[blue] (2.25,-0.4) circle (0.5mm);

\fill[blue] (1.8,0.13) circle (0.5mm);
\fill[blue] (2.55,-0.4) circle (0.5mm);

\fill[blue] (4.2,0.56) circle (0.5mm);
\fill[blue] (3.35,-0.4) circle (0.5mm);

\fill[red] (4.2,1.0) circle (0.5mm);
\fill[red] (4.2,-0.4) circle (0.5mm);
\fill[red] (1.8,1.0) circle (0.5mm);
\fill[red] (1.8,-0.4) circle (0.5mm);

\end{tikzpicture}
\caption{Calculated lines and location of critical points.}\label{fig:cute-b}

\end{figure}
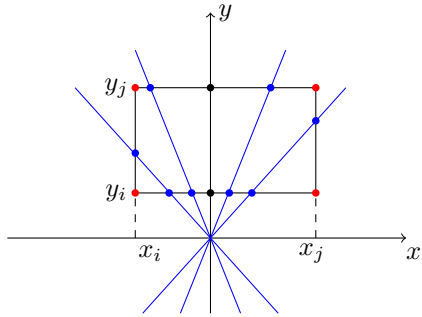

Let us consider the function $f_x$; the second case of $f_y$ is symmetrical. First notice that by solving the system of equations
\begin{eqnarray}
\frac{\partial}{\partial x} f_x(x,y) & = & \frac{x^{a-1}(x^2(a-b) + ay^2)}{r^{b+2}} = 0\label{eqn:diff-x}\\
\frac{\partial}{\partial y} f_x(x,y) & = & - \frac{bx^2y}{r^{b+2}} = 0\label{eqn:diff-y}
\end{eqnarray}
we obtain $(x,y) = (0,0)$, which is impossible in our settings. Since there is no local extremum inside $D$ thus it is attained on the edge of $D$.
Te determine the specific points, where this extremum can be, we explore the equations~(\ref{eqn:diff-x}) and (\ref{eqn:diff-y}) and analogical for $f_y$, and we obtain:
\begin{itemize}
\item border points on the lines $y = \pm x\sqrt{\frac{b-a}{a}}$ and $x = \pm y\sqrt{\frac{b-a}{a}}$ (blue points in Figure~\ref{fig:cute-b})
\item border points for $x = 0$ and for $y = 0$ (black points in Figure~\ref{fig:cute-b})
\end{itemize}
Additionally we consider corners of $D$ (red points in Figure~\ref{fig:cute-b}).
Next we examine all these point to establish the maximum and the minimum value of $f_x$ and $f_y$.

Note, that there is still a lot of room for further optimization, but for now only this version is implemented in the program.

\appendix
\newpage
\section{Results and conclusions from the run of the program}\label{sec:proofs}
Below we give the output (the report files) from the runs of our program  for various number of bodies. These are proofs of Theorems~\ref{thm:main} and~\ref{thm:non-sym}. In the computations the masses were normalized so that $\sum m_i=1$
and we set $\lambda=1$ in (\ref{eq:cc-with-lambda}).

Other authors  have used different normalizations,  both in \cite{MNum} and \cite{F02}
all masses were equal to $1$,  and in \cite{MNum} the condition $I=1$ was demanded, while in \cite{F02} the following normalization has been
used $|q_1| \leqslant |q_2| \leqslant  \dots \leqslant  |q_{n-1}|=1 \leqslant  |q_n|$, with $q_{n-1}=(1,0)$.

It is easy to see that the following quantity
\begin{equation}
J=\frac{U I^{1/2}}{(\sum_i m_i)^{5/2}},  \label{eq:inv}
\end{equation}
is an invariant of the scaling transformations of configuration space and masses. We believe that this is a meaningful invariant which makes comparison between
CC obtained in different works relatively easy, as it gives an obvious formula for scaling  of CCs normalized in various ways.

In order to make a comparison with \cite{MNum,F02} easier we also include
the Moeckel's potential, which equal to $U(q)$ where $q$ is CC scaled so that it satisfies $I=1$ and $m_i=1$. This Moeckel's potential
coincides with $U \sqrt{I}$ contained in table for CC from \cite{F02}.
To compare coordinate-wise the CC given  in \cite{F02,MNum} with ours it would be necessary to scale, rotate and permute the bodies (indices).
We do not do it here.
Our work  have proved  that the lists of CC given in  \cite{F02,MNum} for $n=4,5,6,7$ are complete, by which we mean that we were able to identify
all our CCs with the ones provided in these works and vice-versa.

For each CC we tried to prove that
\begin{itemize}
  \item it has a reflectional symmetry with respect $OX$-axis
  \item it has a reflectional symmetry with respect to some other line. This other line is obtained as the bisector of the angle with the vertex at $(0,0)$ and the  rays passing through  $q_{n-2}$ (the body furthest from $(0,0)$) and through   $q_i$ (different bodies are tested).
\end{itemize}
Each of these reflections act in $\mathbb{R}^2$ and are accompanied with a suitable permutation $\sigma$ of bodies, which is displayed as a list
of $\sigma(0),\dots,\sigma(n-1)$. For example, for $n=3$ the list $2,1,0$ means that $\sigma(0)=2$, $\sigma(1)=1$ and $\sigma(2)=0$.

The meaning of parameters in report files:
\begin{itemize}
  \item {\tt eps} $= \varepsilon$, if the diameter of the box in max norm is smaller than this number then program stops subdividing this box and returns
     the message that it can not conclude whether there is CC in this box. This never happened in our computations
  \item {\tt bias} -- the Krawczyk method is applied to a box only if the diameter of each side of  the box is smaller than bias, as discussed in Sec.~\ref{subsubsec:det-opty}
\end{itemize}

\subsection{Report files}

\subsubsection*{Three bodies}
\begin{verbatim}
Number of bodies = 3
Accuracy eps = 1e-05, bias = 0.01

Input data:
i: 0 X: [-2.001, 2.01] Y: [-2.001, 2.01] mass: [0.3333333333, 0.3333333333]
i: 1 X: [-0.001, 2.01] Y: [0, 0] mass: [0.3333333333, 0.3333333333]
i: 2 X: [-4.02, 2.002] Y: [-2.01, 2.001] mass: [0.3333333333, 0.3333333333]


The number of undecided cubes: 0
The number of zeros in the method: 3
The number of calls the main search function: 467

Program computed 0.00011735 minutes

Tests usage:
checkAprioriBounds -- 0
clusterTest -- 6
distanceTest -- 55
checkZero -- 114
krawczyk: methodFailed  -- 0
krawczyk: zeroIside -- 3
krawczyk: no zero inside -- 14
spreadTest -- 0
U = I -- 42

Different cc:
---------------------
position 0
i: 0 X: [-0.7469007911, -0.7469007911] Y: [-1.141274668e-14, 1.141055616e-14]
i: 1 X: [0.7469007911, 0.7469007911] Y: [0, 0]
i: 2 X: [-1.440514374e-14, 1.432187702e-14] Y: [-1.141055616e-14, 1.141274668e-14]

U = [0.3719071945, 0.3719071945], I = [0.3719071945, 0.3719071945],
U*(I)^(1/2)/(M)^(5/2) = [0.2268046058, 0.2268046058]
Moeckel's potential = [3.535533906, 3.535533906]

collinear solution no 1
permutation: 0, 1, 2,
symmetric with respect to OX no 1
permutation: 1, 0, 2,
reflectional symmetry with respect to other line 1
---------------------
position 1
i: 0 X: [-0.2886751346, -0.2886751346] Y: [-0.5, -0.5]
i: 1 X: [0.5773502692, 0.5773502692] Y: [0, 0]
i: 2 X: [-0.2886751346, -0.2886751346] Y: [0.5, 0.5]

U = [0.3333333333, 0.3333333333], I = [0.3333333333, 0.3333333333],
U*(I)^(1/2)/(M)^(5/2) = [0.1924500897, 0.1924500897]
Moeckel's potential = [3, 3]

permutation: 2, 1, 0,
symmetric with respect to OX no 2
permutation: 1, 0, 2,
reflectional symmetry with respect to other line 2

Number of different cc = 2
\end{verbatim}
\hrule

\subsubsection*{Four bodies}
\begin{verbatim}
Number of bodies = 4
Accuracy eps = 1e-05, bias = 0.01

Input data:
i: 0 X: [-3.001, 0.01] Y: [-0.001, 3.01] mass: [0.25, 0.25]
i: 1 X: [-3.001, 3.01] Y: [-3.001, 0.01] mass: [0.25, 0.25]
i: 2 X: [0.499, 3.01] Y: [0, 0] mass: [0.25, 0.25]
i: 3 X: [-6.03, 5.503] Y: [-3.02, 3.002] mass: [0.25, 0.25]


The number of undecided cubes: 0
The number of zeros in the method: 5
The number of calls of the main search function: 39485

Program computed 0.02834338333 minutes

Tests usage:
checkAprioriBounds -- 0
clusterTest -- 377
distanceTest -- 3506
checkZero -- 11564
krawczyk: methodFailed  -- 5381
krawczyk: zeroIside -- 5
krawczyk: no zero inside -- 3310
spreadTest -- 0
U = I -- 981

Different CC:
---------------------
position 0
i: 0 X: [-0.9051285388, -0.9051285354] Y: [-5.049805694e-09, 5.049805696e-09]
i: 1 X: [-0.2862410122, -0.2862410082] Y: [-1.476723831e-09, 1.476723833e-09]
i: 2 X: [0.9051285343, 0.9051285399] Y: [0, 0]
i: 3 X: [0.2862410037, 0.2862410167] Y: [-6.526529529e-09, 6.526529526e-09]

U = [0.4505957879, 0.4505957967], I = [0.450595789, 0.4505957955],
U*(I)^(1/2)/(M)^(5/2) = [0.3024688757, 0.3024688838]
Moeckel's potential = [9.679004022, 9.679004282]

collinear solution no 1
permutation: 0, 1, 2, 3,
symmetric with respect to OX no 1
permutation: 2, 3, 0, 1,
reflectional symmetry with respect to other line 1
---------------------
position 1
i: 0 X: [-0.6208313565, -0.6208117917] Y: [-1.338276202e-05, 1.338276202e-05]
i: 1 X: [-9.896654491e-06, 9.896654492e-06] Y: [-0.620827361, -0.6208157872]
i: 2 X: [0.6208135881, 0.6208295602] Y: [0, 0]
i: 3 X: [-2.766510061e-05, 2.766510061e-05] Y: [0.6208024044, 0.6208407438]

U = [0.3853953437, 0.3854435136], I = [0.3853938373, 0.3854450185],
U*(I)^(1/2)/(M)^(5/2) = [0.239253801, 0.2392995931]
Moeckel's potential = [7.656121632, 7.65758698]

permutation: 0, 3, 2, 1,
symmetric with respect to OX no 2
permutation: 3, 2, 1, 0,
reflectional symmetry with respect to other line 2
---------------------
position 2
i: 0 X: [-0.3821936947, -0.3821936947] Y: [0.6195346528, 0.6195346528]
i: 1 X: [-0.3821936947, -0.3821936947] Y: [-0.6195346528, -0.6195346528]
i: 2 X: [0.7436490828, 0.7436490828] Y: [0, 0]
i: 3 X: [0.0207383067, 0.02073830672] Y: [-5.828892924e-12, 5.828559857e-12]

U = [0.4033086121, 0.4033086121], I = [0.4033086121, 0.4033086121],
U*(I)^(1/2)/(M)^(5/2) = [0.2561275196, 0.2561275196]
Moeckel's potential = [8.196080629, 8.196080629]

permutation: 1, 0, 2, 3,
symmetric with respect to OX no 3
---------------------
position 3
i: 0 X: [-0.3666565002, -0.3666565001] Y: [0.6350676872, 0.6350676872]
i: 1 X: [-0.3666565002, -0.3666565001] Y: [-0.6350676872, -0.6350676872]
i: 2 X: [0.7333130003, 0.7333130003] Y: [0, 0]
i: 3 X: [-2.157363177e-11, 2.157307666e-11] Y: [-6.115996598e-12, 6.114997397e-12]

U = [0.4033001027, 0.4033218324], I = [0.4033037569, 0.403318178],
U*(I)^(1/2)/(M)^(5/2) = [0.2561205739, 0.256138953]
Moeckel's potential = [8.195858366, 8.196446497]

permutation: 1, 0, 2, 3,
symmetric with respect to OX no 4
permutation: 2, 1, 0, 3,
reflectional symmetry with respect to other line 3

Number of different cc = 4
\end{verbatim}

\subsubsection*{Five bodies}
\begin{verbatim}
Number of bodies = 5
Accuracy eps = 1e-05, bias = 0.01

Input data:
i: 0 X: [-4.001, 0.001] Y: [-0.001, 4.001] mass: [0.2, 0.2]
i: 1 X: [-4.001, 4.001] Y: [-4.001, 0.001] mass: [0.2, 0.2]
i: 2 X: [-4.001, 4.001] Y: [-4.001, 4.001] mass: [0.2, 0.2]
i: 3 X: [0.499, 4.001] Y: [0, 0] mass: [0.2, 0.2]
i: 4 X: [-12.004, 11.504] Y: [-8.003, 8.003] mass: [0.2, 0.2]


The number of undecided cubes: 0
The number of zeros in the method: 25
The number of calls of the main search function: 2722255

Program computed 1.946497467 minutes

Tests usage:
checkAprioriBounds -- 17
clusterTest -- 45125
distanceTest -- 279903
checkZero -- 931717
krawczyk: methodFailed  -- 13549
krawczyk: zeroIside -- 8
krawczyk: no zero inside -- 56614
U = I -- 47744

Different CC:
---------------------
position 0
i: 0 X: [-0.7315026092, -0.7314991044] Y: [-1.22322208e-05, 1.22322208e-05]
i: 1 X: [-4.380645756e-06, 4.380645756e-06] Y: [-0.7315035144, -0.7314981992]
i: 2 X: [-3.486361374e-06, 3.486361374e-06] Y: [-6.548894919e-06, 6.54889492e-06]
i: 3 X: [0.7314993242, 0.7315023894] Y: [0, 0]
i: 4 X: [-1.115205008e-05, 1.115205008e-05] Y: [0.7314794181, 0.7315222955]

U = [0.4280678291, 0.4280817768], I = [0.4280667911, 0.4280828147],
U*(I)^(1/2)/(M)^(5/2) = [0.2800711397, 0.2800855073]
Moeckel's potential = [15.65645267, 15.65725584]

permutation: 0, 4, 2, 3, 1,
symmetric with respect to OX no 1
permutation: 4, 3, 2, 1, 0,
reflectional symmetry with respect to other line 1
---------------------
position 2
i: 0 X: [-1.019255982, -1.019255982] Y: [-4.721349739e-15, 4.724568794e-15]
i: 1 X: [-0.480767439, -0.480767439] Y: [-2.083289376e-15, 2.084403198e-15]
i: 2 X: [0.480767439, 0.480767439] Y: [-4.095121196e-15, 4.092555012e-15]
i: 3 X: [1.019255982, 1.019255982] Y: [0, 0]
i: 4 X: [-9.714451465e-15, 9.575673587e-15] Y: [-1.0901527e-14, 1.089976031e-14]

U = [0.5080080345, 0.5080080345], I = [0.5080080345, 0.5080080345],
U*(I)^(1/2)/(M)^(5/2) = [0.3620811129, 0.3620811129]
Moeckel's potential = [20.24094955, 20.24094955]

collinear solution no 1
permutation: 0, 1, 2, 3, 4,
symmetric with respect to OX no 2
permutation: 3, 2, 1, 0, 4,
reflectional symmetry with respect to other line 2
---------------------
position 4
i: 0 X: [-0.6591405438, -0.659140481] Y: [0.1800138879, 0.1800139591]
i: 1 X: [-0.2807232123, -0.2807230535] Y: [-0.7143993153, -0.7143992478]
i: 2 X: [0.09876565526, 0.09876576604] Y: [-0.1449273191, -0.1449272214]
i: 3 X: [0.767575262, 0.7675752984] Y: [0, 0]
i: 4 X: [0.07352247002, 0.07352283887] Y: [0.6793125102, 0.6793127465]

U = [0.4285688598, 0.4285690028], I = [0.4285688537, 0.4285690089],
U*(I)^(1/2)/(M)^(5/2) = [0.2805633344, 0.2805634788]
Moeckel's potential = [15.68396719, 15.68397526]

permutation: 0, 1, 2, 3, 4,
there is no symmetry with respect to OX
permutation: 4, 3, 2, 1, 0,
reflectional symmetry with respect to other line 3
---------------------
position 5
i: 0 X: [-0.7609899148, -0.7609897315] Y: [0.4684633687, 0.4684636923]
i: 1 X: [-0.1902473651, -0.1902472194] Y: [-0.6719536689, -0.6719535985]
i: 2 X: [0.3050645568, 0.3050646073] Y: [0.02352925245, 0.02352944105]
i: 3 X: [0.8936238356, 0.8936238701] Y: [0, 0]
i: 4 X: [-0.2474515264, -0.2474511125] Y: [0.1799604652, 0.1799610478]

U = [0.4544156383, 0.4544159448], I = [0.4544156667, 0.4544159163],
U*(I)^(1/2)/(M)^(5/2) = [0.3063232187, 0.3063235095]
Moeckel's potential = [17.1239885, 17.12400476]

permutation: 0, 1, 2, 3, 4,
there is no symmetry with respect to OX
permutation: 3, 1, 4, 0, 2,
reflectional symmetry with respect to other line 4
---------------------
position 7
i: 0 X: [-0.5262827178, -0.5262702462] Y: [0.3823582109, 0.3823662804]
i: 1 X: [0.2010155547, 0.2010239026] Y: [-0.6186815315, -0.6186686875]
i: 2 X: [0.2010162837, 0.2010231736] Y: [0.6186717554, 0.6186784636]
i: 3 X: [0.6505104284, 0.6505165851] Y: [0, 0]
i: 4 X: [-0.526293415, -0.526259549] Y: [-0.3823760565, -0.3823484348]

U = [0.4231607357, 0.4231749094], I = [0.423156383, 0.4231792621],
U*(I)^(1/2)/(M)^(5/2) = [0.2752680534, 0.2752847151]
Moeckel's potential = [15.38795199, 15.3888834]

permutation: 4, 2, 1, 3, 0,
symmetric with respect to OX no 3
permutation: 0, 3, 4, 1, 2,
reflectional symmetry with respect to other line 5

Number of different cc = 5
\end{verbatim}

\subsubsection*{Six bodies}
\begin{verbatim}
Number of bodies = 6
Accuracy eps = 1e-05, bias = 0.01

Input data:
i: 0 X: [-5.001, 0.01] Y: [-0.001, 5.01] mass: [0.1666666667, 0.1666666667]
i: 1 X: [-5.001, 5.01] Y: [-5.001, 5.01] mass: [0.1666666667, 0.1666666667]
i: 2 X: [-5.001, 5.01] Y: [-5.001, 5.01] mass: [0.1666666667, 0.1666666667]
i: 3 X: [-5.001, 5.01] Y: [-5.001, 5.01] mass: [0.1666666667, 0.1666666667]
i: 4 X: [0.499, 5.01] Y: [0, 0] mass: [0.1666666667, 0.1666666667]
i: 5 X: [-20.05, 19.505] Y: [-20.04, 15.004] mass: [0.1666666667, 0.1666666667]


The number of undecided cubes: 0
The number of zeros in the method: 17
The number of calls the main search function: 367150795

Program computed 613.6726331 minutes

Tests usage:
checkAprioriBounds -- 847
clusterTest -- 5966128
distanceTest -- 37701237
checkZero -- 127981911
krawczyk: methodFailed  -- 29404399
krawczyk: zeroIside -- 17
krawczyk: no zero inside -- 9330011
spreadTest -- 0
U = I -- 2595247

Different CC:
---------------------
position 0
i: 0 X: [-0.4566157831, -0.4566157831] Y: [0.7908817358, 0.7908817358]
i: 1 X: [-0.4566157831, -0.4566157831] Y: [-0.7908817358, -0.7908817358]
i: 2 X: [0.3632321511, 0.3632321511] Y: [-7.031018693e-14, 7.042146137e-14]
i: 3 X: [-0.1816160756, -0.1816160756] Y: [-0.3145682703, -0.3145682703]
i: 4 X: [0.9132315661, 0.9132315661] Y: [0, 0]
i: 5 X: [-0.1816160756, -0.1816160756] Y: [0.3145682703, 0.3145682703]

U = [0.4829647445, 0.4829647445], I = [0.4829647445, 0.4829647445],
U*(I)^(1/2)/(M)^(5/2) = [0.3356395606, 0.3356395606]
Moeckel's potential = [29.59724379, 29.59724379]

permutation: 1, 0, 2, 5, 4, 3,
symmetric with respect to OX no 1
permutation: 4, 1, 5, 3, 0, 2,
reflectional symmetry with respect to other line 1
---------------------
position 1
i: 0 X: [-0.8729847413, -0.8729699697] Y: [-3.153053147e-05, 3.153053147e-05]
i: 1 X: [-2.066755246e-05, 2.066755246e-05] Y: [-0.7348115058, -0.7347940916]
i: 2 X: [0.2958793068, 0.2958952513] Y: [-2.213796528e-05, 2.213796528e-05]
i: 3 X: [-2.060050505e-05, 2.060050505e-05] Y: [0.7347941366, 0.7348114609]
i: 4 X: [0.8729723781, 0.8729823329] Y: [0, 0]
i: 5 X: [-0.2959488825, -0.2958256756] Y: [-7.103772467e-05, 7.103772467e-05]

U = [0.4631611066, 0.4632214963], I = [0.4631765859, 0.4632060154],
U*(I)^(1/2)/(M)^(5/2) = [0.3152139153, 0.3152650301]
Moeckel's potential = [27.79607709, 27.80058447]

permutation: 0, 3, 2, 1, 4, 5,
symmetric with respect to OX no 2
permutation: 4, 1, 5, 3, 0, 2,
reflectional symmetry with respect to other line 2
---------------------
position 3
i: 0 X: [-1.107397433, -1.107397433] Y: [-2.586401239e-13, 2.586390677e-13]
i: 1 X: [-0.6252338707, -0.6252338707] Y: [-1.414474722e-13, 1.414467659e-13]
i: 2 X: [0.6252338707, 0.6252338707] Y: [-1.478243613e-13, 1.478249714e-13]
i: 3 X: [0.2035751658, 0.2035751658] Y: [-1.828755208e-13, 1.82876261e-13]
i: 4 X: [1.107397433, 1.107397433] Y: [0, 0]
i: 5 X: [-0.2035751658, -0.2035751658] Y: [-7.30787066e-13, 7.307874782e-13]

U = [0.5528964387, 0.5528964387], I = [0.5528964387, 0.5528964387],
U*(I)^(1/2)/(M)^(5/2) = [0.4111172398, 0.4111172398]
Moeckel's potential = [36.25298863, 36.25298863]

collinear solution no 1
permutation: 0, 1, 2, 3, 4, 5,
symmetric with respect to OX no 3
permutation: 4, 2, 1, 5, 0, 3,
reflectional symmetry with respect to other line 3
---------------------
position 5
i: 0 X: [-0.643093934, -0.643093934] Y: [0.0703689589, 0.0703689589]
i: 1 X: [-0.3886253348, -0.3886253348] Y: [-0.5758987967, -0.5758987967]
i: 2 X: [0.3059118435, 0.3059118435] Y: [-0.5700346852, -0.5700346852]
i: 3 X: [0.2374201596, 0.2374201596] Y: [0.3518298267, 0.3518298267]
i: 4 X: [0.7804401193, 0.7804401193] Y: [0, 0]
i: 5 X: [-0.2920528536, -0.2920528536] Y: [0.7237346963, 0.7237346963]

U = [0.4530097317, 0.4530097317], I = [0.4530097317, 0.4530097317],
U*(I)^(1/2)/(M)^(5/2) = [0.3049027193, 0.3049027193]
Moeckel's potential = [26.88681901, 26.88681901]

permutation: 0, 1, 2, 3, 4, 5,
there is no symmetry with respect to OX
permutation: 2, 1, 0, 3, 5, 4,
reflectional symmetry with respect to other line 4
---------------------
position 6
i: 0 X: [-0.6728086906, -0.6728043254] Y: [-3.730382939e-08, 3.730382933e-08]
i: 1 X: [-0.3364044175, -0.3364020906] Y: [-0.5826694694, -0.5826655862]
i: 2 X: [0.3364020849, 0.3364044232] Y: [-0.5826695026, -0.582665553]
i: 3 X: [0.3364021748, 0.3364043332] Y: [0.5826656192, 0.5826694364]
i: 4 X: [0.6728042418, 0.6728087742] Y: [0, 0]
i: 5 X: [-0.3364111146, -0.3363953934] Y: [0.5826616655, 0.5826733901]

U = [0.452665635, 0.4526715596], I = [0.4526640656, 0.4526731289],
U*(I)^(1/2)/(M)^(5/2) = [0.3045548607, 0.3045618957]
Moeckel's potential = [26.85614427, 26.85676463]

permutation: 0, 5, 3, 2, 4, 1,
symmetric with respect to OX no 4
permutation: 5, 3, 4, 1, 2, 0,
reflectional symmetry with respect to other line 5
---------------------
position 7
i: 0 X: [-0.6473628236, -0.6473628236] Y: [-7.162858506e-15, 6.87140444e-15]
i: 1 X: [-0.3486604816, -0.3486604816] Y: [-0.6038976687, -0.6038976687]
i: 2 X: [0.3236814118, 0.3236814118] Y: [-0.5606326507, -0.5606326507]
i: 3 X: [0.3236814118, 0.3236814118] Y: [0.5606326507, 0.5606326507]
i: 4 X: [0.6973209632, 0.6973209632] Y: [0, 0]
i: 5 X: [-0.3486604816, -0.3486604816] Y: [0.6038976687, 0.6038976687]

U = [0.4526675756, 0.4526675756], I = [0.4526675756, 0.4526675756],
U*(I)^(1/2)/(M)^(5/2) = [0.3045573471, 0.3045573471]
Moeckel's potential = [26.85636352, 26.85636352]

permutation: 0, 5, 3, 2, 4, 1,
symmetric with respect to OX no 5
permutation: 3, 4, 2, 0, 1, 5,
reflectional symmetry with respect to other line 6
---------------------
position 8
i: 0 X: [-0.610049239, -0.6100492317] Y: [0.3479737516, 0.347973759]
i: 1 X: [-0.5659045225, -0.5659045059] Y: [-0.6894027824, -0.6894027675]
i: 2 X: [0.3284053635, 0.3284053678] Y: [-0.05937476547, -0.05937475316]
i: 3 X: [0.1180995529, 0.1180995615] Y: [0.6923137261, 0.6923137306]
i: 4 X: [0.8919215792, 0.8919215814] Y: [0, 0]
i: 5 X: [-0.162472773, -0.1624727341] Y: [-0.2915099689, -0.2915099298]

U = [0.4667151069, 0.4667151301], I = [0.4667151097, 0.4667151273],
U*(I)^(1/2)/(M)^(5/2) = [0.3188436625, 0.3188436844]
Moeckel's potential = [28.11615411, 28.11615604]

permutation: 0, 1, 2, 3, 4, 5,
there is no symmetry with respect to OX
permutation: 3, 4, 5, 0, 1, 2,
reflectional symmetry with respect to other line 7
---------------------
position 11
i: 0 X: [-0.9108130342, -0.9108127712] Y: [0.4120979984, 0.4120987212]
i: 1 X: [-0.1510903126, -0.1510899819] Y: [-0.7004640475, -0.7004639356]
i: 2 X: [0.4856471594, 0.4856472373] Y: [0.01314239532, 0.01314267084]
i: 3 X: [0.01360095054, 0.01360106918] Y: [0.06305499544, 0.06305538081]
i: 4 X: [0.9997025282, 0.9997025852] Y: [0, 0]
i: 5 X: [-0.4370481387, -0.4370472913] Y: [0.2121671627, 0.2121686583]

U = [0.4980832164, 0.4980838632], I = [0.4980832935, 0.4980837861],
U*(I)^(1/2)/(M)^(5/2) = [0.3515223115, 0.3515229418]
Moeckel's potential = [30.99781067, 30.99786625]

permutation: 0, 1, 2, 3, 4, 5,
there is no symmetry with respect to OX
permutation: 4, 1, 5, 3, 0, 2,
reflectional symmetry with respect to other line 8
---------------------
position 15
i: 0 X: [-0.5941270252, -0.5941270252] Y: [0.4316585508, 0.4316585508]
i: 1 X: [0.22693633, 0.22693633] Y: [-0.6984382068, -0.6984382068]
i: 2 X: [0.22693633, 0.22693633] Y: [0.6984382068, 0.6984382068]
i: 3 X: [-7.959261057e-12, 7.959595501e-12] Y: [-1.101531975e-11, 1.101519606e-11]
i: 4 X: [0.7343813904, 0.7343813904] Y: [0, 0]
i: 5 X: [-0.5941270252, -0.5941270251] Y: [-0.4316585509, -0.4316585508]

U = [0.4494300221, 0.4494300221], I = [0.4494300221, 0.4494300221],
U*(I)^(1/2)/(M)^(5/2) = [0.3012958294, 0.3012958294]
Moeckel's potential = [26.56875757, 26.56875757]

permutation: 5, 2, 1, 3, 4, 0,
symmetric with respect to OX no 6
permutation: 0, 4, 5, 3, 1, 2,
reflectional symmetry with respect to other line 9

Number of different cc = 9
\end{verbatim}

\subsubsection*{Seven bodies}
\begin{Verbatim}[commandchars=\\\{\}]
Number of bodies = 7
Accuracy eps = 1e-05, bias = 0.01
MAX_ASYNCH_DEPTH = 64

Input data:
i: 0 X: [-6.001, 0.01] Y: [-0.001, 6.01] mass: [0.1428571429, 0.1428571429]
i: 1 X: [-6.001, 6.01] Y: [-6.001, 0.01] mass: [0.1428571429, 0.1428571429]
i: 2 X: [-6.001, 6.01] Y: [-6.001, 6.01] mass: [0.1428571429, 0.1428571429]
i: 3 X: [-6.001, 6.01] Y: [-6.001, 6.01] mass: [0.1428571429, 0.1428571429]
i: 4 X: [-6.001, 6.01] Y: [-6.001, 6.01] mass: [0.1428571429, 0.1428571429]
i: 5 X: [0.499, 6.01] Y: [0, 0] mass: [0.1428571429, 0.1428571429]
i: 6 X: [-30.06, 29.506] Y: [-24.05, 24.005] mass: [0.1428571429, 0.1428571429]


The number of undecided cubes: 0
The number of zeros in the method: 28
Program computed 111391.6327 minutes
Elapsed time = 128:35:58

Tests usage:
checkAprioriBounds -- 58539
clusterTest -- 680911214
distanceTest -- 4826082105
checkZero -- 22522665822
krawczyk: methodFailed  -- 306897559
krawczyk: zeroIside -- 28
krawczyk: no zero inside -- 633293411
spreadTest -- 0
U = I -- 234682187
improvement in clusterTest -- 1856783827
improvement in U = I Test -- 92519128

Different CC:
---------------------
position 0
i: 0 X: [-1.178582197, -1.178579025] Y: [-3.986562131e-06, 3.986562131e-06]
i: 1 X: [0.3591045716, 0.3591056814] Y: [-1.550325936e-06, 1.550325936e-06]
i: 2 X: [0.7386132282, 0.7386142695] Y: [-1.289962647e-06, 1.289962647e-06]
i: 3 X: [-5.257386353e-07, 5.257386355e-07] Y: [-1.386553701e-06, 1.386553701e-06]
i: 4 X: [-0.3591058478, -0.3591044051] Y: [-1.081032011e-06, 1.081032011e-06]
i: 5 X: [1.178580172, 1.17858105] Y: [0, 0]
i: 6 X: [-0.7386180965, -0.7386094012] Y: [-9.294436426e-06, 9.294436426e-06]

U = [0.5895856837, 0.5895908979], I = [0.5895864506, 0.589590131],
U*(I)^(1/2)/(M)^(5/2) = [0.4527106138, 0.4527160305]
Moeckel's potential = [58.6902253, 58.69092753]

collinear solution no 1
permutation: 0, 1, 2, 3, 4, 5, 6,
symmetric with respect to OX no 1
permutation: 5, 4, 6, 3, 1, 0, 2,
reflectional symmetry with respect to other line 1
---------------------
position 5
i: 0 X: [-0.9770244404, -0.9770138282] Y: [-2.316333755e-05, 2.316333755e-05]
i: 1 X: [-1.37252519e-05, 1.37252519e-05] Y: [-0.7445279538, -0.744517024]
i: 2 X: [0.4771703485, 0.4771785238] Y: [-8.739717176e-06, 8.739717176e-06]
i: 3 X: [-1.372264201e-05, 1.372264201e-05] Y: [0.74451701, 0.7445279679]
i: 4 X: [-5.827070318e-06, 5.827070318e-06] Y: [-1.221355635e-05, 1.221355635e-05]
i: 5 X: [0.9770162073, 0.9770220613] Y: [0, 0]
i: 6 X: [-0.4772200319, -0.4771288405] Y: [-5.506044198e-05, 5.506044198e-05]

U = [0.4961438918, 0.496185006], I = [0.4961530484, 0.4961758487],
U*(I)^(1/2)/(M)^(5/2) = [0.349474491, 0.3495114815]
Moeckel's potential = [45.30650705, 45.31130256]

permutation: 0, 3, 2, 1, 4, 5, 6,
symmetric with respect to OX no 2
permutation: 5, 1, 6, 3, 4, 0, 2,
reflectional symmetry with respect to other line 2
---------------------
position 7
i: 0 X: [-0.7391999133, -0.7391952076] Y: [-8.731057635e-06, 8.731057635e-06]
i: 1 X: [0.3695963649, 0.3696011956] Y: [-0.6401665036, -0.640161228]
i: 2 X: [0.3695949245, 0.369602636] Y: [0.6401607875, 0.640166944]
i: 3 X: [-3.330165813e-06, 3.330165813e-06] Y: [-5.553636666e-06, 5.553636666e-06]
i: 4 X: [-0.369605084, -0.3695924765] Y: [-0.640166307, -0.6401614246]
i: 5 X: [0.7391964306, 0.7391986904] Y: [0, 0]
i: 6 X: [-0.369618168, -0.3695793925] Y: [0.6401414238, 0.6401863077]

U = [0.4683460917, 0.468361966], I = [0.4683443213, 0.4683637363],
U*(I)^(1/2)/(M)^(5/2) = [0.3205158639, 0.3205333712]
Moeckel's potential = [41.5522581, 41.55452777]

permutation: 0, 2, 1, 3, 6, 5, 4,
symmetric with respect to OX no 3
permutation: 6, 5, 4, 3, 2, 1, 0,
reflectional symmetry with respect to other line 3
---------------------
position 8
i: 0 X: [-0.4907100136, -0.4907100136] Y: [0.8499346753, 0.8499346753]
i: 1 X: [-0.4907100136, -0.4907100136] Y: [-0.8499346753, -0.8499346753]
i: 2 X: [0.4879772078, 0.4879772078] Y: [-8.542713523e-15, 8.428297931e-15]
i: 3 X: [-2.472525901e-14, 2.503517727e-14] Y: [-3.692451758e-14, 3.674593646e-14]
i: 4 X: [-0.2439886039, -0.2439886039] Y: [-0.4226006584, -0.4226006584]
i: 5 X: [0.9814200271, 0.9814200271] Y: [0, 0]
i: 6 X: [-0.2439886039, -0.2439886039] Y: [0.4226006584, 0.4226006584]

U = [0.5148458678, 0.5148458678], I = [0.5148458678, 0.5148458678],
U*(I)^(1/2)/(M)^(5/2) = [0.3694161239, 0.3694161239]
Moeckel's potential = [47.89177651, 47.89177651]

permutation: 1, 0, 2, 3, 6, 5, 4,
symmetric with respect to OX no 4
permutation: 5, 1, 6, 3, 4, 0, 2,
reflectional symmetry with respect to other line 4
---------------------
position 9
i: 0 X: [-0.5196489643, -0.5196489643] Y: [0.8230510161, 0.8230510161]
i: 1 X: [-0.5196489643, -0.5196489643] Y: [-0.8230510161, -0.8230510161]
i: 2 X: [0.5088972465, 0.5088972465] Y: [-5.66181332e-15, 5.717243875e-15]
i: 3 X: [0.0366546652, 0.0366546652] Y: [-1.401575697e-14, 1.421911205e-14]
i: 4 X: [-0.2515755625, -0.2515755625] Y: [-0.405120033, -0.405120033]
i: 5 X: [0.996897142, 0.996897142] Y: [0, 0]
i: 6 X: [-0.2515755625, -0.2515755625] Y: [0.405120033, 0.405120033]

U = [0.5148349909, 0.5148349909], I = [0.5148349909, 0.5148349909],
U*(I)^(1/2)/(M)^(5/2) = [0.3694044172, 0.3694044172]
Moeckel's potential = [47.89025883, 47.89025883]

permutation: 1, 0, 2, 3, 6, 5, 4,
symmetric with respect to OX no 5
---------------------
position 10
i: 0 X: [-0.7346039147, -0.734603836] Y: [-1.044347371e-07, 1.04434737e-07]
i: 1 X: [-0.2560042753, -0.2560040583] Y: [-0.8539984503, -0.8539984086]
i: 2 X: [0.4114606324, 0.4114606854] Y: [-1.134599719e-07, 1.134599716e-07]
i: 3 X: [-0.04883637883, -0.04883627109] Y: [-0.350789775, -0.3507897264]
i: 4 X: [-0.04883640953, -0.04883624038] Y: [0.3507896809, 0.3507898206]
i: 5 X: [0.9328241797, 0.9328242201] Y: [0, 0]
i: 6 X: [-0.2560044998, -0.2560038338] Y: [0.8539980966, 0.8539987623]

U = [0.4885260791, 0.4885263463], I = [0.488526066, 0.4885263594],
U*(I)^(1/2)/(M)^(5/2) = [0.3414535429, 0.3414538321]
Moeckel's potential = [44.26665678, 44.26669428]

permutation: 0, 6, 2, 4, 3, 5, 1,
symmetric with respect to OX no 6
---------------------
position 15
i: 0 X: [-1.023269575, -1.023269575] Y: [0.3609518293, 0.3609518293]
i: 1 X: [-0.1222116086, -0.1222116086] Y: [-0.7138431876, -0.7138431876]
i: 2 X: [0.6226884363, 0.6226884363] Y: [0.008260404371, 0.008260404371]
i: 3 X: [0.2110036506, 0.2110036506] Y: [0.03062692087, 0.03062692087]
i: 4 X: [-0.1887985548, -0.1887985548] Y: [0.09907398529, 0.09907398529]
i: 5 X: [1.085065365, 1.085065365] Y: [0, 0]
i: 6 X: [-0.5844777135, -0.5844777135] Y: [0.2149300478, 0.2149300478]

U = [0.5351115441, 0.5351115441], I = [0.5351115441, 0.5351115441],
U*(I)^(1/2)/(M)^(5/2) = [0.3914411515, 0.3914411515]
Moeckel's potential = [50.74714105, 50.74714105]

permutation: 0, 1, 2, 3, 4, 5, 6,
there is no symmetry with respect to OX
permutation: 5, 1, 6, 4, 3, 0, 2,
reflectional symmetry with respect to other line 5
---------------------
position 17
i: 0 X: [-0.69631383, -0.6962995768] Y: [0.3878348168, 0.3878543568]
i: 1 X: [-0.4083238701, -0.4082923695] Y: [-0.7074648947, -0.7074490869]
i: 2 X: [0.3463253975, 0.346335413] Y: [-0.5997643849, -0.5997502603]
i: 3 X: [0.2547995791, 0.2548149494] Y: [0.2287497666, 0.2287768897]
i: 4 X: [0.01213209623, 0.01216749871] Y: [0.7969408167, 0.7969460439]
i: 5 X: [0.8168276037, 0.816831987] Y: [0, 0]
i: 6 X: [-0.3255579018, -0.3254469764] Y: [-0.1063779431, -0.1062961206]

U = [0.4741381696, 0.4741838696], I = [0.4741443627, 0.4741776751],
U*(I)^(1/2)/(M)^(5/2) = [0.3264827309, 0.326525669]
Moeckel's potential = [42.32581356, 42.33138012]

permutation: 0, 1, 2, 3, 4, 5, 6,
there is no symmetry with respect to OX
permutation: 4, 5, 2, 6, 0, 1, 3,
reflectional symmetry with respect to other line 6
---------------------
position 18
i: 0 X: [-0.7639543131, -0.7639543131] Y: [0.6312660195, 0.6312660195]
i: 1 X: [0.151096218, 0.151096218] Y: [-0.6987810731, -0.6987810731]
i: 2 X: [0.496401643, 0.496401643] Y: [0.03141515292, 0.03141515292]
i: 3 X: [0.049678093, 0.049678093] Y: [0.1381095023, 0.1381095023]
i: 4 X: [-0.3626530541, -0.3626530541] Y: [0.3404177807, 0.3404177807]
i: 5 X: [0.9910211804, 0.9910211804] Y: [0, 0]
i: 6 X: [-0.5615897672, -0.5615897672] Y: [-0.4424273823, -0.4424273823]

U = [0.5004058891, 0.5004058891], I = [0.5004058891, 0.5004058891],
U*(I)^(1/2)/(M)^(5/2) = [0.3539839884, 0.3539839884]
Moeckel's potential = [45.89112647, 45.89112647]

permutation: 0, 1, 2, 3, 4, 5, 6,
there is no symmetry with respect to OX
permutation: 5, 6, 4, 3, 2, 0, 1,
reflectional symmetry with respect to other line 7
---------------------
position 19
i: 0 X: [-0.8280586249, -0.8280586249] Y: [0.2436948559, 0.2436948559]
i: 1 X: [0.2762636338, 0.2762636338] Y: [-0.6696172516, -0.6696172516]
i: 2 X: [0.2976897702, 0.2976897702] Y: [0.06892828403, 0.06892828403]
i: 3 X: [0.1111253801, 0.1111253801] Y: [0.7712054397, 0.7712054397]
i: 4 X: [-0.2661193233, -0.2661193233] Y: [0.1501692824, 0.1501692824]
i: 5 X: [0.8631733714, 0.8631733714] Y: [0, 0]
i: 6 X: [-0.4540742073, -0.4540742073] Y: [-0.5643806103, -0.5643806103]

U = [0.4762002217, 0.4762002217], I = [0.4762002217, 0.4762002217],
U*(I)^(1/2)/(M)^(5/2) = [0.3286127349, 0.3286127349]
Moeckel's potential = [42.60195114, 42.60195114]

permutation: 0, 1, 2, 3, 4, 5, 6,
there is no symmetry with respect to OX
permutation: 5, 6, 4, 3, 2, 0, 1,
reflectional symmetry with respect to other line 8
---------------------
position 21
i: 0 X: [-0.7078524538, -0.7078524538] Y: [0.0404633721, 0.0404633721]
i: 1 X: [-0.4664599419, -0.4664599419] Y: [-0.7631136127, -0.7631136127]
i: 2 X: [0.3604500559, 0.3604500559] Y: [-0.1275657688, -0.1275657688]
i: 3 X: [0.3346502375, 0.3346502375] Y: [0.6250612765, 0.6250612765]
i: 4 X: [-0.07914736545, -0.07914736545] Y: [-0.3740761456, -0.3740761456]
i: 5 X: [0.894386529, 0.894386529] Y: [0, 0]
i: 6 X: [-0.3360270613, -0.3360270613] Y: [0.5992308785, 0.5992308785]

U = [0.4813750937, 0.4813750937], I = [0.4813750937, 0.4813750937],
U*(I)^(1/2)/(M)^(5/2) = [0.3339838173, 0.3339838173]
Moeckel's potential = [43.29826801, 43.29826801]

permutation: 0, 1, 2, 3, 4, 5, 6,
there is no symmetry with respect to OX
permutation: 3, 5, 4, 0, 2, 1, 6,
reflectional symmetry with respect to other line 9
---------------------
position 23
i: 0 X: [-0.7000615276, -0.7000615276] Y: [0.3020044281, 0.3020044281]
i: 1 X: [-0.151358008, -0.151358008] Y: [-0.8157984487, -0.8157984487]
i: 2 X: [0.327774992, 0.327774992] Y: [0.4028442789, 0.4028442789]
i: 3 X: [0.3190330007, 0.3190330007] Y: [-0.3836688686, -0.3836688686]
i: 4 X: [-0.1692314465, -0.1692314465] Y: [0.7434067087, 0.7434067087]
i: 5 X: [0.8297206491, 0.8297206491] Y: [0, 0]
i: 6 X: [-0.4558776597, -0.4558776597] Y: [-0.2487880983, -0.2487880983]

U = [0.4754117312, 0.4754117312], I = [0.4754117312, 0.4754117312],
U*(I)^(1/2)/(M)^(5/2) = [0.3277968993, 0.3277968993]
Moeckel's potential = [42.49618473, 42.49618473]

permutation: 0, 1, 2, 3, 4, 5, 6,
there is no symmetry with respect to OX
permutation: 4, 5, 6, 3, 0, 1, 2,
reflectional symmetry with respect to other line 10
---------------------
position 24
i: 0 X: [-0.6221369273, -0.6221369273] Y: [0.2996053536, 0.2996053536]
i: 1 X: [-0.1536551318, -0.1536551317] Y: [-0.6732071187, -0.6732071187]
i: 2 X: [0.430532112, 0.430532112] Y: [0.5398701925, 0.5398701925]
i: 3 X: [0.4305321119, 0.430532112] Y: [-0.5398701925, -0.5398701925]
i: 4 X: [-0.1536551318, -0.1536551317] Y: [0.6732071187, 0.6732071187]
i: 5 X: [0.6905198941, 0.6905198942] Y: [0, 0]
i: 6 X: [-0.6221369274, -0.6221369272] Y: [-0.2996053537, -0.2996053535]

U = [0.4768177241, 0.4768177242], I = [0.4768177241, 0.4768177243],
U*(I)^(1/2)/(M)^(5/2) = [0.3292521244, 0.3292521245]
Moeckel's potential = [42.68484275, 42.68484276]

permutation: 6, 4, 3, 2, 1, 5, 0,
symmetric with respect to OX no 7
permutation: 1, 0, 5, 4, 3, 2, 6,
reflectional symmetry with respect to other line 11
---------------------
position 25
i: 0 X: [-0.5939580764, -0.5939580764] Y: [0.3153511693, 0.3153511693]
i: 1 X: [-0.1055438498, -0.1055438498] Y: [-0.8071466615, -0.8071466615]
i: 2 X: [0.280704098, 0.280704098] Y: [0.3443258497, 0.3443258497]
i: 3 X: [0.280704098, 0.280704098] Y: [-0.3443258497, -0.3443258497]
i: 4 X: [-0.1055438498, -0.1055438498] Y: [0.8071466615, 0.8071466615]
i: 5 X: [0.8375956565, 0.8375956565] Y: [0, 0]
i: 6 X: [-0.5939580764, -0.5939580764] Y: [-0.3153511693, -0.3153511693]

U = [0.4751417478, 0.4751417478], I = [0.4751417478, 0.4751417478],
U*(I)^(1/2)/(M)^(5/2) = [0.3275177081, 0.3275177082]
Moeckel's potential = [42.45998988, 42.45998988]

permutation: 6, 4, 3, 2, 1, 5, 0,
symmetric with respect to OX no 8

Number of different cc = 14

\end{Verbatim}

\subsubsection*{Eight bodies}
\begin{Verbatim}[commandchars=\\\{\}]
(1)
i: 0 X: [-0.7166215124, -0.7166215123] Y: [0.2064500901, 0.2064500901]
i: 1 X: [0.3245586105, 0.3245586105] Y: [-0.7128796496, -0.7128796496]
i: 2 X: [0.4049275132, 0.4049275132] Y: [0.063605646, 0.06360564601]
i: 3 X: [0.001837695963, 0.001837695978] Y: [0.402546434, 0.402546434]
i: 4 X: [-0.1438928121, -0.1438928121] Y: [-0.2575373186, -0.2575373185]
i: 5 X: [-0.181586303, -0.1815863029] Y: [0.8894017544, 0.8894017544]
i: 6 X: [0.9158898541, 0.9158898541] Y: [0, 0]
i: 7 X: [-0.6051130464, -0.6051130463] Y: [-0.5915869564, -0.5915869562]

U = [0.4957241905, 0.4957241906], I = [0.4957241905, 0.4957241906],
U*(I)^(1/2)/(M)^(5/2) = [0.3490279194, 0.3490279195]
Moeckel's potential = [63.18080222, 63.18080223]

permutation: 0, 1, 2, 3, 4, 5, 6, 7,
there is no symmetry with respect to OX
\end{Verbatim}

\noindent\hrulefill
\begin{Verbatim}[commandchars=\\\{\}]
(2)
i: 0 X: [-0.758499825, -0.758499825] Y: [0.5229211679, 0.5229211679]
i: 1 X: [-0.4757362357, -0.4757362357] Y: [-0.8641489841, -0.8641489841]
i: 2 X: [0.5297384082, 0.5297384082] Y: [-0.02473539829, -0.02473539829]
i: 3 X: [0.2038227669, 0.2038227669] Y: [0.7147485589, 0.7147485589]
i: 4 X: [0.09156513918, 0.09156513918] Y: [-0.1080754295, -0.1080754295]
i: 5 X: [-0.2344231355, -0.2344231355] Y: [-0.4560136147, -0.4560136147]
i: 6 X: [0.9995383601, 0.9995383601] Y: [0, 0]
i: 7 X: [-0.3560054781, -0.3560054781] Y: [0.2153036998, 0.2153036998]

U = [0.5138288028, 0.5138288028], I = [0.5138288028, 0.5138288028],
U*(I)^(1/2)/(M)^(5/2) = [0.3683220063, 0.3683220063]
Moeckel's potential = [66.67340501, 66.67340501]

permutation: 0, 1, 2, 3, 4, 5, 6, 7,
there is no symmetry with respect to OX

\end{Verbatim}

\subsubsection*{Nine bodies}
\begin{Verbatim}[commandchars=\\\{\}]
(1)
i: 0 X: [-0.8098124055, -0.8098124055] Y: [0.03092225943, 0.03092225943]
i: 1 X: [-0.3079530509, -0.3079530509] Y: [-0.7934317414, -0.7934317414]
i: 2 X: [0.5336942253, 0.5336942253] Y: [0.05249596102, 0.05249596102]
i: 3 X: [0.3749583879, 0.3749583879] Y: [-0.6313253357, -0.6313253357]
i: 4 X: [0.1303879569, 0.1303879569] Y: [0.2002134463, 0.2002134463]
i: 5 X: [-0.1885012036, -0.1885012036] Y: [0.4988431245, 0.4988431245]
i: 6 X: [-0.2878164038, -0.2878164038] Y: [-0.2447179668, -0.2447179668]
i: 7 X: [0.9919453935, 0.9919453935] Y: [0, 0]
i: 8 X: [-0.4369028999, -0.4369028999] Y: [0.8870002527, 0.8870002527]


U = [0.5170744649, 0.5170744649], I = [0.5170744649, 0.5170744649],
U*(I)^(1/2)/(M)^(5/2) = [0.3718173376, 0.3718173376]
Moeckel's potential = [90.35161303, 90.35161303]

permutation: 0, 1, 2, 3, 4, 5, 6, 7, 8,
there is no symmetry with respect to OX
\end{Verbatim}

\noindent\hrulefill
\begin{Verbatim}[commandchars=\\\{\}]
(2)
i: 0 X: [-0.7275085134, -0.7275085134] Y: [0.1152833946, 0.1152833946]
i: 1 X: [0.01159366005, 0.01159366005] Y: [-0.8949186222, -0.8949186222]
i: 2 X: [0.5488693567, 0.5488693567] Y: [0.0194418351, 0.0194418351]
i: 3 X: [0.1519951733, 0.1519951733] Y: [-0.3967802948, -0.3967802948]
i: 4 X: [0.1389256054, 0.1389256054] Y: [0.2109317686, 0.2109317686]
i: 5 X: [-0.1649134702, -0.1649134702] Y: [0.5071400693, 0.5071400693]
i: 6 X: [-0.4146605713, -0.4146605713] Y: [0.8983002319, 0.8983002319]
i: 7 X: [1.004210874, 1.004210874] Y: [0, 0]
i: 8 X: [-0.5485121146, -0.5485121146] Y: [-0.4593983825, -0.4593983825]

U = [0.519240441, 0.519240441], I = [0.519240441, 0.519240441],
U*(I)^(1/2)/(M)^(5/2) = [0.374156044, 0.374156044]
Moeckel's potential = [90.9199187, 90.9199187]

permutation: 0, 1, 2, 3, 4, 5, 6, 7, 8,
there is no symmetry with respect to OX
\end{Verbatim}

\noindent\hrulefill
\begin{Verbatim}[commandchars=\\\{\}]
(3)
i: 0 X: [-0.7030520826, -0.7030520826] Y: [0.2581361323, 0.2581361323]
i: 1 X: [-0.6798530659, -0.6798530659] Y: [-0.8186027881, -0.818602788]
i: 2 X: [0.6455196701, 0.6455196701] Y: [-0.01022986752, -0.01022986751]
i: 3 X: [0.2525581245, 0.2525581245] Y: [-0.0456682048, -0.0456682048]
i: 4 X: [-0.005489709791, -0.005489709788] Y: [0.4378314751, 0.4378314751]
i: 5 X: [-0.09458813233, -0.09458813233] Y: [0.9180310512, 0.9180310512]
i: 6 X: [-0.09787656765, -0.09787656764] Y: [-0.2387132415, -0.2387132415]
i: 7 X: [1.072063843, 1.072063843] Y: [0, 0]
i: 8 X: [-0.3892820798, -0.3892820798] Y: [-0.5007845568, -0.5007845568]

U = [0.5375070652, 0.5375070652], I = [0.5375070652, 0.5375070652],
U*(I)^(1/2)/(M)^(5/2) = [0.394072624, 0.3940726241]
Moeckel's potential = [95.75964764, 95.75964765]

permutation: 0, 1, 2, 3, 4, 5, 6, 7, 8,
there is no symmetry with respect to OX

\end{Verbatim}

\subsubsection*{Ten bodies}
\begin{Verbatim}[commandchars=\\\{\}]
(1)
i: 0 X: [-0.7007408215, -0.7007408215] Y: [0.4693295286, 0.4693295287]
i: 1 X: [-0.07639930863, -0.07639930862] Y: [-0.9073170049, -0.9073170048]
i: 2 X: [0.5132734942, 0.5132734942] Y: [0.5961143132, 0.5961143132]
i: 3 X: [0.4588100327, 0.4588100327] Y: [-0.1317284788, -0.1317284788]
i: 4 X: [0.1456294856, 0.1456294856] Y: [-0.4767952142, -0.4767952142]
i: 5 X: [-0.002441745604, -0.002441745595] Y: [0.2214431849, 0.2214431849]
i: 6 X: [-0.1155249827, -0.1155249827] Y: [0.8058001754, 0.8058001754]
i: 7 X: [-0.4929493421, -0.4929493421] Y: [-0.04332608545, -0.04332608544]
i: 8 X: [0.9269083278, 0.9269083278] Y: [0, 0]
i: 9 X: [-0.6565651399, -0.6565651398] Y: [-0.5335204189, -0.5335204189]

U = [0.5167031585, 0.5167031586], I = [0.5167031585, 0.5167031586],
U*(I)^(1/2)/(M)^(5/2) = [0.3714169116, 0.3714169117]
Moeckel's potential = [117.4523402, 117.4523402]

permutation: 0, 1, 2, 3, 4, 5, 6, 7, 8, 9,
there is no symmetry with respect to OX

\end{Verbatim}

\noindent\hrulefill
\begin{Verbatim}[commandchars=\\\{\}]
(2)
i: 0 X: [-0.8963214942, -0.8963214942] Y: [0.01882076642, 0.01882076642]
i: 1 X: [-0.2944424421, -0.2944424421] Y: [-0.7961815748, -0.7961815748]
i: 2 X: [0.4843083474, 0.4843083474] Y: [-0.6834013227, -0.6834013227]
i: 3 X: [0.4538131131, 0.4538131131] Y: [0.08780569651, 0.08780569651]
i: 4 X: [0.1725919232, 0.1725919232] Y: [0.460572645, 0.460572645]
i: 5 X: [0.07492379596, 0.07492379596] Y: [-0.289147328, -0.289147328]
i: 6 X: [-0.0101214633, -0.0101214633] Y: [0.9213652726, 0.9213652726]
i: 7 X: [-0.4311240011, -0.4311240011] Y: [-0.2030210622, -0.2030210622]
i: 8 X: [0.9304861223, 0.9304861223] Y: [0, 0]
i: 9 X: [-0.4841139013, -0.4841139013] Y: [0.483186907, 0.483186907]

U = [0.5180473051, 0.5180473051], I = [0.5180473051, 0.5180473051],
U*(I)^(1/2)/(M)^(5/2) = [0.3728671543, 0.3728671543]
Moeckel's potential = [117.9109472, 117.9109472]

permutation: 0, 1, 2, 3, 4, 5, 6, 7, 8, 9,
there is no symmetry with respect to OX

\end{Verbatim}

\noindent\hrulefill
\begin{Verbatim}[commandchars=\\\{\}]
(3)
i: 0 X: [-0.7967875604, -0.7967875604] Y: [0.1958277257, 0.1958277257]
i: 1 X: [0.2919426026, 0.2919426026] Y: [-0.7600314103, -0.7600314103]
i: 2 X: [0.5006139893, 0.5006139893] Y: [-0.0465575355, -0.0465575355]
i: 3 X: [0.4399251191, 0.4399251191] Y: [0.6127399142, 0.6127399142]
i: 4 X: [0.06874098874, 0.06874098874] Y: [-0.1558927224, -0.1558927224]
i: 5 X: [-0.1936088183, -0.1936088183] Y: [0.8476831709, 0.8476831709]
i: 6 X: [-0.2311572039, -0.2311572039] Y: [0.3226693268, 0.3226693268]
i: 7 X: [-0.3368118774, -0.3368118774] Y: [-0.367337145, -0.367337145]
i: 8 X: [0.9623051252, 0.9623051252] Y: [0, 0]
i: 9 X: [-0.705162365, -0.705162365] Y: [-0.6491013244, -0.6491013244]

U = [0.5193488033, 0.5193488033], I = [0.5193488033, 0.5193488033],
 U*(I)^(1/2)/(M)^(5/2) = [0.3742731763, 0.3742731763]
Moeckel's potential = [118.3555704, 118.3555704]

permutation: 0, 1, 2, 3, 4, 5, 6, 7, 8, 9,
there is no symmetry with respect to OX

\end{Verbatim}

\noindent\hrulefill
\begin{Verbatim}[commandchars=\\\{\}]
(4)
i: 0 X: [-0.8675091418, -0.8675091417] Y: [0.272352247, 0.2723522471]
i: 1 X: [0.6017849903, 0.6017849904] Y: [-0.6658290959, -0.6658290959]
i: 2 X: [0.4701929967, 0.4701929968] Y: [0.1559215689, 0.1559215689]
i: 3 X: [0.1519788533, 0.1519788533] Y: [-0.476594594, -0.4765945939]
i: 4 X: [0.1008116833, 0.1008116836] Y: [0.894979045, 0.8949790451]
i: 5 X: [0.0754051335, 0.07540513365] Y: [0.3753489653, 0.3753489653]
i: 6 X: [-0.292654655, -0.2926546549] Y: [-0.4085342723, -0.4085342722]
i: 7 X: [-0.3838258452, -0.3838258452] Y: [0.29952652, 0.2995265201]
i: 8 X: [0.9222584288, 0.9222584289] Y: [0, 0]
i: 9 X: [-0.7784424448, -0.7784424439] Y: [-0.4471703844, -0.447170384]

U = [0.5231657667, 0.5231657669], I = [0.5231657667, 0.5231657669],
U*(I)^(1/2)/(M)^(5/2) = [0.3784068394, 0.3784068396]
Moeckel's potential = [119.6627495, 119.6627495]

permutation: 0, 1, 2, 3, 4, 5, 6, 7, 8, 9,
there is no symmetry with respect to OX

\end{Verbatim}

\noindent\hrulefill
\begin{Verbatim}[commandchars=\\\{\}]
(5)
i: 0 X: [-0.7888545233, -0.7888545233] Y: [0.2769354218, 0.2769354218]
i: 1 X: [0.1573661438, 0.1573661438] Y: [-0.8879803722, -0.8879803722]
i: 2 X: [0.566350671, 0.566350671] Y: [0.02117690755, 0.02117690755]
i: 3 X: [0.1809819542, 0.1809819542] Y: [-0.3943353946, -0.3943353946]
i: 4 X: [0.1673309159, 0.1673309159] Y: [0.2071064401, 0.2071064401]
i: 5 X: [-0.08156809045, -0.08156809045] Y: [0.5470717464, 0.5470717464]
i: 6 X: [-0.2490045724, -0.2490045724] Y: [0.9598108791, 0.9598108791]
i: 7 X: [-0.3968333769, -0.3968333769] Y: [-0.1057160826, -0.1057160826]
i: 8 X: [1.005616634, 1.005616634] Y: [0, 0]
i: 9 X: [-0.5613857558, -0.5613857558] Y: [-0.6240695457, -0.6240695457]

U = [0.5266322308, 0.5266322308], I = [0.5266322308, 0.5266322308],
U*(I)^(1/2)/(M)^(5/2) = [0.3821740131, 0.3821740131]
Moeckel's potential = [120.8540344, 120.8540344]

permutation: 0, 1, 2, 3, 4, 5, 6, 7, 8, 9,
there is no symmetry with respect to OX

\end{Verbatim}

\noindent\hrulefill
\begin{Verbatim}[commandchars=\\\{\}]
(6)
position 0
i: 0 X: [-0.8111684557, -0.8111684557] Y: [0.2656551367, 0.2656551367]
i: 1 X: [-0.3059600863, -0.3059600863] Y: [-0.9399118629, -0.9399118629]
i: 2 X: [0.5694847219, 0.5694847219] Y: [-0.02355906915, -0.02355906915]
i: 3 X: [0.2319114215, 0.2319114215] Y: [0.3782046907, 0.3782046907]
i: 4 X: [0.2014838749, 0.2014838749] Y: [-0.2518381211, -0.2518381211]
i: 5 X: [0.1490178487, 0.1490178487] Y: [0.8945560864, 0.8945560864]
i: 6 X: [-0.09081819785, -0.09081819785] Y: [-0.5373390538, -0.5373390538]
i: 7 X: [-0.3293580397, -0.3293580397] Y: [0.4819496801, 0.4819496801]
i: 8 X: [1.01098039, 1.01098039] Y: [0, 0]
i: 9 X: [-0.6255734779, -0.6255734779] Y: [-0.2677174871, -0.2677174871]

U = [0.5276587071, 0.5276587071], I = [0.5276587071, 0.5276587071],
U*(I)^(1/2)/(M)^(5/2) = [0.3832919194, 0.3832919194]
Moeckel's potential = [121.2075474, 121.2075474]

permutation: 0, 1, 2, 3, 4, 5, 6, 7, 8, 9,
there is no symmetry with respect to OX

\end{Verbatim}

\noindent\hrulefill
\begin{Verbatim}[commandchars=\\\{\}]
(7)
i: 0 X: [-0.8644611425, -0.8644610843] Y: [0.2149252766, 0.2149253394]
i: 1 X: [-0.00513249289, -0.005132358042] Y: [-0.9138239369, -0.9138239261]
i: 2 X: [0.538361395, 0.5383614046] Y: [0.09282895599, 0.09282897502]
i: 3 X: [0.46480923, 0.4648092392] Y: [-0.5151381858, -0.5151381172]
i: 4 X: [0.1823759207, 0.182375933] Y: [0.2881651807, 0.2881652273]
i: 5 X: [-0.1446907068, -0.1446907] Y: [0.5315697209, 0.5315697332]
i: 6 X: [-0.2305530662, -0.230553035] Y: [-0.4706217805, -0.4706217273]
i: 7 X: [-0.4362270962, -0.4362270088] Y: [0.8774029291, 0.8774029878]
i: 8 X: [0.9865289623, 0.9865289661] Y: [0, 0]
i: 9 X: [-0.4910113567, -0.4910110035] Y: [-0.1053084921, -0.1053081602]

U = [0.528844359, 0.5288444955], I = [0.5288443794, 0.528844475],
U*(I)^(1/2)/(M)^(5/2) = [0.3845845407, 0.3845846747]
Moeckel's potential = [121.6163101, 121.6163525]

permutation: 0, 1, 2, 3, 4, 5, 6, 7, 8, 9,
there is no symmetry with respect to OX

\end{Verbatim}

\noindent\hrulefill
\begin{Verbatim}[commandchars=\\\{\}]
(8)
i: 0 X: [-0.8620127236, -0.8620127236] Y: [0.5749149667, 0.5749149667]
i: 1 X: [-0.003887282198, -0.003887282198] Y: [-0.9423644458, -0.9423644458]
i: 2 X: [0.6244444875, 0.6244444875] Y: [-0.003510447333, -0.003510447333]
i: 3 X: [0.2389418223, 0.2389418223] Y: [0.7433087445, 0.7433087445]
i: 4 X: [0.2253692474, 0.2253692474] Y: [-0.004873923735, -0.004873923735]
i: 5 X: [0.04037831843, 0.04037831843] Y: [-0.4785847827, -0.4785847827]
i: 6 X: [-0.1523931658, -0.1523931658] Y: [0.1569749357, 0.1569749357]
i: 7 X: [-0.5005974321, -0.5005974321] Y: [0.3480983721, 0.3480983721]
i: 8 X: [1.045810599, 1.045810599] Y: [0, 0]
i: 9 X: [-0.656053871, -0.656053871] Y: [-0.3939634194, -0.3939634194]

U = [0.5341660935, 0.5341660935], I = [0.5341660935, 0.5341660935],
U*(I)^(1/2)/(M)^(5/2) = [0.3904041955, 0.3904041955]
Moeckel's potential = [123.4566466, 123.4566466]

permutation: 0, 1, 2, 3, 4, 5, 6, 7, 8, 9,
there is no symmetry with respect to OX

\end{Verbatim}

\noindent\hrulefill
\begin{Verbatim}[commandchars=\\\{\}]
(9)
i: 0 X: [-0.7437201474, -0.7437201462] Y: [0.05424247592, 0.05424247654]
i: 1 X: [-0.6097575083, -0.6097575063] Y: [-0.8688275537, -0.868827552]
i: 2 X: [0.6456942806, 0.6456942808] Y: [-0.03032211044, -0.03032210977]
i: 3 X: [0.3841057237, 0.3841057254] Y: [0.7124794069, 0.7124794076]
i: 4 X: [0.2719022074, 0.2719022076] Y: [-0.1083651492, -0.108365148]
i: 5 X: [-0.06052210988, -0.06052210964] Y: [-0.2908103139, -0.2908103134]
i: 6 X: [-0.09397879995, -0.09397879682] Y: [0.3792041934, 0.3792041966]
i: 7 X: [-0.3436577475, -0.3436577466] Y: [-0.5452330647, -0.5452330637]
i: 8 X: [1.062255195, 1.062255195] Y: [0, 0]
i: 9 X: [-0.5123211036, -0.5123210936] Y: [0.6976321062, 0.6976321156]

U = [0.5375196597, 0.5375196629], I = [0.5375196593, 0.5375196632],
U*(I)^(1/2)/(M)^(5/2) = [0.3940864744, 0.3940864782]
Moeckel's potential = [124.6210854, 124.6210866]

permutation: 0, 1, 2, 3, 4, 5, 6, 7, 8, 9,
there is no symmetry with respect to OX

\end{Verbatim}

\noindent\hrulefill
\begin{Verbatim}[commandchars=\\\{\}]
(10)
i: 0 X: [-0.7443869099, -0.7443869099] Y: [0.004965386648, 0.004965386648]
i: 1 X: [0.07184450244, 0.07184450244] Y: [-0.910877299, -0.910877299]
i: 2 X: [0.6495486297, 0.6495486297] Y: [0.01650001945, 0.01650001945]
i: 3 X: [0.2638544878, 0.2638544878] Y: [0.08972086605, 0.08972086605]
i: 4 X: [0.1507954307, 0.1507954307] Y: [-0.4264727703, -0.4264727703]
i: 5 X: [-0.0449474494, -0.0449474494] Y: [0.3011792823, 0.3011792823]
i: 6 X: [-0.3224455973, -0.3224455973] Y: [0.5582706437, 0.5582706437]
i: 7 X: [-0.5083858126, -0.5083858126] Y: [-0.5225251633, -0.5225251633]
i: 8 X: [1.065046288, 1.065046288] Y: [0, 0]
i: 9 X: [-0.5809235696, -0.5809235696] Y: [0.8892390345, 0.8892390345]

U = [0.5395865635, 0.5395865635], I = [0.5395865635, 0.5395865635],
U*(I)^(1/2)/(M)^(5/2) = [0.3963617068, 0.3963617068]
Moeckel's potential = [125.3405771, 125.3405771]

permutation: 0, 1, 2, 3, 4, 5, 6, 7, 8, 9,
there is no symmetry with respect to OX

\end{Verbatim}

\noindent\hrulefill
\begin{Verbatim}[commandchars=\\\{\}]
(11)
i: 0 X: [-0.8333907632, -0.8333907632] Y: [0.7644395423, 0.7644395423]
i: 1 X: [-0.03372694114, -0.03372694114] Y: [-0.9240673101, -0.9240673101]
i: 2 X: [0.7379960306, 0.7379960306] Y: [0.009114360497, 0.009114360497]
i: 3 X: [0.386475537, 0.386475537] Y: [0.03160247231, 0.03160247231]
i: 4 X: [0.04528434434, 0.04528434434] Y: [0.1213495026, 0.1213495026]
i: 5 X: [0.018708438, 0.018708438] Y: [-0.4551657158, -0.4551657158]
i: 6 X: [-0.2519844233, -0.2519844233] Y: [0.2944278672, 0.2944278672]
i: 7 X: [-0.5319428417, -0.5319428417] Y: [0.5069057403, 0.5069057403]
i: 8 X: [1.132810121, 1.132810121] Y: [0, 0]
i: 9 X: [-0.6702295016, -0.6702295016] Y: [-0.3486064594, -0.3486064594]

U = [0.559742471, 0.559742471], I = [0.559742471, 0.559742471],
U*(I)^(1/2)/(M)^(5/2) = [0.4187765849, 0.4187765849]
Moeckel's potential = [132.4287839, 132.4287839]

permutation: 0, 1, 2, 3, 4, 5, 6, 7, 8, 9,
there is no symmetry with respect to OX

\end{Verbatim}

\subsection{Listing of central configurations}
\label{sec:listing}

In this and next sections we present graphically all CCs for $n = \{3, \ldots, 7\}$ and asymmetrical CCs for $n = \{8, 9, 10\}$ listed in~\cite{F02} and reproduced in our proofs. For any CC we give the picture of it (with red lines showing symmetry lines), an invariant $J$ (see~(\ref{eq:inv})), its coordinates and Moeckel's potential $P$. The coordinates for CCs are middle points from the interval bounds obtained in our proof, hence
in the rigorous sense these are just an approximation of the true CC. Moreover, despite that we usually give real  numbers with ten decimal digits
we make no claim that all these digits correct, in fact only 4 or 5 first decimals are correct.

In the case of 3 and 4-body problem all CCs are symmetrical with respect to the $0X$ axis (this is the axis passing through the origin and the body most distant from the origin).
In the case of 5, 6 and 7-body problems all CCs are still symmetrical with respect to some line, but there exists CCs which are not symmetrical with respect to the $OX$ axis. The first asymmetrical CCs appear for $n = 8$.

\counterwithin{figure}{subsubsection}
\subsubsection{Three bodies}\label{sec:3-body}
\noindent
\config{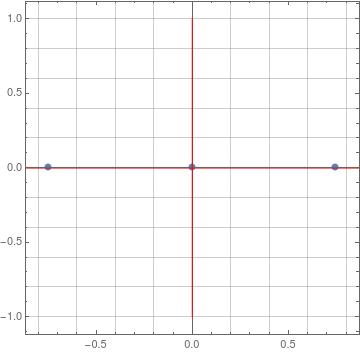}{collinear}{$
    \begin{array}{lr}
    (-0.7469007911, & -1.311948311\times 10^{-15})\\
    (0.7469007911, & 0)\\
    (-2.331468352\times 10^{-15}, & -1.328983689\times 10^{-15})
\end{array}
\mbox{}\\[1ex]
\begin{array}{lcl}
J & = & 0.2268046058\\
P & = & 3.535533906
    \end{array}
    $
}
\par\medskip
\noindent
\config{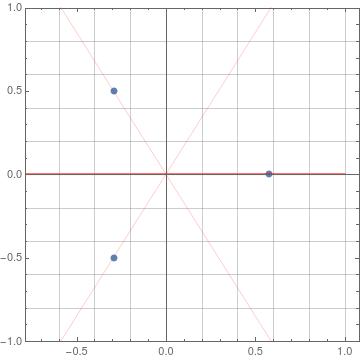}{equilateral triangle}{
    $
\begin{array}{lr}
(-0.2886751346, & -0.5)\\
(0.5773502692, & 0)\\
(-0.2886751346, & 0.5)
\end{array}
\mbox{}\\[1ex]
\begin{array}{lcl}
J & = & 0.1924500897\\
P & = & 3
\end{array}
    $
}

\subsubsection{Four bodies}

\config{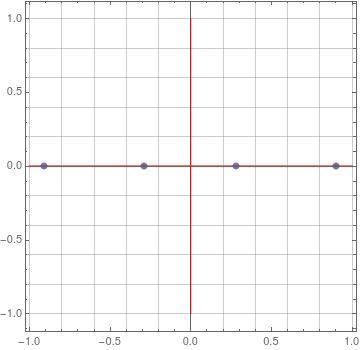}{collinear}{
    $
   \begin{array}{lr}
(-0.9051285388, & -5.049805694\times 10^{-9})\\
(-0.2862410122, & -1.476723831\times 10^{-9})\\
(0.9051285343, & 0)\\
(0.2862410037, & -6.526529529\times 10^{-9})
\end{array}
\mbox{}\\[1ex]
\begin{array}{lcl}
J & = & 0.3024688765\\
P & = & 9.679004022
    \end{array}
    $
    }
\par\medskip
\noindent
\config{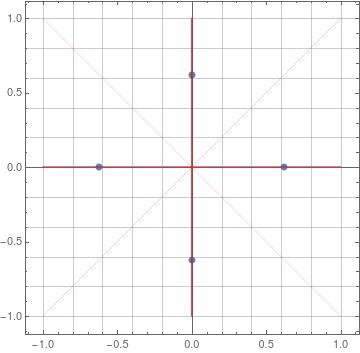}{square (cross)}{
    $
   \begin{array}{lr}
    (-0.6208313565, & -1.338276202\times 10^{-5})\\
    (-9.896654491\times 10^{-6}, & -0.620827361)\\
(0.6208135881, & 0)\\
(-2.766510061\times 10^{-5}, & 0.6208024044)
\end{array}
\mbox{}\\[1ex]
\begin{array}{lcl}
J & = & 0.2392648356\\
P & = & 7.65647474
    \end{array}
    $
}
\par\medskip
\noindent
\config{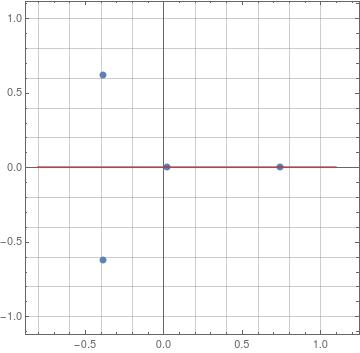}{isoseles triangle}{
    $
   \begin{array}{lr}
(-0.3821936947, & 0.6195346528)\\
(-0.3821936947, & -0.6195346528)\\
(0.7436490828, & 0)\\
(0.0207383067, & -5.828892924\times 10^{-12})
\end{array}
\mbox{}\\[1ex]
\begin{array}{lcl}
J & = & 0.2561261996\\
P & = & 8.196080629
    \end{array}
    $
}
\par\medskip
\noindent
\config{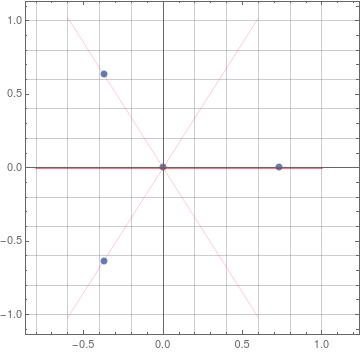}{equilateral}{
    $
   \begin{array}{lr}
(-0.3666565002, & 0.6350676872)\\
(-0.3666565002, & -0.6350676872)\\
(0.7333130003, & 0)\\
(-2.157363177\times 10^{-11}, & -6.115996598\times 10^{-12})
\end{array}
\mbox{}\\[1ex]
\begin{array}{lcl}
J & = & 0.2561297548\\
P & = & 8.196152422
    \end{array}
    $
}

\subsubsection{Five bodies}
\config{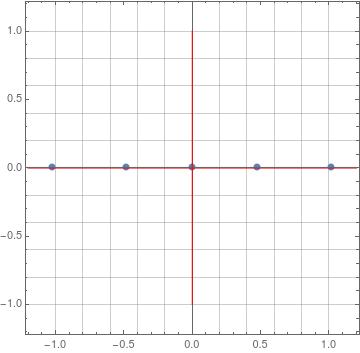}{collinear}{
    $
   \begin{array}{lr}
(-1.019255982, & -4.721349739\times 10^{-15})\\
(-0.480767439, & -2.083289376\times 10^{-15})\\
(0.480767439, & -4.095121196\times 10^{-15})\\
(1.019255982, & 0)\\
(-9.714451465\times 10^{-15}, & -1.0901527\times 10^{-14})
\end{array}
\mbox{}\\[1ex]
\begin{array}{lcl}
J & = & 0.3620810967\\
P & = & 20.24094955
    \end{array}
    $
}
\par\medskip
\noindent
\config{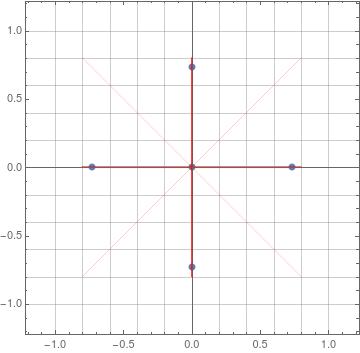}{cross}{
    $
   \begin{array}{lr}
(-0.7315026092, & -1.22322208\times 10^{-5})\\
(-4.380645756\times 10^{-6}, & -0.7315035144)\\
(-3.486361374\times 10^{-6}, & -6.548894919\times 10^{-6})\\
(0.7314993242, & 0)\\
(-1.115205008\times 10^{-5}, & 0.7314794181)
\end{array}
\mbox{}\\[1ex]
\begin{array}{lcl}
J & = &  0.2800671822)\\
P & = & 15.65645267
    \end{array}
    $
}
\par\medskip
\noindent
\config{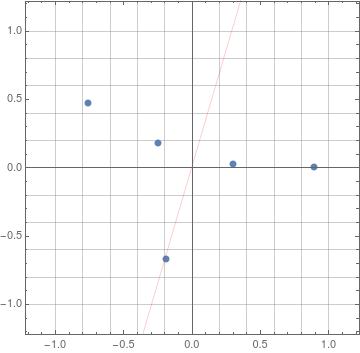}{two isosceles triangles}{
    $
   \begin{array}{lr}
(-0.7609899148, & 0.4684633687)\\
(-0.1902473651, & -0.6719536689)\\
(0.3050645568, & 0.02352925245)\\
(0.8936238356, & 0)\\
(-0.2474515264, & 0.1799604652)
\end{array}
\mbox{}\\[1ex]
\begin{array}{lcl}
J & = & 0.306302609\\
P & = & 17.1239885
    \end{array}
    $
}
\par\medskip
\noindent
\config{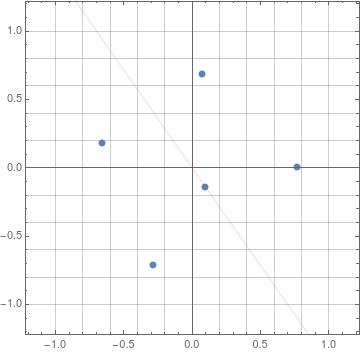}{trapezium}{
    $
   \begin{array}{lr}
(-0.6591405438, & 0.1800138879)\\
(-0.2807232123, & -0.7143993153)\\
(0.09876565526, & -0.1449273191)\\
(0.767575262, & 0)\\
(0.07352247002, & 0.6793125102)
\end{array}
\mbox{}\\[1ex]
\begin{array}{lcl}
J & = & 0.280563341\\
P & = & 15.68396719
    \end{array}
    $
}
\par\medskip
\noindent
\config{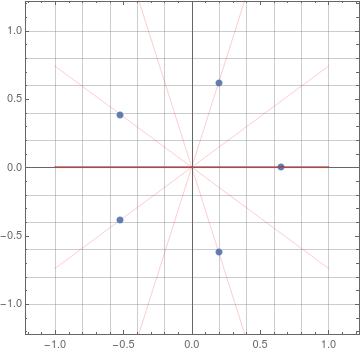}{regular pentagon}{
    $
   \begin{array}{lr}
(-0.5262827178, & 0.3823582109)\\
(0.2010155547, & -0.6186815315)\\
(0.2010162837, & 0.6186717554)\\
(0.6505104284, & 0)\\
(-0.526293415, & -0.3823760565)
\end{array}
\mbox{}\\[1ex]
\begin{array}{lcl}
J & = & 0.2752763841\\
P & = & 15.38795199
    \end{array}
    $
}

\subsubsection{Six bodies}

\config{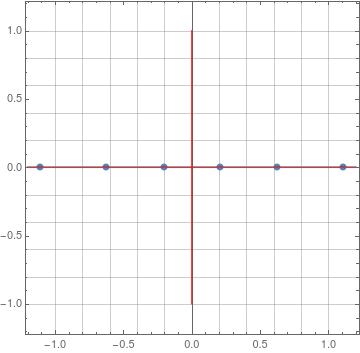}{collinear}{
    $
    \begin{array}{ll}
(-1.107397812, & -1.475347758\times 10^{-6})\\
(-0.6252343397, & -8.060205112\times 10^{-7})\\
(0.6252335566, & -8.499503536\times 10^{-7})\\
(0.2035747752, & -1.050058972\times 10^{-6})\\
(1.107397183, & 0)\\
(-0.2035769692, & -4.181377595\times 10^{-6})
\end{array}
\mbox{}\\[1ex]
\begin{array}{lcl}
J & = & 0.41111724\\
P & = & 36.25288592
    \end{array}
    $
}
\par\medskip
\noindent
\config{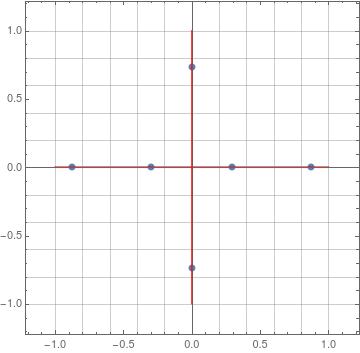}{cross}{
    $
   \begin{array}{lr}
(-0.8729803868, & -1.299986789\times 10^{-5})\\
(-8.494473347\times 10^{-6}, & -0.734806382)\\
(0.2958840069, & -9.097289147\times 10^{-6})\\
(-8.468261278\times 10^{-6}, & 0.7347992334)\\
(0.8729753137, & 0)\\
(-0.295912587, & -2.924581908\times 10^{-5})
\end{array}
\mbox{}\\[1ex]
\begin{array}{lcl}
J & = & 0.315239471\\
P & = & 27.7974043
    \end{array}
    $
}
\par\medskip
\noindent
\config{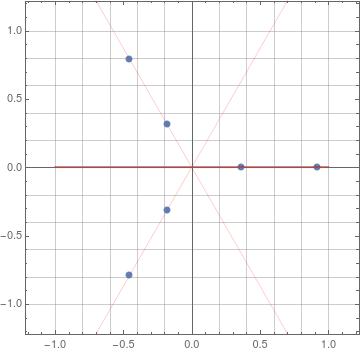}{two equilaterals}{
    $
   \begin{array}{lr}
(-0.4566158049, & 0.7908817243)\\
(-0.4566157923, & -0.7908817398)\\
(0.3632321467, & -1.247155263\times 10^{-8})\\
(-0.1816160799, & -0.3145682734)\\
(0.9132315629, & 0)\\
(-0.1816161185, & 0.3145682393)
\end{array}
\mbox{}\\[1ex]
\begin{array}{lcl}
J & = & 0.335639561\\
P & = & 29.59724201
    \end{array}
    $
}
\par\medskip
\noindent
\config{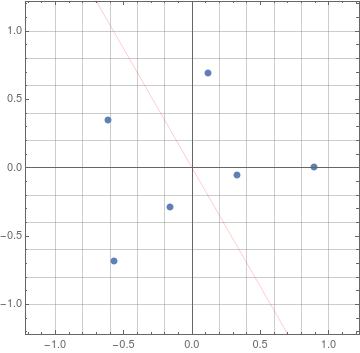}{`trapezium'}{
    $
   \begin{array}{lr}
(-0.6100492386, & 0.3479737521)\\
(-0.5659045216, & -0.6894027816)\\
(0.3284053638, & -0.05937476467)\\
(0.1180995535, & 0.6923137264)\\
(0.8919215793, & 0)\\
(-0.1624727707, & -0.2915099666)
\end{array}
\mbox{}\\[1ex]
\begin{array}{lcl}
J & = & 0.318843673\\
P & = & 28.11615423
    \end{array}
    $
}
\par\medskip
\noindent
\config{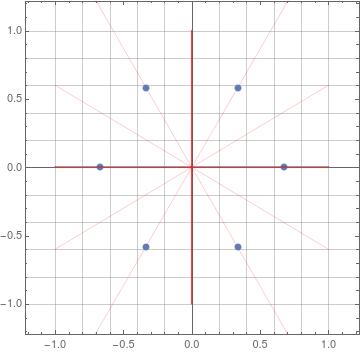}{regular hexagon}{
    $
   \begin{array}{lr}
(-0.672806508, & -4.272099752\times 10^{-15})\\
(-0.336403254, & -0.5826675278)\\
(0.336403254, & -0.5826675278)\\
(0.336403254, & 0.5826675278)\\
(0.672806508, & 0)\\
(-0.336403254, & 0.5826675278)
\end{array}
\mbox{}\\[1ex]
\begin{array}{lcl}
J & = & 0.304558377\\
P & = & 26.85645445
    \end{array}
    $
}
\par\medskip
\noindent
\config{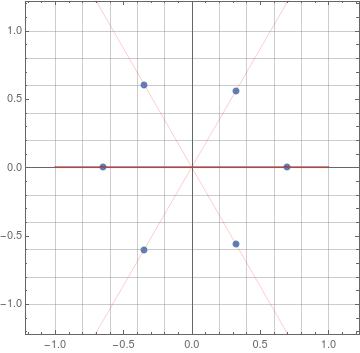}{two equilaterals}{
    $
   \begin{array}{lr}
(-0.6473628823, & -1.536934675\times 10^{-9})\\
(-0.3486605058, & -0.6038977133)\\
(0.3236813862, & -0.5606326965)\\
(0.3236813816, & 0.5606326008)\\
(0.6973209126, & 0)\\
(-0.3486606711, & 0.603897527)
\end{array}
\mbox{}\\[1ex]
\begin{array}{lcl}
J & = & 0.304557347\\
P & = & 26.85635599
    \end{array}
    $
}
\par\medskip
\noindent
\config{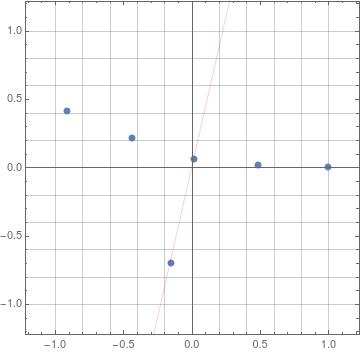}{two isosceles triangles}{
    $
   \begin{array}{lr}
(-0.9108129056, & 0.4120983519)\\
(-0.1510901508, & -0.7004639928)\\
(0.4856471975, & 0.01314253007)\\
(0.01360100857, & 0.06305518391)\\
(0.9997025561, & 0)\\
(-0.4370477243, & 0.2121678941)
\end{array}
\mbox{}\\[1ex]
\begin{array}{lcl}
J & = & 0.351522627\\
P & = & 30.99783785
    \end{array}
    $
}
\par\medskip
\noindent
\config{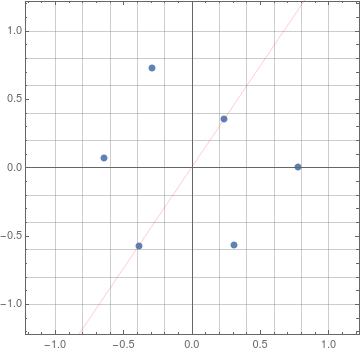}{two isosceles triangles}{$
\begin{array}{lcl}
(-0.643093934, & 0.0703689589)\\
(-0.3886253348, & -0.5758987967)\\
(0.3059118435, & -0.5700346852)\\
(0.2374201596, & 0.3518298267)\\
(0.7804401193, & 0)\\
(-0.2920528536, & 0.7237346963)
\end{array}
\mbox{}\\[1ex]
\begin{array}{lcl}
J & = & 0.304902719\\
P & = & 26.88681901
\end{array}
$
}
\par\medskip
\noindent
\config{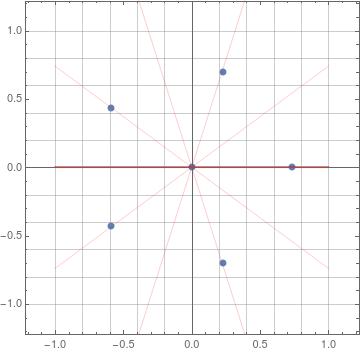}{regular pentagon}{
$
\begin{array}{lcl}
(-0.5941270252, & 0.4316585508)\\
(0.2269363299, & -0.6984382068)\\
(0.22693633, & 0.6984382067)\\
(-2.679815408\times 10^{-11}, & -3.713999317\times 10^{-11})\\
(0.7343813904, & 0)\\
(-0.5941270253, & -0.4316585509)
\end{array}
\mbox{}\\[1ex]
\begin{array}{lcl}
J & = & 0,301295829)\\
P & = & 26.56875757
\end{array}$
}

\subsubsection{Seven bodies}

\config{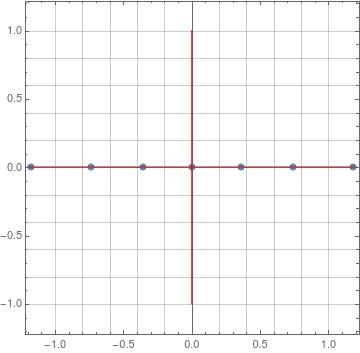}{collinear}{$
\begin{array}{lcl}
(-1.178582197, & -3.986562131\times 10^{-6})\\
(0.3591045716, & -1.550325936\times 10^{-6})\\
(0.7386132282, & -1.289962647\times 10^{-6})\\
(-5.257386353\times 10^{-7}, & -1.386553701\times 10^{-6})\\
(-0.3591058478, & -1.081032011\times 10^{-6})\\
(1.178580172, & 0)\\
(-0.7386180965, & -9.294436426\times 10^{-6})
\end{array}
\mbox{}\\[1ex]
\begin{array}{lcl}
J &  = & 0.452713322)\\
P & = & 58.6902253
\end{array}
$}
\par\medskip
\noindent
\config{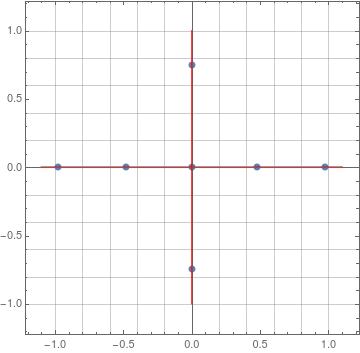}{cross}{
$
\begin{array}{lcl}
(-0.9770244404, & -2.316333755\times 10^{-5})\\
(-1.37252519\times 10^{-5}, & -0.7445279538)\\
(0.4771703485, & -8.739717176\times 10^{-6})\\
(-1.372264201\times 10^{-5}, & 0.74451701)\\
(-5.827070318\times 10^{-6}, & -1.221355635\times 10^{-5})\\
(0.9770162073, & 0)\\
(-0.4772200319, & -5.506044198\times 10^{-5})
\end{array}
\mbox{}\\[1ex]
\begin{array}{lcl}
J & = & 0.3493930453\\
P & = & 45.30650705
\end{array}
$
}
\par\medskip
\noindent
\config{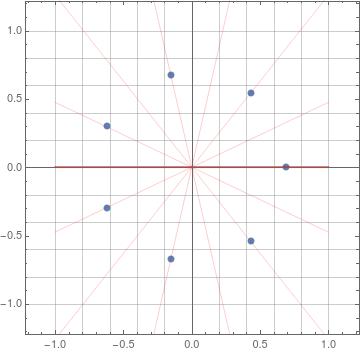}{regular heptagon}{$
\begin{array}{lr}
(-0.6221369273, & 0.2996053536)\\
(-0.1536551318, & -0.6732071187)\\
(0.430532112, & 0.5398701925)\\
(0.4305321119, & -0.5398701925)\\
(-0.1536551318, & 0.6732071187)\\
(0.6905198941, & 0)\\
(-0.6221369274, & -0.2996053537)
\end{array}
\mbox{}\\[1ex]
\begin{array}{lcl}
J & = & 0.329252124\\
P & = & 42.68484275
\end{array}
$}
\par\medskip
\noindent
\config{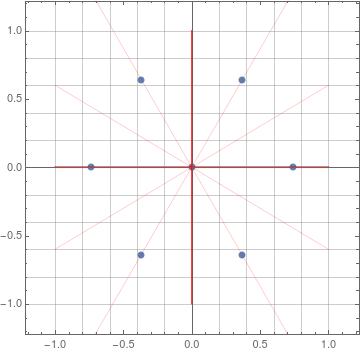}{regular hexagon}{$
\begin{array}{lr}
(-0.7391999133, & -8.731057635\times 10^{-6})\\
(0.3695963649, & -0.6401665036)\\
(0.3695949245, & 0.6401607875)\\
(-3.330165813\times 10^{-6}, & -5.553636666\times 10^{-6})\\
(-0.369605084, & -0.640166307)\\
(0.7391964306, & 0)\\
(-0.369618168, & 0.6401414238)
\end{array}
\mbox{}\\[1ex]
\begin{array}{lcl}
J & = & 0.3205246151\\
P & = & 41.5522581
\end{array}
$}
\par\medskip
\noindent
\config{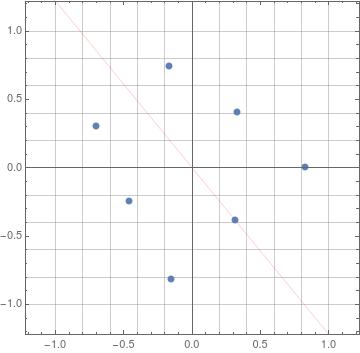}{`trapezium'}{$
\begin{array}{lcl}
(-0.7000615276, & 0.3020044281)\\
(-0.151358008, & -0.8157984487)\\
(0.327774992, & 0.4028442789)\\
(0.3190330007, & -0.3836688686)\\
(-0.1692314465, & 0.7434067087)\\
(0.8297206491, & 0)\\
(-0.4558776597, & -0.2487880983)
\end{array}
\mbox{}\\[1ex]
\begin{array}{lcl}
J & = & 0.3277968993\\
P & = & 42.49618473
\end{array}$}
\par\medskip
\noindent
\config{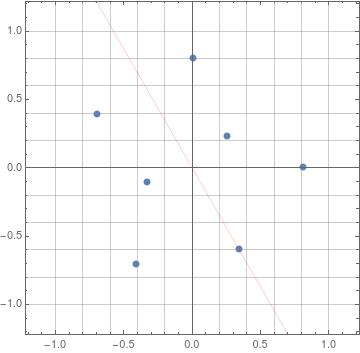}{}{$
\begin{array}{lr}
(-0.69631383, & 0.3878348168)\\
(-0.4083238701, & -0.7074648947)\\
(0.3463253975, & -0.5997643849)\\
(0.2547995791, & 0.2287497666)\\
(0.01213209623, & 0.7969408167)\\
(0.8168276037, & 0)\\
(-0.3255579018, & -0.1063779431)
\end{array}
\mbox{}\\[1ex]
\begin{array}{lcl}
J & = & 0.326504\\
P & = & 42.32581356
\end{array}
$}
\par\medskip
\noindent
\config{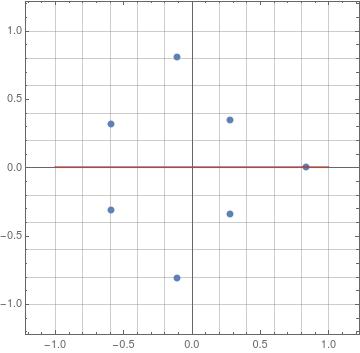}{}{$
\begin{array}{lr}
(-0.5939580764, & 0.3153511693)\\
(-0.1055438498, & -0.8071466615)\\
(0.280704098, & -0.3443258497)\\
(0.280704098, & 0.3443258497)\\
(-0.1055438498, & 0.8071466615)\\
(0.8375956565, & 0)\\
(-0.5939580764, & -0.3153511693)
\end{array}
\mbox{}\\[1ex]
\begin{array}{lcl}
J & = & 0.327518\\
P & = & 42.45998988
\end{array}
$}
\par\medskip
\noindent
\config{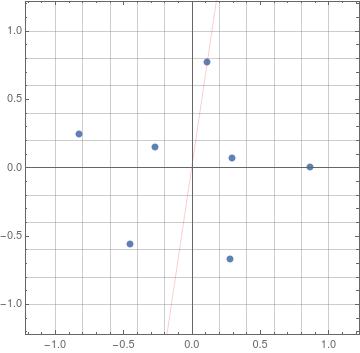}{}{$
\begin{array}{lr}
(-0.8280586249, & 0.2436948559)\\
(0.2762636338, & -0.6696172516)\\
(0.2976897702, & 0.06892828403)\\
(0.1111253801, & 0.7712054397)\\
(-0.2661193233, & 0.1501692824)\\
(0.8631733714, & 0)\\
(-0.4540742073, & -0.5643806103)
\end{array}
\mbox{}\\[1ex]
\begin{array}{lcl}
J & = & 0.3286127349\\
P & = & 42.60195114
\end{array}
$}
\par\medskip
\noindent
\config{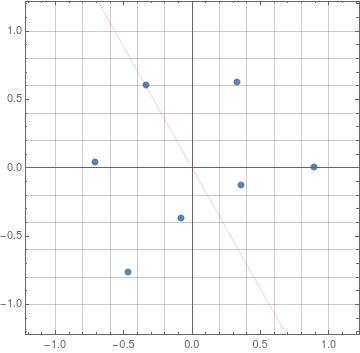}{}{$
\begin{array}{lr}
(-0.7078524538, & 0.0404633721)\\
(-0.4664599419, & -0.7631136127)\\
(0.3604500559, & -0.1275657688)\\
(0.3346502375, & 0.6250612765)\\
(-0.07914736545, & -0.3740761456)\\
(0.894386529, & 0)\\
(-0.3360270613, & 0.5992308785)
\end{array}
\mbox{}\\[1ex]
\begin{array}{lcl}
J & = & 0.3339838173\\
P & = & 43.29826801
\end{array}
$}
\par\medskip
\noindent
\config{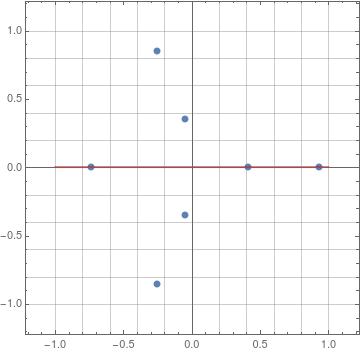}{}{$
\begin{array}{lr}
(-0.7346039147, & -1.044347371\times 10^{-7})\\
(-0.2560042753, & -0.8539984503)\\
(0.4114606324, & -1.134599719\times 10^{-7})\\
(-0.04883637883, & -0.350789775)\\
(-0.04883640953, & 0.3507896809)\\
(0.9328241797, & 0)\\
(-0.2560044998, & 0.8539980966)
\end{array}
\mbox{}\\[1ex]
\begin{array}{lcl}
J & = & 0.3414535429\\
P & = & 44.26665678
\end{array}
$}
\par\medskip
\noindent
\config{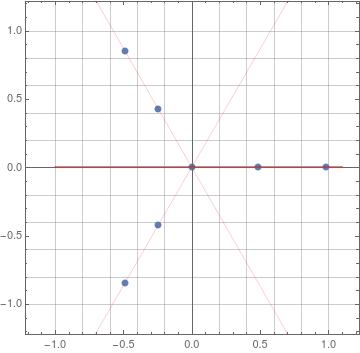}{}{$
\begin{array}{lr}
(-0.4907100136, & 0.8499346753  )\\
(-0.4907100136, & -0.8499346753 )\\
(0.4879772078, & -8.542713523\times 10^{-15})\\
(-2.472525901\times 10^{-14} & -3.692451758\times 10^{-14})\\
(-0.2439886039, & -0.4226006584 )\\
(0.9814200271, & 0 )\\
(-0.2439886039, & 0.4226006584  )
\end{array}
\mbox{}\\[1ex]
\begin{array}{lcl}
J & = & 0.3694161239\\
P & = & 47.89177651
\end{array}
$}
\par\medskip
\noindent
\config{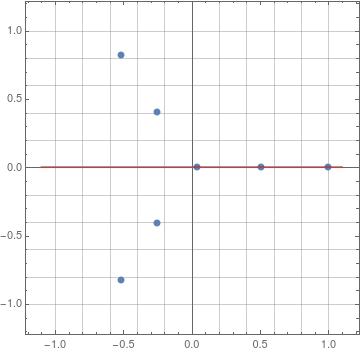}{}{$
\begin{array}{lr}
(-0.5196489643, & 0.8230510161)\\
(-0.5196489643, & -0.8230510161)\\
(0.5088972465, & -5.66181332\times 10^{-15})\\
(0.0366546652, & -1.401575697\times 10^{-14})\\
(-0.2515755625, & -0.405120033)\\
(0.996897142, & 0)\\
(-0.2515755625, & 0.405120033)
\end{array}
\mbox{}\\[1ex]
\begin{array}{lcl}
J & = & 0.3694044172\\
P & = & 47.89025883
\end{array}
$}
\par\medskip
\noindent
\config{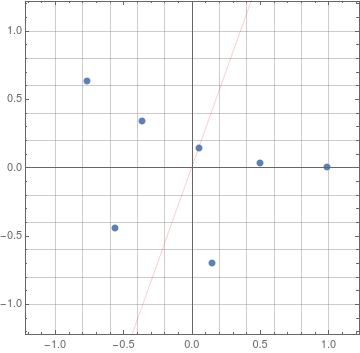}{}{$
\begin{array}{lr}
(-0.7639543131, & 0.6312660195)\\
(0.151096218, & -0.6987810731)\\
(0.496401643, & 0.03141515292)\\
(0.049678093, & 0.1381095023)\\
(-0.3626530541, & 0.3404177807)\\
(0.9910211804, & 0)\\
(-0.5615897672, & -0.4424273823)
\end{array}
\mbox{}\\[1ex]
\begin{array}{lcl}
J & = & 0.3539839884\\
P & = & 45.89112647
\end{array}
$}
\par\medskip
\noindent
\config{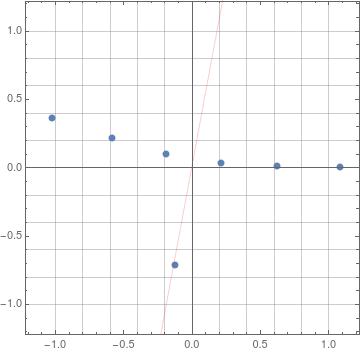}{}{$
\begin{array}{lr}
(-1.023269575, & 0.3609518293)\\
(-0.1222116086, & -0.7138431876 )\\
(0.6226884363, & 0.008260404371 )\\
(0.2110036506, & 0.03062692087 )\\
(-0.1887985548, & 0.09907398529 )\\
(1.085065365, & 0 )\\
(-0.5844777135, & 0.2149300478 )
\end{array}
\mbox{}\\[1ex]
\begin{array}{lcl}
J & = & 0.3914411515\\
P & = & 50.74714105
\end{array}
$}

\subsection{Asymmetrical CCs}\label{sec:asymm-cc}
\subsubsection{Eight bodies}
\config{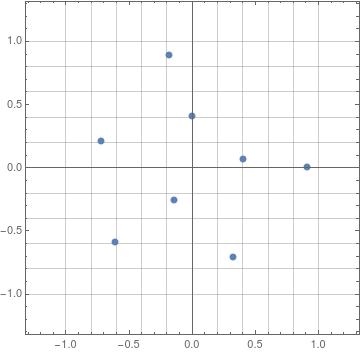}{}{$
\begin{array}{lr}
(-0.7166215124, & 0.2064500901 )\\
(0.3245586105, & -0.7128796496)\\
(0.4049275132, & 0.063605646)\\
(0.001837695963, & 0.402546434)\\
(-0.1438928121, & -0.2575373186)\\
(-0.181586303, & 0.8894017544)\\
(0.9158898541, & 0)\\
(-0.6051130464, & -0.5915869564 )
\end{array}
\mbox{}\\[1ex]
\begin{array}{lcl}
J & = & 0.349028\\
P & = & 63.18080222
\end{array}
$}
\par\medskip
\noindent
\config{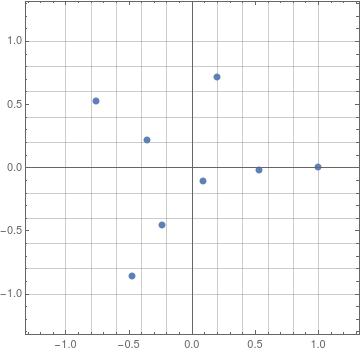}{}{$
\begin{array}{lr}
(-0.758499825, & 0.5229211679)\\
(-0.4757362357, & -0.8641489841)\\
(0.5297384082, & -0.02473539829)\\
(0.2038227669, & 0.7147485589)\\
(0.09156513918, & -0.1080754295)\\
(-0.2344231355, & -0.4560136147)\\
(0.9995383601, & 0 )\\
(-0.3560054781, & 0.2153036998 )
\end{array}
\mbox{}\\[1ex]
\begin{array}{lcl}
J & = & 0.3683220063\\
P & = & 66.67340501
\end{array}
$}

\subsubsection{Nine bodies}
\config{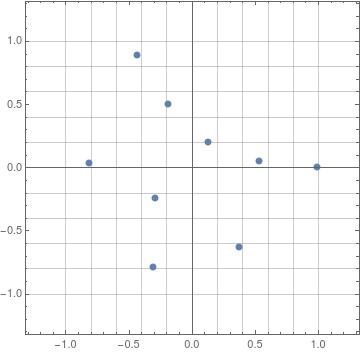}{}{$
\begin{array}{lr}
-0.8098124055, & 0.03092225943)\\
-0.3079530509, & -0.7934317414)\\
0.5336942253, & 0.05249596102)\\
0.3749583879, & -0.6313253357)\\
0.1303879569, & 0.2002134463)\\
-0.1885012036, & 0.4988431245 )\\
-0.2878164038, & -0.2447179668 )\\
0.9919453935, & 0)\\
-0.4369028999, & 0.8870002527 )
\end{array}
\mbox{}\\[1ex]
\begin{array}{lcl}
J & = & 0.3718173376\\
P & = & 90.35161303
\end{array}
$}
\par\medskip
\noindent
\config{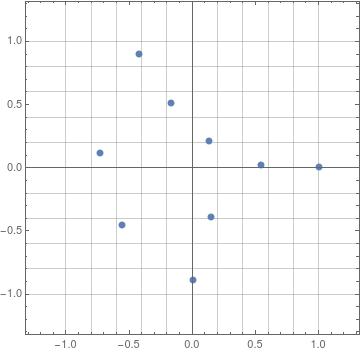}{}{$
\begin{array}{lr}
(-0.7275085134, & 0.1152833946)\\
(0.01159366005, & -0.8949186222)\\
(0.5488693567, & 0.0194418351)\\
(0.1519951733, & -0.3967802948)\\
(0.1389256054, & 0.2109317686)\\
(-0.1649134702, & 0.5071400693)\\
(-0.4146605713, & 0.8983002319)\\
(1.004210874, & 0)\\
(-0.5485121146, & -0.4593983825 )
\end{array}
\mbox{}\\[1ex]
\begin{array}{lcl}
J & = & 0.374156044\\
P & = & 90.9199187
\end{array}
$}
\par\medskip
\noindent
\config{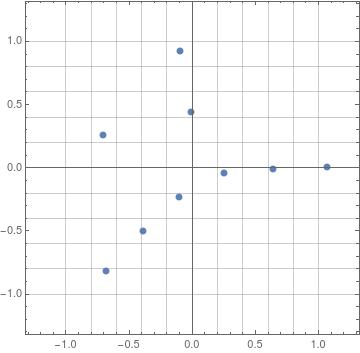}{}{$
\begin{array}{lr}
(-0.7030520826, & 0.2581361323)\\
(-0.6798530659, & -0.8186027881)\\
(0.6455196701, & -0.01022986752)\\
(0.2525581245, & -0.0456682048)\\
(-0.005489709791, & 0.4378314751)\\
(-0.09458813233, & 0.9180310512)\\
(-0.09787656765, & -0.2387132415 )\\
(1.072063843, & 0 )\\
(-0.3892820798, & -0.5007845568 )
\end{array}
\mbox{}\\[1ex]
\begin{array}{lcl}
J & = & 0.3940726241\\
P & = & 95.75964764
\end{array}
$}

\subsubsection{Ten bodies}
\config{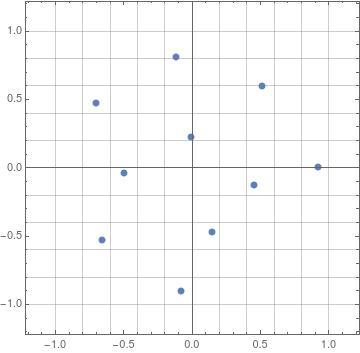}{}{$
\begin{array}{lr}
(-0.7007408215, & 0.4693295286 )\\
(-0.07639930863, & -0.9073170049 )\\
(0.5132734942, & 0.5961143132 )\\
(0.4588100327, & -0.1317284788 )\\
(0.1456294856, & -0.4767952142 )\\
(-0.002441745604, & 0.2214431849 )\\
(-0.1155249827, & 0.8058001754 )\\
(-0.4929493421, & -0.04332608545 )\\
(0.9269083278, & 0 )\\
(-0.6565651399, & -0.5335204189 )\\
\end{array}
\mbox{}\\[1ex]
\begin{array}{lcl}
J & = & 0.3714169116\\
P & = & 117.4523402
\end{array}
$}
\par\medskip
\noindent
\config{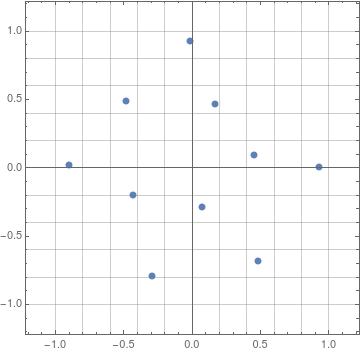}{}{$
\begin{array}{lr}
(-0.8963214942, & 0.01882076642 )\\
(-0.2944424421, & -0.7961815748 )\\
(0.4843083474, & -0.6834013227 )\\
(0.4538131131, & 0.08780569651 )\\
(0.1725919232, & 0.460572645 )\\
(0.07492379596, & -0.289147328 )\\
(-0.0101214633, & 0.9213652726 )\\
(-0.4311240011, & -0.2030210622 )\\
(0.9304861223, & 0 )\\
(-0.4841139013, & 0.483186907 )\\
\end{array}
\mbox{}\\[1ex]
\begin{array}{lcl}
J & = & 0.3728671543\\
P & = & 117.9109472
\end{array}
$}
\par\medskip
\noindent
\config{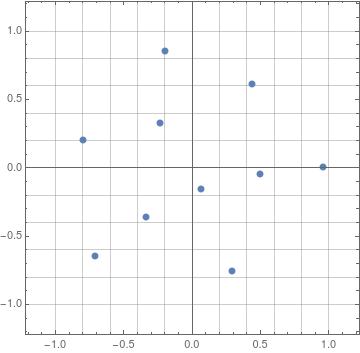}{}{$
\begin{array}{lr}
(-0.7967875604, & 0.1958277257 )\\
(0.2919426026, & -0.7600314103 )\\
(0.5006139893, & -0.0465575355 )\\
(0.4399251191, & 0.6127399142 )\\
(0.06874098874, & -0.1558927224 )\\
(-0.1936088183, & 0.8476831709 )\\
(-0.2311572039, & 0.3226693268 )\\
(-0.3368118774, & -0.367337145 )\\
(0.9623051252, & 0 )\\
(-0.705162365, & -0.6491013244 )\\
\end{array}
\mbox{}\\[1ex]
\begin{array}{lcl}
J & = & 0.3742731763\\
P & = & 118.3555704
\end{array}
$}
\par\medskip
\noindent
\config{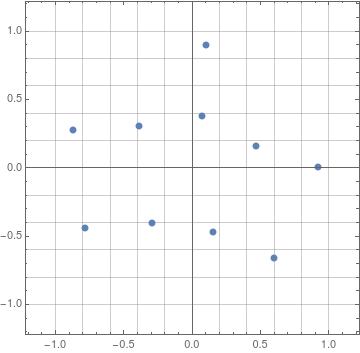}{}{$
\begin{array}{lr}
(-0.8675091418, & 0.272352247 )\\
(0.6017849903, & -0.6658290959 )\\
(0.4701929967, & 0.1559215689 )\\
(0.1519788533, & -0.476594594 )\\
(0.1008116833, & 0.894979045 )\\
(0.0754051335, & 0.3753489653 )\\
(-0.292654655, & -0.4085342723 )\\
(-0.3838258452, & 0.29952652 )\\
(0.9222584288, & 0 )\\
(-0.7784424448, & -0.4471703844 )\\
\end{array}
\mbox{}\\[1ex]
\begin{array}{lcl}
J & = & 0.3784068394\\
P & = & 119.6627495
\end{array}
$}
\par\medskip
\noindent
\config{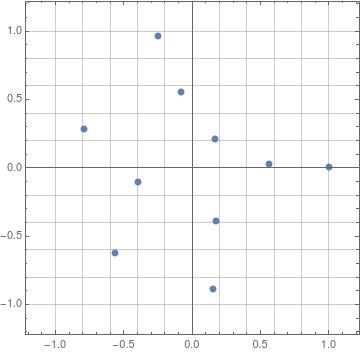}{}{$
\begin{array}{lr}
(-0.7888545233, & 0.2769354218 )\\
(0.1573661438, & -0.8879803722 )\\
(0.566350671, & 0.02117690755 )\\
(0.1809819542, & -0.3943353946 )\\
(0.1673309159, & 0.2071064401 )\\
(-0.08156809045, & 0.5470717464 )\\
(-0.2490045724, & 0.9598108791 )\\
(-0.3968333769, & -0.1057160826 )\\
(1.005616634, & 0 )\\
(-0.5613857558, & -0.6240695457 )
\end{array}
\mbox{}\\[1ex]
\begin{array}{lcl}
J & = & 0.3821740131\\
P & = & 120.8540344
\end{array}
$}
\par\medskip
\noindent
\config{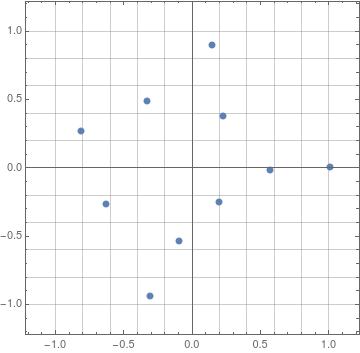}{}{$
\begin{array}{lr}
(-0.8111684557, & 0.2656551367 )\\
(-0.3059600863, & -0.9399118629 )\\
(0.5694847219, & -0.02355906915 )\\
(0.2319114215, & 0.3782046907 )\\
(0.2014838749, & -0.2518381211 )\\
(0.1490178487, & 0.8945560864 )\\
(-0.09081819785, & -0.5373390538 )\\
(-0.3293580397, & 0.4819496801 )\\
(1.01098039, & 0 )\\
(-0.6255734779, & -0.2677174871 )
\end{array}
\mbox{}\\[1ex]
\begin{array}{lcl}
J & = & 0.3832919194\\
P & = & 121.2075474
\end{array}
$}
\par\medskip
\noindent
\config{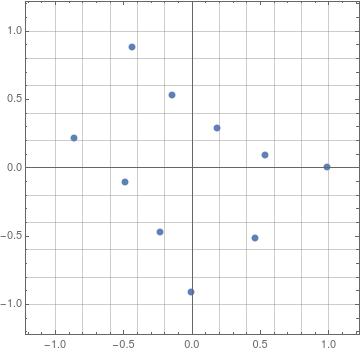}{}{$
\begin{array}{lr}
(-0.8644611425, & 0.2149252766 )\\
(-0.00513249289, & -0.9138239369 )\\
(0.538361395, & 0.09282895599 )\\
(0.46480923, & -0.5151381858 )\\
(0.1823759207, & 0.2881651807 )\\
(-0.1446907068, & 0.5315697209 )\\
(-0.2305530662, & -0.4706217805 )\\
(-0.4362270962, & 0.8774029291 )\\
(0.9865289623, & 0 )\\
(-0.4910113567, & -0.1053084921 )
\end{array}
\mbox{}\\[1ex]
\begin{array}{lcl}
J & = & 0.3845845407\\
P & = & 121.6163101
\end{array}
$}
\par\medskip
\noindent
\config{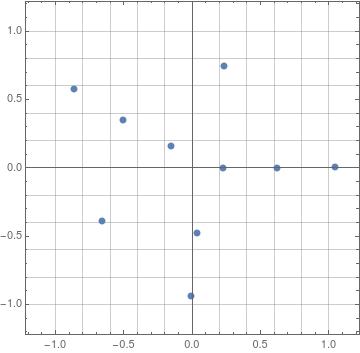}{}{$
\begin{array}{lr}
(-0.8620127236, & 0.5749149667 )\\
(-0.003887282198, & -0.9423644458 )\\
(0.6244444875, & -0.003510447333 )\\
(0.2389418223, & 0.7433087445 )\\
(0.2253692474, & -0.004873923735 )\\
(0.04037831843, & -0.4785847827 )\\
(-0.1523931658, & 0.1569749357 )\\
(-0.5005974321, & 0.3480983721 )\\
(1.045810599, & 0 )\\
(-0.656053871, & -0.3939634194 )
\end{array}
\mbox{}\\[1ex]
\begin{array}{lcl}
J & = & 0.3904041955\\
P & = & 123.4566466
\end{array}
$}
\par\medskip
\noindent
\config{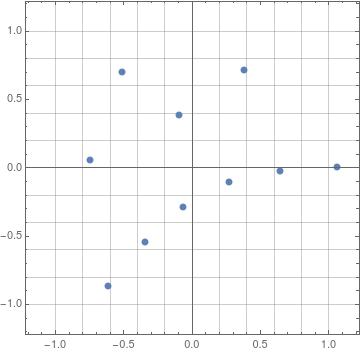}{}{$
\begin{array}{lr}
(-0.7437201474, & 0.05424247592 )\\
(-0.6097575083, & -0.8688275537 )\\
(0.6456942806, & -0.03032211044 )\\
(0.3841057237, & 0.7124794069 )\\
(0.2719022074, & -0.1083651492 )\\
(-0.06052210988, & -0.2908103139 )\\
(-0.09397879995, & 0.3792041934 )\\
(-0.3436577475, & -0.5452330647 )\\
(1.062255195, & 0 )\\
(-0.5123211036, & 0.6976321062 )
\end{array}
\mbox{}\\[1ex]
\begin{array}{lcl}
J & = & 0.3940864744\\
P & = & 124.6210854
\end{array}
$}
\par\medskip
\noindent
\config{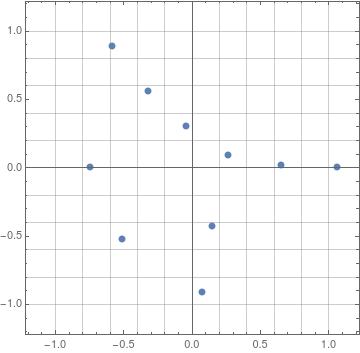}{}{$
\begin{array}{lr}
(-0.7443869099, & 0.004965386648 )\\
(0.07184450244, & -0.910877299 )\\
(0.6495486297, & 0.01650001945 )\\
(0.2638544878, & 0.08972086605 )\\
(0.1507954307, & -0.4264727703 )\\
(-0.0449474494, & 0.3011792823 )\\
(-0.3224455973, & 0.5582706437 )\\
(-0.5083858126, & -0.5225251633 )\\
(1.065046288, & 0 )\\
(-0.5809235696, & 0.8892390345 )
\end{array}
\mbox{}\\[1ex]
\begin{array}{lcl}
J & = & 0.3963617068\\
P & = & 125.3405771
\end{array}
$}
\par\medskip
\noindent
\config{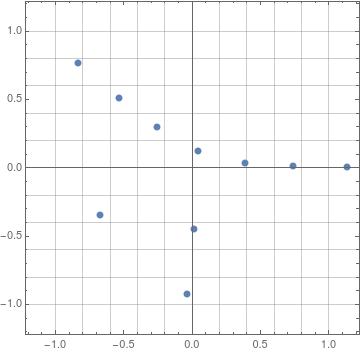}{}{$
\begin{array}{lr}
(-0.8333907632, & 0.7644395423 )\\
(-0.03372694114, & -0.9240673101 )\\
(0.7379960306, & 0.009114360497 )\\
(0.386475537, & 0.03160247231 )\\
(0.04528434434, & 0.1213495026 )\\
(0.018708438, & -0.4551657158 )\\
(-0.2519844233, & 0.2944278672 )\\
(-0.5319428417, & 0.5069057403 )\\
(1.132810121, & 0 )\\
(-0.6702295016, & -0.3486064594 )
\end{array}
\mbox{}\\[1ex]
\begin{array}{lcl}
J & = & 0.4187765849\\
P & = & 132.4287839
\end{array}
$}

\end{document}